\documentclass[reqno]{amsart}

\usepackage{enumerate}
\usepackage[colorlinks=true]{hyperref}
\hypersetup{citecolor=blue}

\usepackage{bbm}

\numberwithin{equation}{section}

\newtheorem{theorem}{Theorem}[section]
\newtheorem{lemma}[theorem]{Lemma}
\newtheorem{proposition}[theorem]{Proposition}
\theoremstyle{definition}
\newtheorem{remark}[theorem]{Remark}

\newcommand\R{{\mathbb R}}

\newcommand\N{{\mathbb N}}

\newcommand\Comp{{\mathrm{c}}}

\newcommand\Cz{{C_0(\R^N )}}

\newcommand\Czo{{C_0(\Omega )}}

\newcommand\wN{{\widetilde {N}}}

\newcommand\Lt{{\mathcal L}_2}
\newcommand\Ltu{{\mathcal L}_1}

\newcommand\Loc{{\mathrm{loc}}}

\newcommand\goto{\mathop{\longrightarrow}}

\newcommand\MScN[1]{\href{http://www.ams.org/mathscinet-getitem?mr=#1}{\nolinkurl{(#1)}}}
\newcommand\DOI[1]{\href{http://dx.doi.org/#1}{(doi: \nolinkurl{#1})}}
\newcommand\LINK[1]{\href{#1}{(link: \nolinkurl{#1})}}

\newcommand\DI{u_0 }

\newcommand\NUU{{\eta}}
\newcommand\MUU{\mu }
\newcommand\MUUz{1 }
\newcommand\Phimu{W}

\begin{document}

\title{Sign-changing self-similar solutions of the nonlinear heat equation with positive initial value}

\def\shorttitle{Sign-changing self-similar solutions}

\author[T. Cazenave]{Thierry Cazenave$^1$}
\email{\href{mailto:thierry.cazenave@sorbonne-universite.fr}{thierry.cazenave@sorbonne-universite.fr}}

\author[F. Dickstein]{Fl\'avio Dickstein$^{1,2}$}
\email{\href{mailto:flavio@labma.ufrj.br}{flavio@labma.ufrj.br}}

\author[I. Naumkin]{Ivan Naumkin$^3$}
\email{\href{mailto:ivan.naumkin@iimas.unam.mx}{ivan.naumkin@iimas.unam.mx}}

\author[F. B.~Weissler]{Fred B.~Weissler$^4$}
\email{\href{mailto:weissler@math.univ-paris13.fr}{weissler@math.univ-paris13.fr}}

\address{$^1$Sorbonne Universit\'e \& CNRS, Laboratoire Jacques-Louis Lions,
B.C. 187, 4 place Jussieu, 75252 Paris Cedex 05, France}

\address{$^2$Instituto de Matem\'atica, Universidade Federal do Rio de Janeiro, Caixa Postal 68530, 21944--970 Rio de Janeiro, R.J., Brazil}

\address{$^3$Departamento de F\'\i sica Matem\'atica,
Instituto de Investigaciones en Matem\'aticas Aplicadas y en Sistemas, 
Universidad Nacional Aut\'onoma de M\'exico, 
Apartado Postal 20-126, M\'exico DF 01000, M\'exico}

\address{$^4$Universit\'e Sorbonne Paris Nord, LAGA CNRS UMR 7539, 99 Avenue J.-B. Cl\'e\-ment, F-93430 Villetaneuse, France}



\thanks{Research supported by the ``Brazilian-French Network in Mathematics"}
\thanks{Fl\'avio Dickstein was partially supported by CNPq (Brasil).}
\thanks{Ivan Naumkin is a Fellow of Sistema Nacional de Investigadores, and he was supported in part by project PAPIIT IA101820.
This work was prepared while he was visiting the Laboratoire J.A. Dieudonn\'e of the Universit\'e de Nice Sophia-Antipolis. 
He thanks the project ERC-2014-CdG 646.650 SingWave for its financial support, and the Laboratoire J.A. Dieudonn\'e for its kind hospitality. }

\begin{abstract}

We consider the nonlinear heat equation $u_t - \Delta u =  |u|^\alpha u$ on ${\mathbb R}^N$, where $\alpha >0$ and $N\ge 1$.  We prove that in the range $0 < \alpha <\frac {4} {N-2}$, for every $\mu  >0$, there exist infinitely many sign-changing, self-similar solutions to the Cauchy problem with initial value $u_0 (x)= \mu   |x|^{-\frac {2} {\alpha }}$. The construction is based on the analysis of the related inverted profile equation.
In particular, we construct (sign-changing) self-similar solutions for positive initial values for which it is known that there does not exist any local, nonnegative solution.
\end{abstract}

\maketitle

\section{Introduction} \label{Intro} 

In this paper we prove the existence of radially symmetric self-similar solutions
of the nonlinear heat equation on $\R^N $
\begin{equation} \label{NLHE}
u_t = \Delta u +  |u|^\alpha u
\end{equation}
with initial value
$ \MUU  |x|^{-\frac {2} {\alpha }}$, where 
$0 < \alpha <\frac {4} {N-2}$ ($0< \alpha < \infty $ if $N=1,2$)
and $ \MUU  \in \R$, $ \MUU  \neq 0$.  These solutions are classical solutions, in $C((0,\infty); \Cz )$,
and the initial value is realized in the sense of $L^1(\R^N) + L^\infty(\R^N)$
if $\frac{2}{N} < \alpha < \frac {4} {N-2}$, and in the sense of
$\mathcal{D}'(\R^N\setminus\{0\})$ if $0 < \alpha \le \frac {2} {N}$.
In all these cases, as $t \to 0$, the solution approaches the initial value uniformly
on the exterior of any ball around $0$.  In fact, for every $\alpha$
in the range $0 < \alpha < \frac {4} {N-2}$, and every $ \MUU  \in \R$, $ \MUU \not = 0$,
we prove the existence of {\it infinitely} many self-similar solutions
of \eqref{NLHE} with initial value $ \MUU  |x|^{-\frac {2} {\alpha }}$.

This result is significant for several reasons.  First, while it is already known that for every $\alpha > 0$,  the equation \eqref{NLHE} has at least one self-similar solution
with initial value $ \MUU  |x|^{-\frac {2} {\alpha }}$ if $ \MUU  \in \R$ is sufficiently small, including $\MUU =0$, we establish this fact for {\it all} $ \MUU  \in \R$.  Our restriction to the range $0 < \alpha < \frac {4} {N-2}$ is sharp; for $\alpha \ge \frac {4} {N-2}$,
 \eqref{NLHE} does not admit a radially symmetric self-similar solution with initial value
$ \MUU  |x|^{-\frac {2} {\alpha }}$ if $| \MUU |$ is large. See Remark~\ref{zRem2}~\eqref{eSE4}, \eqref{eSE5} and~\eqref{eSE6}, as well as~\cite[Theorem~6.1]{CW}.

Second, in the range $0 < \alpha < \frac {4} {N-2}$, it is known (Remark~\ref{zRem2}~\eqref{eSE4} and~\eqref{eSE5}) that  there exists an arbitrarily large (finite) number of self-similar
solutions with initial value $ \MUU  |x|^{-\frac {2} {\alpha }}$ if $| \MUU |$ is sufficiently
small, and only if $ \MUU  = 0$ is the existence of infinitely many such solutions known (Remark~\ref{zRem2}~\eqref{eSE2}, \eqref{eSE3} and~\eqref{eSE3:1}).
Here we establish the existence of infinitely many self-similar solutions
with initial value $ \MUU  |x|^{-\frac {2} {\alpha }}$, for all $ \MUU  \in \R$. 

Third, for all $ \MUU \in \R$ 
 in the case $0 < \alpha \le \frac {2} {N}$,
and for all $ \MUU  \in \R$ with $| \MUU |$ sufficiently large
 in the case $\frac{2}{N} < \alpha < \frac {4} {N-2}$,
 the self-similar solutions we construct are all sign-changing.
  In other words, the positive initial value 
$ \MUU  |x|^{-\frac {2} {\alpha }}$, with $ \MUU  > 0$,
gives rise to solutions which assume both signs for every $t > 0$.
This is particularly striking since
equation \eqref{NLHE} does not have {\it any nonnegative} local in time solution
with initial value $ \MUU  |x|^{-\frac {2} {\alpha }}$
 with $ \MUU  > 0$ sufficiently large in the case $\alpha >\frac {2} {N} $. 
The same is true for all $ \MUU  > 0$ in the case $0 < \alpha \le \frac {2} {N}$ since $| \cdot |^{-\frac {2} {\alpha }} \not \in L^1_\Loc (\R^N ) $. (See Proposition~\ref{eNEX1} in Appendix~\ref{sNEX} for details on these nonexistence properties.)
Thus, we construct sign-changing solutions of~\eqref{NLHE}
with a positive initial value for which there is no
local in time nonnegative solution.

This last point merits further explanation. Consider the initial value problem
\begin{equation} \label{Inlh} 
\begin{cases} 
u_t = \Delta u +  |u|^\alpha u \\
u ( 0, \cdot ) = \DI (\cdot )
\end{cases} 
\end{equation} 
where $u= u(t,x)$, $t\ge 0$, $x \in \Omega $ where $\Omega $ is a domain in $\R^N $ (possibly $\Omega =\R^N $), and $\alpha >0$.
In the case where $\Omega \not = \R^N $, we impose Dirichlet boundary conditions. 
It is well-known that this problem is locally well-posed in $\Czo $,  and also in $L^p (\Omega ) $ for $p\ge 1$, $p>\frac {N\alpha } {2}$. See \cite{Weissler1,Weissler2,BrezisC}.
Moreover, if $\DI \ge 0$ in $\Omega $, then the resulting solution satisfies $u\ge 0$ in $\Omega $. This is a consequence of the iterative method used to construct solutions, based on Duhamel's formula.

On the other hand, this problem is not well-posed in $L^p (\Omega ) $ if $\alpha >\frac {2} {N}$ and $1\le p < \frac {N\alpha } {2}$. For example, regular initial values can yield multiple solutions which are continuous into these spaces.  See \cite{HarauxW, Baras, QuittnerS}.  Also, it was observed in \cite{Weissler2, Weissler4}  (see also \cite{QuittnerS,LaisterRSVL, FujishimaI})  that there are nonnegative  $\DI \in L^p (\Omega )$ for which there is no local-in-time nonnegative solution in the weakest possible sense. 
This last fact has often been considered as a proof of the non-existence of solutions of~\eqref{Inlh} associated to those initial values.
However, the possibility remains that positive initial data can give rise to local solutions which assume both positive and negative values. This of course would seemingly violate the maximum principle. Such a possibility is not completely unknown. Indeed, in~\cite{HarauxW} both positive and negative solutions were constructed with initial value $0$. Even though these solutions are regular for $t>0$, they are too singular as $(t,x)\to (0,0)$ for any maximum principle to apply.

The present paper gives the first result, to our knowledge, of the existence of a sign-changing solution of \eqref{Inlh} with an initial
value $\DI \in L^1_\Loc (\R^N ) $, $\DI \ge 0$ for which {\it no} local in time nonnegative solution
exists.  While the initial value $\DI(x) =  \MUU  |x|^{-\frac {2} {\alpha }}$
is not in any $L^p$ space, we show in a subsequent article~\cite{CDNW2} that the solutions
constructed here can be perturbed to give solutions to \eqref{Inlh}
with nonnegative initial value in  $L^p$, $1\le p < \frac {N\alpha } {2}$, with the same property.

In order to state our results more precisely, we recall some known facts about self-similar solutions of \eqref{NLHE}.
A self-similar solution of~\eqref{NLHE} is a solution of the form
\begin{equation} \label{fpr1} 
u(t, x) = t^{-\frac {1} {\alpha }} f \Bigl( \frac {x} {\sqrt t} \Bigr),
\end{equation}
where $f  : \R^N \to \R$ is  the profile of the
self-similar solution $u$ given by \eqref{fpr1}.  In order for $u$ given by 
\eqref{fpr1} to be a classical solution of  \eqref{NLHE} for $t > 0$,
the profile $f$ must be of class $C^2$ and satisfy the elliptic equation
\begin{equation} \label{fpr2} 
\Delta f + \frac {1} {2} x\cdot \nabla f + \frac {1} {\alpha } f +  |f|^\alpha f  = 0.
\end{equation} 
In our investigations, we look only for radially symmetric profiles (except if $N = 1$),
and so we write, by abuse of notation, $f(r) = f(x)$ where $r= |x|$, so that
$f : [0, \infty) \to \R$ is of class $C^2$, and satisfies the following
initial value ODE problem,
\begin{gather}
\displaystyle f''(r) + \Big(\frac{N-1}{r} +  \frac{r}{2}\Big)f'(r) + \frac{1}{\alpha}f(r) +  |f(r)|^{\alpha}f(r) = 0 \label{fpr3}\\
f(0) = a, \quad f'(0) = 0 \label{fpr3:1}
\end{gather}
for some $a\in \R$.
Of course, if $N = 1$, it is not necessary that $f'(0) = 0$ in order
to obtain a regular profile on $\R$.  See Theorem~\ref{eMain3} below.

It is known \cite{HarauxW} that the problem~\eqref{fpr3}-\eqref{fpr3:1} is well-posed.
More precisely, given $a \in \R$, there exists a unique solution
$f_a \in C^2([0,\infty);\R)$. Furthermore the limit
\begin{equation} \label{prlm}
L(a) = \lim_{r\to\infty}r^{\frac{2}{\alpha}}f_a(r) \in \R
\end{equation}
exists, and is a locally Lipschitz function of $a \in \R$.
If $L(a) = 0$, then $f _a $ decays exponentially.  See~\cite{PTW} for more
precise information on the asymptotic behavior of solutions
to~\eqref{fpr3}-\eqref{fpr3:1}.  Moreover, if $a \neq 0$, then $f_a$ has at most finitely
many zeros and we set
\begin{equation}
\label{nz}
N(a)  = \text{the number of zeros of the function $f_a$} .
\end{equation}
Finally, given a radially symmetric self-similar solution $u$
of the form \eqref {fpr1}, with $f = f_a$, its initial value can be easily determined.
Indeed, if $x \neq 0 $,
\begin{equation} \label{ssiv}
\lim_{t \to 0}u(t,x) = \lim_{t \to 0}t^{-\frac {1} {\alpha }} f_a \Bigl( \frac {|x|} {\sqrt t} \Bigr)
= |x|^{-\frac {2} {\alpha }}\lim_{r\to\infty}r^{\frac{2}{\alpha}}f_a(r)
= L(a)|x|^{-\frac {2} {\alpha }}.
\end{equation}
It follows that every radially symmetric regular self-similar solution of \eqref{NLHE}
with initial value $\DI= \MUU |x|^{-\frac {2} {\alpha }}$ is given, {\it via} the formula
\eqref{fpr1}, by a profile $f = f_a$ which is a solution
of~\eqref{fpr3}-\eqref{fpr3:1} such that $L(a) =  \MUU $.
Moreover, \eqref{fpr1} and~\eqref{prlm} imply that 
\begin{equation} \label{fEP1} 
 |u(t, x)| \le C (t+  |x|^2 )^{- \frac {1} {\alpha }} . 
\end{equation} 
One deduces easily from~\eqref{ssiv} and~\eqref{fEP1}  that 
\begin{equation} \label{fEP2} 
u(t) \goto _{ t\downarrow 0 } \DI  \text{ in } L^q  ( \{  |x|>\varepsilon  \} )  \text{ for all } \varepsilon >0  \text{ and } q \ge 1,  q >\frac {N\alpha } {2} ,
\end{equation} 
and that if  $\alpha > \frac {2} {N}$, then 
\begin{equation} \label{fEP3} 
u(t) \goto _{ t\downarrow 0 } \DI  \text{ in } L^p (\R^N ) +L^q  (\R^N )   \text{ for all }  1\le p < \frac {N\alpha } {2} <q .
\end{equation}

We are now able to state our main result. It concerns the case $N\ge 2$. In dimension $1$, we have a similar result (see Theorem~\ref{eMain3} below), whose proof is somewhat different.

\begin{theorem}  \label{eMain1}
Assume
\begin{equation} \label{fNZ1} 
N\ge 2, \quad \alpha >0,\quad  \text{and}\quad (N-2) \alpha < 4 ,
\end{equation}
 and  let $\MUU \in \R$, $ \MUU  \not =  0$. 
There exists $m_0\ge 0$ such that for all $m \ge m_0$ there exist at least two different, radially symmetric regular self-similar solutions $u$ of \eqref{NLHE} with initial value $\DI=  \MUU  |x|^{-\frac {2} {\alpha }}$ in the sense~\eqref{fEP2}, and also~\eqref{fEP3} if $\alpha >\frac {2} {N}$, and whose profiles have exactly $m$ zeros. These solutions are such that $u\in C^1((0,\infty ), L^r (\R^N ) )$ and $\Delta u,  |u|^\alpha u\in C((0,\infty ), L^r (\R^N ) )$ for all $r\ge 1$, $r> \frac {N\alpha } {2}$. 

Furthermore, if  $\alpha >\frac {2} {N}$,  the solutions satisfy the integral equation 
\begin{equation} \label{fNZ2} 
u(t) = e^{t\Delta } \DI + \int _0^t e^{ (t-s) \Delta }  |u(s)|^\alpha u(s)\, ds
\end{equation} 
where each term is in $C((0,\infty ), L^r (\R^N ) )$ for all $r> \frac {N\alpha } {2} > 1$. Moreover, the map $t \mapsto u (t) - e^{t \Delta } \DI $ is in $ C([0,\infty ), L^r (\R^N ) )$ for all $r \ge 1$ such that $\frac {N\alpha } {2(\alpha +1) } < r < \frac {N\alpha } {2}$.
\end{theorem} 

To prove this theorem, we need to show, among other things, that the function $L$
takes on every value $ \MUU \in\R$ infinitely often.  More precisely, given $ \MUU  \in \R$, $ \MUU  \neq 0$, we need to show that for all sufficiently large integers $m$, there exist at least two values of $a$ such that $L(a) =  \MUU $ and $N(a) = m$. Note that since $L(-a) = -L(a)$, it suffices to consider $ \MUU  > 0$. To put these assertions in the appropriate historical context, and describe our approach to the proof, we recall some of the known results about  the solutions $f_a$ of~\eqref{fpr3}-\eqref{fpr3:1} and the function $L$ defined in~\eqref{prlm}. We let
\begin{equation}
\label{bta}
\beta = \frac {2} {\alpha }  \Bigl( N -2 - \frac {2} {\alpha } \Bigr).
\end{equation}
If $\beta > 0$, i.e. if $\alpha > \frac{2}{N-2}$, then
\begin{equation} \label{statsingsol}
u(x) = \beta^{\frac{1}{\alpha}}|x|^{-\frac{2}{\alpha}}
\end{equation}
is a singular stationary solution of \eqref{NLHE}. It is also self-similar, with singular profile $f(x) = \beta^{\frac{1}{\alpha}}|x|^{-\frac{2}{\alpha}}$.  We will see later (see Theorem~\ref{eMain2}) that \eqref{NLHE} has
other singular self-similar solutions when $\beta \ge 0$.  Hence the need to specify
that the self-similar solutions in Theorem~\ref{eMain1} are regular.

The detailed study of the profiles $f_a$ is based on the numbers
\begin{equation}
\label{am}
a_m = \inf\{a > 0; N(a) \ge m+1\}
\end{equation}
for $m\ge 0$,  first defined in~\cite{Weissler6}. The following remark recalls some of the important properties of these numbers.

\begin{remark} \label{zRem2} 
\begin{enumerate}[{\rm (i)}] 
\item \label{eSE1}  
If $0 < a_m < \infty$, then $L(a_m) = 0$ and $N(a_m) = m$.
Furthermore, if $0 < a_m < \infty$ for all large $m$, then
$a_m \to \infty$ as $m\to \infty $. See~\cite[Theorem~1]{Weissler6}.

\item \label{eSE2}
If $\frac {2} {N} < \alpha < \frac{4}{N-2}$, then  $0 < a_m < \infty$ for all
$m \ge 0$.  In particular, there are infinitely many radially symmetric self-similar
solutions of \eqref{NLHE} with initial value $0$, including one which is positive. See~\cite{HarauxW} and~\cite[Theorem~1]{Weissler6}. 

\item \label{eSE3}
 If $0 < \alpha \le \frac {2} {N}$, then there exists $m_0 \ge 0$
such that $0 < a_m < \infty$ for all $m \ge m_0$.  In particular, there are infinitely many radially symmetric self-similar solutions of \eqref{NLHE} with initial value $0$. See~\cite[Theorems~1 and~2]{Weissler6} and~\cite[Proposition~2]{Yanagida}. (In this case, $a_0 =0$ by~\cite[Theorem~5~(a)]{HarauxW}.)

\item \label{eSE3:1}
Parts~\eqref{eSE2} and~\eqref{eSE3} above show that Theorem~\ref{eMain1} is already known in the case $\MUU =0$. It suffices to consider the profiles $f_a$ with $a= \pm a_m$.

\item \label{eSE4} If $0 < \alpha  < \frac{2}{N-2}$, then $N(a)$ is a nondecreasing
function of $a > 0$ and the numbers $a_m$ are uniquely determined
by the property that $L(a_m) = 0$ and $N(a_m) = m$. Furthermore, if
$0 < a_m < a < a_{m+1}$, then $N(a) = m+1$ and $L(a) > 0$. See~\cite{Yanagida}.

\item \label{eSE5} If $\frac{2}{N-2} \le \alpha < \frac{4}{N-2}$, then $a_0$ is  uniquely
determined by the property $L(a_0) = 0$ and $N(a_0) = 0$ (i.e. $f_{a_0}(r) > 0$ for all $r > 0$). Moreover,  if $0 < a < a_0$, then
$f_a(r) > 0$ for all $r > 0$ and $L(a) > 0$. If $a > a_0$, then $N(a) \ge 1$.
(See~\cite{DohmenH}.)  In this range of $\alpha $, it is also true that
for every $m \ge 0$, there exists $a_m < a < a_{m+1}$ such that
$N(a) = m+ 1$ and $L(a) > 0$.  This is not explicitly proved anywhere,
as far as we know.  It does follow from the results in \cite{Weissler6} along
with a slight improvement of Proposition 3.7 in that paper which is straightforward to prove.

\item \label{eSE6} If $\alpha \ge \frac{4}{N-2}$ and $a > 0$, then $f_a(r) > 0$ for all
$r > 0$ and $L(a) > 0$ (see \cite{HarauxW}).
Moreover,  $L:\R \to \R$ is bounded. 
Indeed, there is no local in time positive solution of~\eqref{NLHE} with the initial value $\MUU  |x|^{-\frac {2} {\alpha }}$ for $\MUU $ large, see~\cite[Theorem~1]{Weissler4} (and Lemma~\ref{eQSol}). 

\end{enumerate}
\end{remark} 

In light of these properties  and in view of \eqref{ssiv}, to prove Theorem~\ref{eMain1}, it would suffice to  show
(at least in the case $0 < \alpha  < \frac{2}{N-2}$) that the successive maxima of $ | L (a)| $ on the intervals $[ a_m , a _{ m+1 } ]$ tend to infinity.
This is in fact our strategy. See formulas \eqref{eNZ3:1} and \eqref{eNZ3:3} below.
To accomplish this, however, we do not directly study
 $L(a)$ using equation~\eqref{fpr3}-\eqref{fpr3:1}.
Instead, we specify an arbitrary value of $L(a)$ and construct
a solution of~\eqref{fpr3} satifying \eqref{prlm} by a fixed point argument at infinity.  This idea was
previously used in \cite{SoupletW}, which introduced the
{\it inverted} profile equation  for this purpose.

More precisely, if $f$ is the profile of a radially symmetric self-similar solution, we set
\begin{equation} \label{fIP1} 
w(s) = s^{-\frac{1} {\alpha}}f  \Bigl(  \frac{1} {\sqrt{s}}  \Bigr),
\end{equation}
for $0 < s < \infty$.
The profile equation for $f$, i.e.~\eqref{fpr3},
is equivalent to the following equation for $w$,
\begin{equation} \label{IPE}
4s^2w''(s) + 4\gamma sw'(s) -w'(s) - \beta w(s) + |w(s)|^\alpha w(s) = 0,
\end{equation}
where $\beta$ is given by \eqref{bta} and
\begin{equation} \label{gma} 
\gamma = \frac{2}{\alpha} - \frac{N-4}{2}.
\end{equation} 
If $f = f_a$, we see that 
\begin{equation}\label{fILA}
\begin{cases}
\displaystyle a = \lim_{r \to 0}f_a(r) = \lim_{s\to\infty}s^{\frac{1}{\alpha}}w(s),\\
\displaystyle L(a) = \lim_{r\to\infty}r^{\frac{2}{\alpha}}f_a(r) = \lim_{s\to 0}w(s) = w(0).
\end{cases}
\end{equation}
For an arbitrary solution $w: (0,\infty) \to \R$ of \eqref{IPE}, it is not clear that
 $\lim_{s\to 0}w(s) = w(0)$ exists.
In spite of the highly singular nature of equation \eqref{IPE} at $s = 0$,
global regular solutions $w: [0,\infty) \to \R$ can be constructed
given any fixed values of $w(0)$ and $w'(s_0)$ for sufficiently small
$s_0 > 0$.  See \cite{SoupletW} and Section~\ref{sIPE} below for the details.

Consequently, in order to show that the successive maxima of $ | L (a)| $ 
on the intervals $[ a_m , a _{ m+1 } ]$ tend to infinity,
it suffices to show the existence of a sequence  $\MUU_m \to \infty$ such that for every $m$
 there exists a
solution $w_m$ of \eqref{IPE} with $w(0) =  \MUU_m $, which has precisely
$m$ zeros, and such that $\lim_{s\to\infty}s^{\frac{1}{\alpha}}w_m(s)$
exists and is finite.

It is therefore important to understand the asymptotic behavior
as $s \to \infty$
of solutions to \eqref{IPE}.  It was shown in \cite{SoupletW}
that if $\alpha < \frac{4}{N-2}$, then $\lim_{s\to \infty}w(s)$
always exists and is a stationary solution of \eqref{IPE}.
If $\beta > 0$, which corresponds to $\frac{2}{N-2} < \alpha < \frac{4}{N-2}$,
then $w(s)$ must tend to either $0$ or $\pm\beta^{\frac{1}{\alpha}}$.
In the case $\beta \le 0$, i.e.
$0 < \alpha \le \frac{2}{N-2}$, all solutions must satisfy $w(s) \to 0$
as $s \to \infty$.

The possible asymptotic
behaviors of $w(s)$ can be determined from the following partially formal argument.
Setting $w(s) = z(t)$ where $t = \log s$, we obtain
\begin{equation}
4z''(t) + 4(\gamma - 1)z'(t) - e^{-t}z'(t) - \beta z(t) + |z(t)|^\alpha z(t) = 0.
\end{equation}
Deleting the term $e^{-t}z'(t)$ on the (intuitive) basis that this
will be negligible as $t \to \infty$ yields the autonomous differential equation
\begin{equation}
\label{autoode}
4z''(t) + 4(\gamma - 1)z'(t) - \beta z(t) + |z(t)|^\alpha z(t) = 0.
\end{equation}
Note that $\gamma -1>0$ if $(N-2) \alpha < 4$.
In the case $\beta < 0$, a standard
phase plane stability analysis  shows that all solutions
of~\eqref{autoode} have one of the following two asymptotic
behaviors for large $t > 0$:
\begin{equation}
\label{fastslow1}
\begin{cases}
z(t) \sim  z_1(t) = e^{-\frac{t}{\alpha}}\\
z(t) \sim  z_2(t) =
\begin{cases}
 e^{-t\big(\frac{1}{\alpha}-\frac{N-2}{2}\big)} & N \ge 3\\
 te^{-\frac{t}{\alpha}} & N = 2 .
\end{cases}
\end{cases}
\end{equation}
This translates, for solutions of \eqref{IPE}, into the following
two possible asymptotic behaviors: 
\begin{equation}
\label{fastslow2}
\begin{cases}
w(s) \sim \phi _1(s) = s^{-\frac{1}{\alpha}}\\
w(s) \sim \phi _2(s) = 
\begin{cases}
s^{-\frac{1}{\alpha}+\frac{N-2}{2}} & N \ge 3\\
s^{-\frac{1}{\alpha}}\log s & N = 2 .
\end{cases}
\end{cases}
\end{equation}
If $\beta =0$, $w$ exhibits analogous behaviors with $\phi _1 (s) = s^{-\frac {1} {\alpha }}$ and $\phi _2 (s) = (\log s)^{-\frac {1} {\alpha }}$. If $\beta >0$, the corresponding asymptotic behaviors of $w$
are given by  $\phi _1(s) = s^{-\frac{1}{\alpha}}$ and
$\phi _2(s) = \beta^{\frac{1}{\alpha}}$.  

It is now relatively clear what must be done.  These asymptotic
behaviors must be proved for solutions $w(s)$ of~\eqref{IPE}, and the existence of solutions
with asymptotic behavior like $\phi _1(s)$ with an arbitrarily large
number of zeros must be established. It turns out that for $N\ge 2$ the  long time behavior determined by $\phi _2$, which represents slower decay than the one determined by  $\phi _1$, is stable, and that the desired asymptotic behavior is unstable.
Our approach to proving the existence of solutions with the behavior determined by  $\phi _1$ is a shooting argument. More precisely, we show the existence of solutions with the asymptotic behavior given by   $\phi _2$ with arbitrarily large number of zeros. Solutions with  long time behavior given by $\phi _1$ are found in the transitional regions where the number of zeros increases. 

Solutions of~\eqref{IPE} with the asymptotic behavior $\phi _2$ give rise formally to self-similar
solutions of~\eqref{NLHE} with singular profiles.
It turns out that these are genuine weak solutions of~\eqref{NLHE} when $\beta \ge 0$.  

\begin{theorem} \label{eMain2} 
Assume $N\ge 3$ and let $ \frac {2} {N-2} \le \alpha < \frac {4} {N-2}$. There exist an integer $ \overline{m}  \ge 0$ and an increasing sequence $(  \overline{\MUU}  _m)  _{ m \ge  \overline{m}  } \subset (0, \infty )$, $\overline{\MUU}  _m \to \infty $ as $m\to \infty $, such that for all $m\ge  \overline{m} $ there exists a radially symmetric profile $ h _m$ which is a solution of~\eqref{fpr2}  in the sense of distributions, has exactly $m$ zeros, is regular for $r>0$,  has the singularity
\begin{equation} \label{feMain2} 
\begin{cases} 
r^{\frac {2} {\alpha }}  h _m (r) \to (-1)^m \beta^{\frac {1} {\alpha }} &  \frac {2} {N-2} < \alpha < \frac {4} {N-2} \\
 |\log r|^{\frac {1} {\alpha }}r^{\frac {2} {\alpha }}  h _m (r) \to (-1)^m ( {\frac {2} {\alpha }} ) ^{\frac {2} {\alpha }} & \alpha = \frac {2} {N-2} 
\end{cases} 
\end{equation} 
as $r\downarrow 0$, and the asymptotic behavior $r^\frac {2} {\alpha } h_m (r) \to   \overline{\MUU}  _m$ as $r \to \infty $.

If $u$ is given by~\eqref{fpr1} with $f= h_m$, and $\DI=  \overline{\MUU}  _m |x|^{-\frac {2} {\alpha }}$, then $u$ is a solution of~\eqref{NLHE}, and also of~\eqref{fNZ2}, in $C((0,\infty ), L^p (\R^N ) + L^q (\R^N ) ) $ where $p, q$ satisfy
\begin{equation} \label{eSP1:5} 
\begin{cases} 
1\le p< \frac {N\alpha } {2(\alpha +1)} < \frac {N\alpha } {2} <q &  \text{if }   \alpha >\frac {2} {N-2} \\
1= p < \frac {N\alpha } {2} < q &  \text{if }  \alpha =\frac {2} {N-2} . \\
\end{cases} 
\end{equation} 
Moreover,  $u$ satisfies~\eqref{fEP2} and~\eqref{fEP3}.
 In addition, the map $t \mapsto u (t) - e^{t \Delta } \DI $ is in $ C([0,\infty ), L^r (\R^N ) )$ for all   $\frac {N\alpha } {2(\alpha +1) } < r < \frac {N\alpha } {2}$ if $\alpha > \frac {2} {N-2}$, and  in $ C([0,\infty ), L^1 (\R^N ) + L^r (\R^N ) )$ for all  $r>1$ if $\alpha = \frac {2} {N-2}$.
\end{theorem} 

If $\alpha <\frac {2} {N-2}$, the singular profiles do not give rise to weak solutions of~\eqref{NLHE}, see Remark~\ref{eRemz}.

In the case $N=1$, the conclusion of Theorem~\ref{eMain1} holds. However, the asymptotic analysis of  $w(s) $ is much simpler, since the profile equation~\eqref{fpr3} has no singularity at $r=0$. 
This last property also allows us to construct odd profiles. 
Our result in this case is the following.

\begin{theorem} \label{eMain3} 
Assume $N =1$ and let $\alpha >0$ and  $ \MUU \not =  0$. 
\begin{enumerate}[{\rm (i)}] 

\item \label{eMain3:1}  
There exists $m_0\ge 0$ such that for all $m \ge m_0$ there exist at least two different, even, regular self-similar solutions of \eqref{NLHE} with initial value $\DI=  \MUU  |x|^{-\frac {2} {\alpha }}$ in the sense~\eqref{fEP2}, and also~\eqref{fEP3} if $\alpha > 2$,  and whose profiles have exactly $2 m$ zeros.

\item \label{eMain3:2}  
There exists $m_0\ge 0$ such that for all $m \ge m_0$ there exist at least two different, odd, regular self-similar solutions of \eqref{NLHE} with initial value $\DI=  \MUU  |x|^{-\frac {2} {\alpha } -1} x$ in the sense~\eqref{fEP2}, and also~\eqref{fEP3} if $\alpha >2$,  and whose profiles have exactly $2m +1$ zeros.
\end{enumerate} 

\noindent These solutions satisfy $u\in C^1((0,\infty ), L^r (\R ) )$ and  $u _{ xx },  |u|^\alpha u\in C((0,\infty ), L^r (\R ) )$ for all $r\ge 1$, $r> \frac {\alpha } {2}$. 
Furthermore, if  $\alpha >2$,  they are solutions of the integral equation~\eqref{fNZ2},  where each term is in $C((0,\infty ), L^r (\R ) )$ for all $r> \frac {\alpha } {2} > 1$. Moreover, the map $t \mapsto u (t) - e^{t \Delta } \DI $ is in $ C([0,\infty ), L^r (\R ) )$ for all $r \ge 1$ such that $\frac {\alpha } {2(\alpha +1) } < r < \frac {\alpha } {2}$.

\end{theorem} 

The rest of the paper is organized as follows. The proof of Theorem~\ref{eMain1} is given in Section~\ref{sProof}. 
This proof depends on several intermediate results, which are stated in Section~\ref{sProof}, but whose proofs are deferred until Sections~\ref{sIPE} to~\ref{sFNL}.
Theorems~\ref{eMain2} and~\ref{eMain3} are proved in Sections~\ref{sSinPro} and~\ref{sDim1}, respectively.
Finally, in Appendix~\ref{sZEROS} we make a remark about the smallest value of $m_0$ possible in Theorems~\ref{eMain1} and~\ref{eMain3}.

\section{Proof of Theorem~$\ref{eMain1}$} \label{sProof} 

In this section, we prove Theorem~\ref{eMain1}. The bulk of the work is the construction of the appropriate profiles, as stated in the following theorem.

\begin{theorem} \label{eMT2} 
Suppose~\eqref{fNZ1}. 
For every $\MUU >0$, there exist an integer $m_\MUU$ and  four sequences $(a _{ \MUU , m }^\pm ) _{ m \ge  m _ \MUU  }\subset (0,\infty )$ and $(b_{ \MUU , m }^\pm ) _{ m \ge  m _\MUU   }\subset (-\infty , 0)$ such that the following properties hold, with the notation~\eqref{prlm}-\eqref{nz}.
\begin{enumerate}[{\rm (i)}] 

\item \label{eMT2:1} 
$ a_{ \MUU , m }^\pm \to \infty $ and $ b_{ \MUU , m }^\pm \to - \infty $ as $m\to \infty $. 

\item \label{eMT2:2} 
$L (a_{ \MUU , m }^\pm ) =L (b_{ \MUU , m }^\pm ) =\MUU $ for all $m\ge  m _ \MUU $.

\item \label{eMT2:3} 
For all $m\ge  m _ \MUU  $, $a_{ \MUU , m }^- < a_{ \MUU , m }^+$ and $ N (a_{ \MUU , m }^\pm ) =2m$. 

\item \label{eMT2:4} 
For all $m\ge  m _ \MUU $, $b_{ \MUU , m }^+ < b_{ \MUU , m }^-$ and $ N (b_{ \MUU , m }^\pm ) =2m +1$. 
\end{enumerate} 
\end{theorem} 

The proof of Theorem~\ref{eMT2} depends on a series of propositions, which are stated below and proved in the subsequent sections. In the last part of this section, we give the proof of Theorem~\ref{eMT2} assuming these propositions and, at the very end of this section, obtain Theorem~\ref{eMain1} as a consequence.

As mentioned in Section~\ref{Intro}, our approach is based on the study of the inverted profile equation~\eqref{IPE}. It is convenient to set
\begin{equation} \label{fNot1} 
g(s) = -\beta s+  |s|^\alpha s, \quad s\in \R
\end{equation} 
and 
\begin{equation} \label{fNot2} 
G(s) = -\frac {\beta } {2} s^2+ \frac {1} {\alpha +2} |s|^{\alpha + 2}, \quad s\in \R
\end{equation} 
We collect in the following proposition an existence result for solutions of~\eqref{IPE} with appropriate initial conditions, as well as several properties of these solutions. These results are in part taken from~\cite{SoupletW}, and the detailed proof is given in Section~\ref{sIPE}. 

\begin{proposition} \label{eEX1:b1}  
Assume $N\ge 1$ and $\alpha >0$. 
There exists $\NUU >0$ such that for all ${\MUU}\ge \MUUz$, there exists a  unique solution $w_{\MUU} \in C^1([0,\infty )) \cap C^2(0,\infty )$ of~\eqref{IPE} satisfying
\begin{equation}  \label{fTT1:b1}
w_{\MUU} (0)= \MUU , \quad {\textstyle{ w_{\MUU}'(  \NUU \MUU^{- \alpha } ) }} =0  , \quad  \| w_\MUU  \| _{ L^\infty  (0,  \NUU \MUU^{- \alpha } ) } \le 10 \MUU  .
\end{equation} 
Moreover, $w _\MUU$ satisfies the  following properties.
\begin{enumerate}[{\rm (i)}] 

\item  \label{eEX1:6:b1} Given any $T>0$, $w_{\MUU}$ depends continuously on ${\MUU}$ in $C^1([0,T])$.

\item  \label{eEX1:2:b1} If
\begin{equation} \label{fENE1:b2} 
H_{\MUU} (s)= 2s^2 w_{\MUU} '(s)^2 + G (w_\MUU (s))
\end{equation} 
with the notation~\eqref{fNot1}-\eqref{fNot2},  then
\begin{equation} \label{fENE3:b1} 
H_{\MUU} '(s) = w_{\MUU} '(s)^2 [1 - 4(\gamma -1) s]
\end{equation} 
for all $s\ge 0$. 

\item  \label{eEX1:1:b1} If $(N-2) \alpha < 4$, then $w_{\MUU}(s)$ and $sw_{\MUU}'(s)$ are bounded as $s\to \infty $. 

\item  \label{eEX1:2:b2} If $(N-2) \alpha < 4$, then $H_{\MUU} (s) $ defined by~\eqref{fENE1:b2}  has a finite limit as $s\to \infty $.

\item  \label{eEX1:3:b1} 
If $(N-2) \alpha < 4$, then $s w_{\MUU}' (s) \to 0$ as $s\to \infty $.
Moreover, if $\beta \le 0$, then $w_{\MUU}(s) \to 0$ as $s\to \infty $; and if $\beta >0$, then $w_{\MUU}(s)$ converges  to either $0$ or $\pm \beta ^{\frac {1} {\alpha  }}$, for $\beta >0$. 

\item  \label{eEX1:2:6} $w _\MUU$ has a finite number of zeros. 

\item  \label{eEX1:2:7} If $w_\MUU (s) \to 0$ as $s\to \infty $, then $w _\MUU ' $ has a finite number of zeros, and $w_\MUU w_\MUU ' <0$ for $s$ large.

\end{enumerate} 
\end{proposition} 

\begin{remark} 
There is an arbitrary choice in Proposition~\ref{eEX1:b1}. Indeed, as the proof of the proposition shows, $\NUU$ can be any sufficiently small positive number.   In fact, all that is needed in the subsequent arguments is a collection of solutions $(w_{\MUU} ) _{ \MUU \ge 1 }\subset  C^1([0,\infty )) \cap C^2(0,\infty )$ of~\eqref{IPE} such that $w_\MUU (0) =\MUU$ and Property~\eqref{eEX1:6:b1} of Proposition~\ref{eEX1:b1} is true. 
\end{remark} 

The  key  ingredient in the proof of Theorem~\ref{eMT2} is that the number of zeros of $w_\MUU$ is arbitrarily large for large $\MUU$. 
This is stated in the following proposition. 

\begin{proposition} \label{fFIn1} 
Assume $N\ge 1$, $\alpha >0$, and let  $( w_\MUU ) _{ \MUU \ge 1 }$ be the collection of solutions of~\eqref{IPE} given by Proposition~$\ref{eEX1:b1}$.
Given any $T>0$ and $m\in \N$, there exists $ \overline{\MUU } \ge  \MUUz$ such that if $\MUU \ge  \overline{\MUU} $, then $w_\MUU$ has at least $m$ zeros on $[0,T]$. 
In particular, if 
\begin{equation} \label{fFIn1:1} 
 \text{$ \wN (\MUU ) $ is the number of zeros of $w_\MUU $ on $(0,\infty )$, }
\end{equation}  
then $ \wN (\MUU ) \to \infty  $ as $\MUU \to \infty $.
\end{proposition} 

The proof of Proposition~\ref{fFIn1} is given in Section~\ref{sAB}.
The second ingredient we use is a classification of the possible asymptotic behaviors of $w_\MUU (s)$ as $s\to \infty $.
We let
\begin{equation}  \label{fPR1}
\lambda _1 = \frac {1} {\alpha } , \quad \lambda _2 = \frac {1} {\alpha } - \frac {N-2} {2} . 
\end{equation} 
Note that $\lambda _1> \lambda _2$ if $N\ge 3$ and $\lambda _1= \lambda _2$ if $N=2$. Moreover,
\begin{equation} \label{fNZE3} 
\beta = - 4 \lambda _1 \lambda _2,\quad \gamma = 1 + \lambda _1 + \lambda _2
\end{equation} 
where $\beta $ and $\gamma  $ are defined by~\eqref{bta} and~\eqref{gma}. 
We also define
\begin{equation} \label{fPR2}
\phi _1 (s) = s^{- \lambda _1}
\end{equation} 
and
\begin{equation} \label{fPR3}
\phi _2 (s) = 
\begin{cases} 
1 & N\ge 3  \text{ and }\alpha > \frac {2} {N-2}  \\
(\log s )^{-\frac {1} {\alpha }} & N\ge 3  \text{ and }\alpha = \frac {2} {N-2} \\
s^{- \lambda _2} & N\ge 3 \text{ and }\alpha < \frac {2} {N-2} \\
s^{- \lambda _1} \log s & N=2.
\end{cases} 
\end{equation} 
The following result shows that  as $s\to \infty $, the solutions $w_\MUU (s)$ of~\eqref{IPE}  given by Proposition~\ref{eEX1:b1} decay either as $\phi _2$ (slow decay), or else as $\phi _1$ (fast decay).

\begin{proposition} \label{eAB1} 
Assume~\eqref{fNZ1}, and let  $( w_\MUU ) _{ \MUU \ge 1 }$ be the collection of solutions of~\eqref{IPE} given by Proposition~$\ref{eEX1:b1}$.
If $\phi _1$ and $\phi _2$ are defined by~\eqref{fPR2} and~\eqref{fPR3}, respectively, then 
the following properties hold.
\begin{enumerate}[{\rm (i)}] 

\item \label{eAB1:1} 
The limit
\begin{equation} \label{eAB2} 
\Lt(\MUU)=\lim_{s\to \infty} \frac {w_\MUU(s) } {\phi _2(s)}
\end{equation} 
exists and is finite. 

\item \label{eAB1:2} 
If $\Lt( \MUU )=0$,  then
\begin{equation} \label{eAB4} 
\Ltu( \MUU )=\lim_{s\to \infty} \frac {w_ \MUU (s) } {\phi _1(s)}
\end{equation} 
exists and is finite, and $\Ltu( \MUU )\ne 0$. 
\end{enumerate} 
\end{proposition}

Proposition~\ref{eAB1} in the case $\beta >0$ follows from the results in~\cite{SoupletW}. For $\beta \le 0$, the proof requires the separate study of the various cases $N\ge 3$ and $\beta <0$; $N\ge 3$ and $\beta =0$; $N=2$. The proof is given in Section~\ref{sec:auxiliary}.

Proposition~\ref{eAB1} is relevant for the following reason. 
To prove Theorem~\ref{eMT2} we need to find $a>0$ such that $L (a) = \MUU$. 
For the solution $w$ of the inverted profile equation~\eqref{IPE} defined by~\eqref{fIP1} in terms of the profile $f_a$, this means, by~\eqref{fILA},  that $w (0)= \MUU $. The existence of solutions $w$ of~\eqref{IPE} such that $w (0)= \MUU $ follows from Proposition~\ref{eEX1:b1}. However, in order that these solutions correspond to a profile $f_a$, we must ensure that $f$ defined by~\eqref{fIP1} in terms of $w$ has a finite limit as $r\to 0$. This last condition is equivalent to the fact that $ \frac {w(s)} {\phi _1(s)}=s^{\frac {1} {\alpha }} w(s)$ has a finite limit as $s\to \infty $. By Proposition~\ref{eAB1}, we see that we should look for solutions $w_\MUU(s)$ of~\eqref{IPE} such that  $\Lt( \MUU )=0$.

The last ingredient needed in the proof of Theorem~\ref{eMT2} concerns the local behavior of $ \wN ( \MUU  ) $, the number of zeros of $w_\MUU$. 

\begin{proposition}  \label{eEX2} 
Assume~\eqref{fNZ1}. Let $\Lt $ be defined by~\eqref{eAB2} and $\wN$ by~\eqref{fFIn1:1}, and consider $  \overline{\MUU}  \ge  \MUUz$.
\begin{enumerate}[{\rm (i)}] 

\item \label{eEX2:1}
If $\Lt( \overline{\MUU} )\ne 0$, then  there exists $\delta >0$ such that  if $\MUU \ge \MUUz$ and  $ | \MUU -  \overline{\MUU } | \le \delta $, then   $ \wN  ( \MUU ) = \wN  (  \overline{\MUU}  ) $.

\item \label{eEX2:2} 
If $\Lt( \overline{\MUU} )= 0$, then there exists $\delta >0$ such that  if $\MUU \ge \MUUz$ and $ | \MUU -  \overline{\MUU } | \le \delta $, then   either $ \wN  ( \MUU ) = \wN  (  \overline{\MUU}  ) $ or else $ \wN  ( \MUU ) = \wN  (  \overline{\MUU}   ) +1 $ and $\Lt ( \MUU) \not = 0$.
\end{enumerate} 
\end{proposition} 

The proof of Proposition~\ref{eEX2} is carried out in Section~\ref{sFNL}. 

We are now in a position to complete the proof of Theorem~\ref{eMT2}, assuming Propositions~\ref{eEX1:b1}, \ref{fFIn1}, \ref{eAB1} and~\ref{eEX2}. 
We first give a lemma, which will also be used in the proofs of Theorems~\ref{eMain2} and~\ref{eMain3}. 

\begin{lemma} \label{eTPA15} 
Assume~\eqref{fNZ1}.  Let $\Lt $ be defined by~\eqref{eAB2} and $\wN$ by~\eqref{fFIn1:1},  and set $ \overline{m} =  \wN  (2) +1$. It follows that there exists an increasing sequence $( \MUU _m)_{m\ge \overline{m}} \subset [ \MUUz , \infty )$ with the following properties.
\begin{enumerate}[{\rm (i)}] 

\item \label{eTPA15:1} 
 $ \MUU  _m \to \infty $ as $m\to \infty $,  $\Lt( \MUU  _m)=0$, $ \wN ( \MUU _m) =m$.
 
 \item \label{eTPA15:2} 
 For every $m\ge  \overline{m} $, there exists $\MUU \in ( \MUU_m, \MUU  _{ m+1 } )$ such that $ \wN ( \MUU ) =m +1$ and $\Lt( \MUU  ) \not = 0$.
\end{enumerate} 
\end{lemma}  

\begin{proof} 
Given $m\ge  \overline{m} $, let
\begin{equation} \label{eNZ1:9}
E_m= \{ \MUU \ge 2 ;\,  \wN ( \MUU )\ge  m+1  \}.
\end{equation} 
It follows from Proposition~\ref{fFIn1} that $E_m \not = \emptyset$, and we define
\begin{equation} \label{eNZ1:10}
\MUU _m = \inf  E _m \in [2, \infty ). 
\end{equation} 
By Proposition~\ref{eEX2}, $\wN ( \MUU) \le  \wN (2)+1=\overline{m} \le m$ for $\MUU\ge 2$ close to $2$. Thus we see that $\MUU_m > 2$.   
Suppose that $ \wN ( \MUU_m ) \le m-1$. It follows from Proposition~\ref{eEX2} that $\wN (\MUU) \le m$ for $\MUU $ close to $\MUU_m$, contradicting~\eqref{eNZ1:10}. 
Suppose now that $\wN ( \MUU_m ) \ge m+1$. Applying again Proposition~\ref{eEX2}, we deduce that $\wN ( \MUU) \ge m+1$  for $\MUU $ close to $\MUU_m$, contradicting again~\eqref{eNZ1:10}. 
So we must have $ \wN  (\MUU _m) =m$. Moreover, if $\Lt( \MUU  _m)\not = 0$, then $\wN ( \MUU) = m$  for $\MUU $ close to $\MUU_m$, contradicting once more~\eqref{eNZ1:10}; and so, $\Lt( \MUU  _m) = 0$.
We finally prove that $ \MUU  _m \to \infty $ as $m\to \infty $. Otherwise,  there exist $ \overline{\MUU}\in  [2, \infty )$ and a sequence $m_k \to \infty $ such that $\MUU _{m_k} \to  \overline{\MUU} $ as $k \to \infty $ and  $ \wN ( \MUU _{m_k}) =m_k $. 
It follows from Proposition~\ref{eEX2} that $ \wN  (\MUU _{m_k}) \le  \wN  ( \overline{\MUU} ) +1$ for all sufficiently large $k$, which is absurd. This proves Property~\eqref{eTPA15:1}. 

Moreover, $E _{ m+1 } \subset E_m$, so that $\MUU  _{ m+1 } \ge \MUU _m$. Since $ \wN ( \MUU  _{ m+1 }) \not =  \wN ( \MUU _m) $, we see that $\MUU  _{ m+1 } > \MUU _m $, so that $( \MUU _m)_{m\ge \overline{m}} $ is  an increasing sequence.

We obtain Property~\eqref{eTPA15:2} by the following argument. (The authors thank the referee for having improved the original argument.) It follows from Proposition~\ref{eEX2}~\eqref{eEX2:2} 
that for $\MUU $ close to $  {\MUU} _m$, we have either $\wN ( \MUU) =m$ or else $\wN ( \MUU) = m+1$ and $ \Lt ( \MUU ) \not = 0$. 
Since, by~\eqref{eNZ1:10}, there exists $\MUU > \MUU_m$ arbitrarily close to $ {\MUU}_m $ (hence less than $ {\MUU} _{ m+1 } $) such that $ \wN ( \MUU ) =m+1$,  Property~\eqref{eTPA15:2} follows.
\end{proof} 

\begin{proof} [Proof of Theorem~$\ref{eMT2}$]
We first prove that there exists a sequence $(\alpha  _m) _{ m\ge  \overline{m}  } \subset (0,\infty ) $ satisfying (with the notation~\eqref{prlm}-\eqref{nz}) 
\begin{gather} 
\alpha _m \goto _{ m\to \infty  }\infty  \text{ and } (-1)^m  L (\alpha _m ) \goto _{ m\to \infty  }\infty   \label{eNZ3:1} \\
N (\alpha _m ) = m . \label{eNZ3:3} 
\end{gather} 
To see this, we consider the collection  $( w_\MUU ) _{ \MUU \ge 1 }$  of solutions of~\eqref{IPE} given by Proposition~\ref{eEX1:b1},  and the sequence $(\MUU _m) _{ m\ge  \overline{m}  }$ of Lemma~\ref{eTPA15}. 
Since $\Lt ( \MUU _m )=0$, it follows from Proposition~\ref{eAB1}~\eqref{eAB1:2} that $\Ltu ( \MUU _m )$ is well defined and $\Ltu ( \MUU _m )\not = 0$. On the other hand, $ \wN  (\MUU _m) =m$ and $w _{ \MUU _m } (0)= \MUU_m >0$, so that $(-1)^m w _{ \MUU _m } (s) >0$ for $s$ large.  Thus we see that  $ (-1)^m \Ltu ( \MUU _m ) > 0$. 
Setting $\alpha   _m=  (-1)^m  \Ltu ( \MUU _{ m } )$, it follows  from \eqref{fILA} and \eqref{eAB4} that 
\begin{equation} \label{fTPA15} 
f _{ \alpha _m }(r)=  (-1)^m  r^{-\frac {2} {\alpha }} w _{ \MUU _m } ( r^{-2} )
\end{equation} 
for all $r>0$.  Moreover $N (\alpha _m) =  \wN ( \MUU  _m ) =m $, $L (\alpha _m) = (-1)^m w _{ \MUU _m  } (0)= (-1)^m  \MUU_m $.
In particular, $(-1)^m  L (\alpha _m) \to \infty $ as $m\to \infty $, and since $L$ is continuous $[0, \infty ) \to \R$, we see that $\alpha _m \to \infty $ as $m\to \infty $. This proves~\eqref{eNZ3:1} and~\eqref{eNZ3:3}.

Fix $\MUU >0$. By~\eqref{eNZ3:1}, we may choose an integer $m_ \MUU  \ge   \overline{m} $ sufficiently large so that
\begin{equation} \label{fMT2:1} 
L ( \alpha _{2m}) >  \MUU  \text{ and } L ( \alpha _{2m -1}) <-  \MUU    \text{ for all } m\ge m_ \MUU  .
\end{equation} 
Gven $m\ge m_\MUU$, $L$ has a zero to the right of $\alpha  _{ 2m }$ by~\eqref{eNZ3:1}. Moreover, $L(0) =0$. By continuity of $L : \R \to \R$, there exist 
\begin{equation}  \label{fMT2:2} 
0\le  a_m^- < \alpha _{2m} < a_m ^+  
\end{equation} 
such that 
\begin{equation}  \label{fMT2:3} 
 L ( a_m^- )= L (a_m ^+ ) =0  \text{ and } L(a) >0  \text{ for all }   a_m^- < a < a_m^+ .
\end{equation} 
We observe that necessarily
\begin{equation} \label{fSSP1} 
a_m^- \goto  _{ m\to \infty  } \infty .
\end{equation} 
Indeed, if not there exist $K>0$ and a  subsequence $a_{m_j }^-$ such that $ a _{ m_j  }^- \le K $. 
By~\eqref{fMT2:3} we have $L(a) >0$  for all $K<a< \alpha  _{ 2m_j }$. Since $ \alpha  _{ 2m_j } \to \infty $ as $j\to \infty $, this implies that $L(a) > 0$ for all $a>K$. This contradicts~\eqref{eNZ3:1} and establishes~\eqref{fSSP1}. 
Next, since $N$ is locally constant where $L \not = 0$ (by~\cite[Proposition~3.7.b]{Weissler6}), it follows that $N$ is constant on $(a_m^- , a_m^+)$. 
Next,  $N (\alpha _{2m}) =2m  $ (by~\eqref{eNZ3:3}) and $\alpha _{2m} \in (a_m^- , a_m^+) $ (by~\eqref{fMT2:2}), so that
\begin{equation}  \label{fMT2:4} 
N ( a) = 2m  \text{ for all } a_m^- < a < a_m^+ .
\end{equation} 
From~\eqref{fMT2:1} and~\eqref{fMT2:3}, we deduce that there exist 
\begin{equation} \label{fMT2:6} 
a_m^- < a _{  \MUU  , m }^- < \alpha _{2m} < a _{  \MUU  ,m }^+ < a_m^+
\end{equation} 
such that $L ( a _{  \MUU  , m }^\pm ) =  \MUU  $.
Moreover, it follows from~\eqref{fMT2:4} and~\eqref{fMT2:6} that $N ( a _{  \MUU  ,m }^\pm ) = 2 m$,
 and from~\eqref{fSSP1} and~\eqref{fMT2:6} that $a _{  \MUU  ,m }^\pm \to \infty $ as $m\to \infty $. 
Thus we see that the sequences $(a _{  \MUU  ,m }^\pm ) _{ m \ge m_ \MUU   }$ satisfy all the conclusions of the theorem. 

By considering $\alpha  _{ 2m+1 } \to \infty $ (instead of  $\alpha _{2m}$), one constructs as above a sequence $( \widetilde{b} _{  \MUU  , m }^\pm ) _{ m \ge  m_ \MUU   }\subset (0, \infty )$ such that $  \widetilde{b} _{  \MUU  , m }^\pm \to  \infty $ as $m\to \infty $, $L (  \widetilde{b} _{  \MUU  , m }^\pm ) = -  \MUU  $, $ \widetilde{b} _{  \MUU  , m }^- <  \widetilde{b} _{  \MUU  , m }^+$ and $ N (  \widetilde{b} _{  \MUU  , m }^\pm ) =2m +1$. 
The result follows by letting $b  _{  \MUU  , m }^\pm = -  \widetilde{b}  _{  \MUU  , m }^\pm$.
\end{proof} 

With the completion of the proof of Theorem~\ref{eMT2}, the existence of the appropriate profiles has been established. It remains to show that the resulting solutions given by formula~\eqref{fpr1} have the properties required by Theorem~\ref{eMain1}. To this end, we prove the following lemma, which includes the case of non-radially symmetric profiles.

\begin{lemma} \label{eQSol} 
Let $\alpha >0$ and let $f\in C^2 (\R^N )$ be a solution of~\eqref{fpr2} such that $ |f(x)| +  |x\cdot \nabla f(x)|\le C (1+ |x|^2)^{-\frac {1} {\alpha }}$ and $ |x|^\frac {2} {\alpha } f(x)- \omega (x) \to 0$ as $ |x| \to \infty $, for some $\omega \in C (\R^N \setminus \{0\})$ homogeneous of degree $0$. If $u (t, x) $ is defined by~\eqref{fpr1} for all $t>0$ and $x\in \R^N $, then  $u\in C^1((0,\infty ), L^r (\R^N ) )$ and $\Delta u,  |u|^\alpha u\in C((0,\infty ), L^r (\R^N ) )$ for all $r\ge 1$, $r> \frac {N\alpha } {2}$, and $u$ satisfies~\eqref{fEP2} with $\DI= \omega (x) |x|^{-\frac {2} {\alpha }}$.
If in addition $\alpha >\frac {2} {N}$, then $u$ also satisfies~\eqref{fEP3}, and $u $ is a solution of~\eqref{fNZ2} where each term is in $C((0,\infty ), L^r (\R^N ) )$ for  $r> \frac {N\alpha } {2}$. Moreover, the map $t \mapsto u (t) - e^{t \Delta } \DI $ is in $ C([0,\infty ), L^r (\R^N ) )$ for all $r \ge 1$ such that $\frac {N\alpha } {2(\alpha +1) } < r < \frac {N\alpha } {2}$.
\end{lemma} 

\begin{remark} 
If $\alpha \le \frac {2} {N}$, then $\DI \not \in L^1 _\Loc (\R^N ) $ (except if $\omega \equiv 0$), so that the first term on the right hand side of~\eqref{fNZ2} does not make sense, hence~\eqref{fNZ2} altogether, does not make sense.
\end{remark} 

\begin{proof} [Proof of Lemma~$\ref{eQSol}$]
Since $f\in C^2 (\R^N )$ is a solution of~\eqref{fpr2}, it follows from formula~\eqref{fpr1} and elementary calculations that $u\in C^2((0,\infty )\times \R^N )$ is a solution of~\eqref{NLHE} on   $(0,\infty )\times \R^N$. Moreover, it is not difficult to show by dominated convergence that $u\in C ((0,\infty ), L^r (\R^N ) )$  for all $r\ge 1$, $r>\frac {N\alpha } {2}$. Differentiating~\eqref{fpr1} with respect to $t$, we see that
\begin{equation}  \label{feQSol4} 
u_t= - \frac {1} {\alpha t}u - \frac {1} {2t} x\cdot \nabla u ,
\end{equation} 
and since $ |f(x)| +  |x\cdot \nabla f(x)|\le C (1+ |x|^2)^{-\frac {1} {\alpha }}$, we have $u\in C^1 ((0,\infty ), L^r (\R^N ) )$  for all $r\ge 1$, $r>\frac {N\alpha } {2}$.
Moreover, $ |f|^{\alpha +1}\le C (1+ |x|^2)^{-\frac {\alpha +1} {\alpha }}\le C (1+ |x|^2)^{-\frac {1 } {\alpha }}$, from which it follows easily that $ |u|^\alpha u \in C((0,\infty ), L^r (\R^N ) )$ for all $r\ge 1$, $r> \frac {N\alpha } {2}$. The regularity of $\Delta u$ now follows from equation~\eqref{NLHE}. 
Next, $ |f(x)|\le C (1+ |x|^2)^{-\frac {1 } {\alpha }}\le C  |x|^{-\frac {2 } {\alpha }}$, so that $ |u(t,x) |\le C  |x|^{-\frac {2} {\alpha }}$. Moreover given $x\not = 0$,
\begin{equation} \label{fTR1} 
\lim_{t \to 0}u(t,x) = \lim_{t \to 0}t^{-\frac {1} {\alpha }} f \Bigl( \frac { x } {\sqrt t} \Bigr)
= |x|^{-\frac {2} {\alpha }}\lim_{r\to\infty}r^{\frac{2}{\alpha}} f  \Bigl( \frac {x} { |x|} r \Bigr)
= \omega (x) |x|^{-\frac {2} {\alpha }},
\end{equation} 
and property~\eqref{fEP2} follows by dominated convergence.

If $\alpha >\frac {2} {N}$, then $ |\cdot |^{-\frac {2} {\alpha }}\in L^p (\R^N ) +L^q (\R^N ) $ for all $1\le p<\frac {N\alpha } {2} < q $, which implies property~\eqref{fEP3}. Moreover, the regularity of $u$ on $(0,\infty )$ ensures that
\begin{equation}  \label{feQSol1} 
u(t) = e^{(t- \varepsilon )\Delta } u(\varepsilon ) + \int _\varepsilon ^t e^{ (t-s) \Delta }  |u(s)|^\alpha u(s)\, ds
\end{equation} 
for all $0< \varepsilon < t$. Consider now $r \ge 1$ such that $\frac {N\alpha } {2(\alpha +1) } < r < \frac {N\alpha } {2}$. In particular $f\in L^{(\alpha +1) r} (\R^N ) $, therefore,
\begin{equation}  \label{feQSol3} 
\begin{split} 
\| e^{(t-s) \Delta }  |u(s)|^\alpha u(s) \| _{ L^r } &  \le \|   |u(s)|^\alpha u(s) \| _{ L^r }  =   \| u(s)\| _{ L^{ ( \alpha +1) r } }^{\alpha +1}
\\ &=   s^{- \frac {\alpha +1} {\alpha } + \frac {N} {2 r }}  \| f  \| _{ L^{(\alpha +1) r } }^{\alpha +1} .
\end{split} 
\end{equation} 
Note that $- \frac {\alpha +1} {\alpha } + \frac {N} {2 r } >-1$ because $r  < \frac {N\alpha } {2}$.
Applying~\eqref{feQSol3}, one easily passes to the limit in~\eqref{feQSol1} as $\varepsilon \downarrow 0$ and obtain equation~\eqref{fNZ2}. Since the first two terms in~\eqref{fNZ2} are in $C((0,\infty ), L^r (\R^N ) )$ for  $r> \frac {N\alpha } {2}$, so is the integral term.
Finally, that $ u (t) - e^{t \Delta } \DI \in  C([0,\infty ), L^r (\R^N ) )$ for all $r \ge 1$ such that $\frac {N\alpha } {2(\alpha +1) } < r < \frac {N\alpha } {2}$ easily follows  from~\eqref{fNZ2} and~\eqref{feQSol3}.
\end{proof} 

We finally complete the proof of Theorem~\ref{eMain1}. We may assume $\MUU >0$ without loss of generality, and we apply Theorem~\ref{eMT2}. The profiles $f= f_{a^\pm }$ with $a^\pm = a^\pm  _{ \MUU, \frac {m} {2} }$ if $m\ge m_\MUU $ is even and $a^\pm =  b^\pm  _{ \MUU, \frac {m-1} {2} }$ if $m\ge m_\MUU$ is odd, are two different radially symmetric solutions of~\eqref{fpr2} with $m$ zeros. Moreover, $r^{\frac {2} {\alpha }} f_{a^\pm} (r) \to \MUU$ as $r\to \infty $, and it follows from~\cite[Proposition~3.1]{HarauxW} that $ |f_{a^\pm}  (r)| + r  |f'_{a^\pm}  (r)|\le C( 1+ r^2)^{-\frac {1} {\alpha }}$. Theorem~\ref{eMain1} is now an immediate consequence of Lemma~\ref{eQSol}, where $f (x) =f_{a^\pm} ( |x|)$ and $\omega (x) \equiv \MUU$. 

\section{The inverted profile equation} \label{sIPE} 

This section is devoted to the proof of Proposition~\ref{eEX1:b1}, which is based on the results and methods in~\cite{SoupletW}. 

We begin the proof of Proposition~\ref{eEX1:b1} with a fixed-point argument adapted from the proof of~\cite[Proposition~2.5]{SoupletW}. In~\cite{SoupletW},  local existence of solutions to equation~\eqref{IPE} is established for solutions such that $w(0)= \mu $ and $w'(T) = \nu$ where $\mu ,\nu$ and $T$ satisfy certain conditions. In particular, continuous dependence on $\mu $ and $\nu$ is established in $C^1 ([0,T])$ for a fixed $T$. Since this formulation is  different from the formulation of Proposition~\ref{eEX1:b1}, we need to adapt the arguments in~\cite{SoupletW} to our present situation. In doing so, we use  the notation and many of the intermediate results of~\cite{SoupletW}. In particular, we consider the functions 
\begin{equation*} 
K_1 (t)= t^{\gamma -2} e^{\frac {1} {4t}} \int _0^t s ^{- \gamma } e^{- \frac {1} {4s }} ds ,\quad 
K_2 (t)= \sup _{ 0<s<t } \frac {1} {4} s ^{- \gamma } e^{- \frac {1} {4s }} \int _s^t \sigma ^{\gamma -2} e^{\frac {1} {4 \sigma }} d\sigma ,
\end{equation*} 
both continuous on $[0,\infty )$, with $K_1(0)= 4$. 
We recall that if 
\begin{equation}  \label{fLDC2:bz} 
R>0, \quad M>5R,\quad T>0,
\end{equation} 
and
\begin{equation}  \label{fLDC2} 
\frac {R} {M} (1+ K_1 (T) ) + T K_2(T) \sup  _{  | y |\le M }  |g'(  y )| \le 1 ,
\end{equation} 
where $g$ is given by \eqref{fNot1},
then for every $\MUU, \nu \in \R$ such that $ | \MUU |\le R$ and $T^2  | \nu |\le R$, there exists a unique solution $w \in C^1([0, T ]) \cap C^2((0, T ])$ of~\eqref{IPE} satisfying $w (0)= \MUU$, $w ' (T) = \nu$ and $ \| w \| _{ L^\infty  (0,  T) } \le M$.  
See~~\cite[Theorem~2.5, formula~(2.9), Proposition~2.4]{SoupletW}. 

We  fix $t_0 >0$ sufficiently small so that 
\begin{equation} \label{fLDC0} 
\sup _{ 0\le s\le t_0  } K_1(s) \le 5
\end{equation} 
and we set
\begin{equation} \label{fLDC0:1} 
k _0 =   \sup _{ 0\le s\le t_0  } K_2(s) .
\end{equation} 
Note that for all $M>0$
\begin{equation}  \label{fLDC00} 
\sup  _{  |y|\le M }  |g'(y)| \le |\beta |+ (\alpha +1) M^\alpha .
\end{equation} 
We fix  $\NUU >0$ sufficiently small so that 
\begin{gather} 
\NUU \le t_0 \label{fTPA12} \\
 5 \NUU k _0  (  |\beta |  + (\alpha +1) 20^\alpha  )  \le 2 . \label{fLDC1} 
\end{gather} 
(Condition~\eqref{fLDC1} is stronger than what we need at this point of the argument.)
It follows from~\eqref{fLDC0}--\eqref{fLDC1}  that if ${\MUU} \ge \MUUz$ and 
\begin{equation*}
R =  \MUU , \quad M= 10 \MUU , \quad T= \NUU \MUU^{- \alpha } 
\end{equation*} 
then~\eqref{fLDC2:bz}-\eqref{fLDC2}  hold;
and so there exists a unique solution $w_{\MUU} \in C^1([0, \NUU \MUU^{- \alpha } ]) \cap C^2((0, \NUU \MUU^{- \alpha } ])$ of~\eqref{IPE} satisfying the conditions~\eqref{fTT1:b1}.  
Moreover, \cite[Proposition~3.1]{SoupletW} implies that the solution $w_{\MUU} $ can be extended to $[0, \infty )$. This proves the first part of the statement.
For future reference, we note that by~\cite[formula~(2.1)]{SoupletW}, since $w_\MUU '( \NUU \MUU^{-\alpha } )=0$,
\begin{equation} \label{fLF1} 
4 s^\gamma e^{\frac {1} {4s}}  w _\MUU ' (s) = - \int  _{ \NUU \MUU^{-\alpha }  } ^s \sigma ^{\gamma -2} e^{\frac {1} {4 \sigma }} g ( w_\MUU (\sigma )) \,d\sigma  .
\end{equation} 
for all $s>0$. 

We now prove continuous dependence. 
The result of~\cite[Theorem~2.5]{SoupletW} cannot immediately be applied since the value of $T$ where $w_\MUU '(T)=0$ depends on $\MUU$. 
Fix $ \overline{\MUU} \ge 1$,  let $w _{  \overline{\MUU}  }$ be as above, and set
\begin{equation*} 
 \overline{R}  =  2 \overline{\MUU } , \quad  \overline{M} = 20  \overline{\MUU}  , \quad  \overline{T} = \NUU  \overline{\MUU }^{-\alpha }  .
\end{equation*} 
It follows from~\eqref{fLDC0}--\eqref{fLDC1}  that~\eqref{fLDC2:bz}-\eqref{fLDC2}  hold, with $R, M, T$ replaced by $ \overline{R} ,  \overline{M} ,  \overline{T} $. 
Let $\delta >0$ to be chosen sufficiently small, let $\MUU \ge 1$ satisfy $  |\MUU -  \overline{\MUU} | \le \delta $.
Let $w_\MUU$ be as above, so that $ |w_\MUU |\le 10\MUU$ on $[0, \NUU \MUU^{- \alpha }]$. 
In particular, if $\delta \le  \frac {1} {2}  \overline{\MUU} $, then 
\begin{equation*} 
\MUU \le  \frac {3} {2}  \overline{\MUU} < \overline{R}  , 
\end{equation*} 
so that 
\begin{equation}  \label{fCCDD2:b1} 
 |w_\MUU (s)| \le  15  \overline{\MUU} <  \overline{M}  ,\quad 0\le s\le  \NUU \MUU^{- \alpha } .  
\end{equation} 
We claim that if $\delta >0$ is sufficiently small, then 
\begin{equation}  \label{fCCDD2} 
 |w_\MUU (s)| \le 20  \overline{\MUU} =  \overline{M}  ,\quad 0\le s\le  \overline{T}.  
\end{equation} 
This is immediate if $\MUU \le  \overline{\MUU} $, since then $ \overline{T} \le \NUU \MUU^{-\alpha }$.
To see this in the case $\MUU >  \overline{\MUU} $, 
we define $\tau _\MUU \in ( \NUU \MUU^{- \alpha }, \infty ]$  by
\begin{equation}  \label{fCCDD11} 
\tau _\MUU = \sup \{ t>0; \,  |w_\MUU (s)| \le 20  \overline{\MUU}  \text{ on }[0,t]  \}.
\end{equation} 
It follows from~\eqref{fLF1} and~\eqref{fCCDD11} that
\begin{equation*} 
4 s^\gamma e^{\frac {1} {4s}}  |w _\MUU ' (s)| \le  \Bigl(  \int  _{ \NUU \MUU^{-\alpha }  } ^s \sigma ^{\gamma -2} e^{\frac {1} {4 \sigma }} d\sigma  \Bigr)\sup _{  |y| \le 20  \overline{\MUU}  }  |g(y)| ,
\end{equation*} 
 for all $s \in [  \NUU \MUU^{- \alpha }, \tau _\MUU)$ hence, using~\eqref{fCCDD2:b1}, 
\begin{equation}  \label{fCCDD13} 
 |w _\MUU (t)| \le 15  \overline{\MUU}  + \frac {1} {4}  \Bigl( \int _{ \NUU \MUU^{-\alpha }  } ^t s^{ -\gamma } e^{ - \frac {1} {4s}}  \int  _{ \NUU \MUU^{-\alpha }  } ^s \sigma ^{\gamma -2} e^{\frac {1} {4 \sigma }} d\sigma ds \Bigr)\sup _{  |y| \le 20  \overline{\MUU}  }  |g(y)| ,
\end{equation} 
 for all $t \in [  \NUU \MUU^{- \alpha }, \tau _\MUU)$. 
Since the right-hand side of~\eqref{fCCDD13} is increasing in $t$, to show that $\tau _\MUU \ge  \overline{T} $, we need only verify that 
\begin{equation*} 
 \frac {1} {4}  \Bigl( \int _{ \NUU \MUU^{-\alpha }  } ^{ \NUU  \overline{\MUU} ^{-\alpha }  } s^{ -\gamma } e^{ - \frac {1} {4s}}  \int  _{ \NUU \MUU^{-\alpha }  } ^s \sigma ^{\gamma -2} e^{\frac {1} {4 \sigma }} d\sigma ds \Bigr)\sup _{  |y| \le 20  \overline{\MUU}  }  |g(y)| \le 5  \overline{\MUU}, 
\end{equation*} 
 which is clearly true if $\delta >0$ is sufficiently small. Hence~\eqref{fCCDD2} holds. It now follows from~\eqref{fLF1} with $s= \overline{T} $, \eqref{fCCDD2:b1}  and~\eqref{fCCDD2} that 
\begin{equation*} 
4  \overline{T} ^\gamma e^{\frac {1} {4 \overline{T} }}  |w _\MUU ' (  \overline{T} )| \le  \Bigl|  \int  _{ \NUU \MUU^{-\alpha }  } ^{  \overline{T}  } \sigma ^{\gamma -2} e^{\frac {1} {4 \sigma }} d\sigma  \Bigr| \sup _{  |y| \le 20  \overline{\MUU}  }  |g(y)| ,
\end{equation*} 
and we deduce that there exists a constant $C$ such that 
 \begin{equation*} 
  |w _\MUU ' (  \overline{T} )| \le C  | \MUU -  \overline{\MUU} |
 \end{equation*} 
 if $\delta >0$ is sufficiently small. 
Therefore, for small $\delta >0$ we see that $ \| w _{ \MUU } \| _{ L^\infty (0,  \overline{T} ) }\le  \overline{M} $, $ | w_\MUU (0) |\le  \overline{R} $, and $  \overline{T}^2   | w_\MUU '(  \overline{T}  ) | \le  \overline{R}  $. 
The last statement of~\cite[Theorem~2.5]{SoupletW} shows that $w_\MUU \to w _{  \overline{\MUU}  }$ in $C^1 ([0,  \overline{T}]) $ as $\MUU \to  \overline{\MUU} $. 
Since equation~\eqref{IPE}  is not degenerate on $[  \overline{T} , T]$ for $T>  \overline{T}  $, continuous dependence in $C^1([0,T])$ follows easily.
 
Identity~\eqref{fENE3:b1} follows from elementary calculations.
 
 Moreover, \cite[Proposition~3.1~(i)]{SoupletW} implies that the solution $w_{\MUU} $ satisfies~\eqref{eEX1:1:b1}.

If $(N-2) \alpha \le 4$, then $\gamma >1$, see \eqref{gma}. It follows that $H_{\MUU} (s)$ defined by~\eqref{fENE1:b2} is nonincreasing for $s$ large, and bounded by~\eqref{eEX1:1:b1}, so it has a finite limit as $s\to \infty $. This proves Property~\eqref{eEX1:2:b2}. 
 
 The second statement of Property~\eqref{eEX1:3:b1} is a consequence of~\cite[Propositions~3.3 and~3.2]{SoupletW}. In particular, $G (w_\MUU (s))$ defined by~\eqref{fNot2}  has a limit as $s\to \infty $. Since $H_{\MUU} (s)$ also has a limit, we see that $s^2 w_{\MUU} '(s)^2$ must have a limit, which is necessarily $0$. This proves the first statement of Property~\eqref{eEX1:3:b1}.

We now turn to the proof of Property~\eqref{eEX1:2:6}.
If $\beta >0$ and $w_\MUU (s)\to \pm \beta ^{\frac {1} {\alpha }}$ then the result is clearly true. 
Otherwise, if $\beta >0$ and $w_\MUU (s)\not\to \pm \beta ^{\frac {1} {\alpha }}$, or if
$\beta \le 0$, then $w_\MUU (s)\to 0$, by Property~\eqref{eEX1:3:b1}. 
We first consider the case $N \ge 3$, so that $\lambda _1 > \lambda _2$ by~\eqref{fPR1}. 
Given $\sigma \in \R$,  let
\begin{equation} \label{fBN5a0} 
z(s)=s^\sigma w_\MUU (s).
\end{equation} 
It follows from~\eqref{IPE} and~\eqref{fNZE3} that
\begin{equation} \label{fBN5a} 
4s^2z''+ 4 (1 + \lambda _1 + \lambda _2 -2\sigma )sz'-z'+ 4 (\sigma -\lambda _1) (\sigma -\lambda _2) z+\frac {\sigma } {s}z+ \frac {1} {s^{ \alpha \sigma }} |z|^\alpha z=0.
\end{equation}  
If we fix $\sigma >0$, $\sigma \in ( \lambda _2, \lambda _1)$, then $ 4 (\sigma -\lambda _1) (\sigma -\lambda _2) < 0$.
Since $s^{-\alpha \sigma }|z|^\alpha=|w|^\alpha \to 0$ as $s\to \infty$, we deduce from~\eqref{fBN5a} that if $s_0>0$ is sufficiently large, then for $s\ge s_0$, $z(s) z''(s)>0$ whenever $z' (s) =0$.
In other words, at any point $s\ge s_0$ where $z'(s)=0$, $z^2$ has a local minimum. Thus $z$ can vanish at most once for $s\ge s_0$. 

If $N=1$, so that $\lambda _1< \lambda _2$ by~\eqref{fPR1}, the above argument works with $\sigma \in ( \lambda _1, \lambda _2)$.

The case $N=2$, where $\lambda _1 = \lambda _2$, requires a more delicate argument. 
We let $\sigma =\lambda _1$, so that equation~\eqref{fBN5a} becomes
\begin{equation} \label{fBN5a1} 
4s^2z''+ 4 sz'-z'  + \frac {\lambda _1 } {s}z+ \frac {1} {s } |z|^\alpha z=0.
\end{equation} 
Multiplying~\eqref{fBN5a1} by $z' $ yields the energy identity
\begin{equation*} 
 \Bigl( 2 s^2 z'\null ^2 + \frac {\lambda _1} {2s} z^2 + \frac {1} {(\alpha +2) s}  |z|^{\alpha +2} \Bigr) '= z' \null ^2 -  \frac {\lambda _1} {2s^2} z^2- \frac {1} {(\alpha +2) s^2}  |z|^{\alpha +2} \le z' \null ^2
\end{equation*} 
so that
\begin{equation*} 
 \Bigl( 2 s^2 z'\null ^2 + \frac {\lambda _1} {2s} z^2 + \frac {1} {(\alpha +2) s}  |z|^{\alpha +2} \Bigr) ' \le  \frac {1} {2s^2} \Bigl( 2 s^2 z'\null ^2 + \frac {\lambda _1} {2s} z^2 + \frac {1} {(\alpha +2) s}  |z|^{\alpha +2} \Bigr) .
\end{equation*} 
Integrating the above inequality yields
\begin{equation} \label{fBN5a3} 
2 s^2 z' (s) ^2 \le e^{\frac {1} {2}- \frac {1} {2s}}  \Bigl( 2  z'(1) ^2 + \frac {\lambda _1} {2} z(1) ^2 + \frac {1} {\alpha +2}  |z (1) |^{\alpha +2} \Bigr) 
\end{equation} 
for $s\ge 1$. 
Therefore, $ | z' (s)| \le C (1+s)^{-1}$, so that
\begin{equation} \label{fBN5a4} 
 |z (s)| \le C \log (2 +s)
\end{equation} 
for $s\ge 1$. 
We now  define for $s> 1$
\begin{equation} \label{fBN5a5} 
t^{\frac {1} {2}} v (t)=z (s), \quad t = \log s .
\end{equation} 
It follows in particular that
\begin{equation} \label{fbw1} 
\frac {dt} {ds}= \frac {1} {s}.
\end{equation} 
We deduce from~\eqref{fBN5a5} and~\eqref{fbw1} that
\begin{equation} \label{fbw2} 
z'(s)= \frac {1} {s}  \Bigl( t^{\frac {1} {2}} v'(t) + \frac {1} {2} t^{- \frac {1} {2}} v(t)  \Bigr)
\end{equation} 
so that
\begin{equation} \label{fbw3} 
s z'(s)=  t^{\frac {1} {2}} v'(t) + \frac {1} {2} t^{- \frac {1} {2}} v(t)  .
\end{equation} 
Next, differentiating~\eqref{fbw2} with respect to $s$ and applying again~\eqref{fbw1}, we obtain 
\begin{equation} \label{fbw4} 
s^2 z''(s) =  -  \Bigl( t^{\frac {1} {2}} v'(t) + \frac {1} {2} t^{- \frac {1} {2}} v(t)  \Bigr) + \Bigl( t^{\frac {1} {2}} v''(t) +   t^{- \frac {1} {2}} v '(t) - \frac {1} {4}  t^{- \frac {3} {2}} v (t)  \Bigr) .
\end{equation} 
It follows from~\eqref{fBN5a1}, \eqref{fBN5a5}, \eqref{fbw3}  and~\eqref{fbw4}  that
\begin{equation} \label{eNZE3:3}
4v ''+\Bigl ( \frac { 4 } {t}-e^{-t}\Bigr ) v '-\Bigl (\frac { 1 } {t^2}+e^{-t}\Bigl (\frac {1} {2 t}-\frac {1} {\alpha } - ( t^{\frac {1} {2}} |v |)^\alpha  \Bigr)\Bigr ) v  =0.
\end{equation}
Note that by~\eqref{fBN5a4} and~\eqref{fBN5a5},
\begin{equation*} 
 t^{\frac {1} {2}} |v | =  |z (s)| \le  C (\log (2 +s)) \le C t 
\end{equation*} 
for $t$ large, so that
\begin{equation} \label{fTPA19} 
t^2 e^{-t}\Bigl (\frac {1} {2 t}-\frac {1} {\alpha } - ( t^{\frac {1} {2}} |v |)^\alpha  \Bigr) \goto  _{ t\to \infty  } 0 .
\end{equation} 
It follows from~\eqref{eNZE3:3} and~\eqref{fTPA19}  that for $t$ large, if $v' $ vanishes, then $v '' $ has the sign of $v$. 
Arguing as  in the case $N\ge 3$, we conclude that $v$ has a finite number of zeros.

We finally prove Property~\eqref{eEX1:2:7}, so we suppose $w_\MUU (s) \to 0$ as $s\to \infty $. By Property~\eqref{eEX1:2:6}, $w_\MUU (s) \not = 0$ for $s$ large, and we deduce from~\eqref{IPE} that if $w_\MUU ' (s)= 0$, then 
\begin{equation*} 
4s^2 w_\MUU '' (s)= (\beta - |w_\MUU (s)|^\alpha ) w_\MUU (s) .
\end{equation*} 
It easily follows that if $s$ is large, then either $w_\MUU '' (s) $ has the sign of $w_\MUU (s)$ (if $\beta >0$), or else $w_\MUU '' (s) >0$ has the opposite sign of $w_\MUU (s)$ (if $\beta <0$ or if $\beta =0$) whenever $w_\MUU ' (s)=0$. This shows that $w_\MUU ' $ cannot vanish for $s$ large.  Since $w_\MUU (s) \to 0$ as $s\to \infty $, we conclude that $w_\MUU w_\MUU ' <0$ for $s$ large. 
 
\section{Arbitrarily many zeros} \label{sAB} 
This section, in its entirety, constitutes the proof of Proposition~\ref{fFIn1}. 
We set 
\begin{equation}
\label{frd}
b = \frac {1} {4\gamma +8}.
\end{equation}
It will be shown that there exists a function
$\iota :(1, \infty) \to (0,\infty)$  such that
$\iota (\MUU) \to 0$ as $\MUU \to \infty$ with the property that
if $I \subset (0, b )$ is an interval on which $w_\MUU(s) \neq 0$, then $|I| \le \iota(\MUU)$. This implies the proposition. Indeed, if 
$\iota(\MUU) < \frac {T} {m}$, where $T < b$, then $w_\MUU$ has at least $m$ zeros on $[0,T]$.

Thus we fix an interval $I \subset (0, b )$ for which $w_\MUU(s) > 0$ for all $s \in I$.  (The case $w_\MUU(s) < 0$ can be handled
analogously.) We need to estimate $|I|$ as a function of
$\MUU$, and it suffices to do so for $\MUU > 1$ sufficiently large.

The crucial observation is that $H_\MUU $ given by~\eqref{fENE1:b2} is increasing on $[0,b]$, by \eqref{fENE3:b1}, so that 
\begin{equation} \label{fTPA01} 
 2s^2w_\MUU'(s)^2 + G(w_\MUU(s)) = H _\MUU (s) > H _\MUU (0) = G(  \MUU   ) 
\end{equation} 
 for $ s\in (0, b ] $.  The first largeness condition we impose
 on $\MUU$ is that
 \begin{equation}
 \label{frd1}
 \MUU > 2\big(2(\alpha + 2)|\beta|\big)^{\frac{1}{\alpha}}.
 \end{equation}
This condition guarantees, in particular, that there is at most one value of $s \in I$
where $w_\MUU'(s) = 0$.  To see this, suppose $s \in I$ and $w_\MUU'(s) = 0$. 
It follows by \eqref{fTPA01} and \eqref{frd1} that
\begin{equation*}
G(w_\MUU(s))  > G(  \MUU   ) > 0,
\end{equation*}
and so it must be that
$w_\MUU(s) > \MUU$.   This implies, again by \eqref{frd1},  that $g(w_\MUU(s)) > 0$,
which in turn implies, by equation \eqref{IPE}, that $w_\MUU''(s) < 0$.
In other words, at any point in $I$ where $w_\MUU'(s) = 0$ the solution 
$w_\MUU$ must have
a local maximum. Hence there is at most one such point in $I$.

It follows that the interval $I$ can be partitioned into
four pairwise disjoint subintervals, which depend on $\MUU$,  
 \begin{equation}
 \label{frd3}
 I = I_\MUU^1 \cup I_\MUU^2 \cup I_\MUU^3 \cup I_\MUU^4,
 \end{equation}
given by
\begin{gather*} 
I_\MUU^1 = \{s \in I: w_\MUU(s) \le \frac{\MUU}{2}, w_\MUU'(s) > 0\}, \quad 
I_\MUU^2 = \{s \in I: w_\MUU(s) > \frac{\MUU}{2}, w_\MUU'(s) > 0\}, \\
I_\MUU^3 = \{s \in I: w_\MUU(s) > \frac{\MUU}{2}, w_\MUU'(s) \le 0\}, \quad 
I_\MUU^4 = \{s \in I: w_\MUU(s) \le \frac{\MUU}{2}, w_\MUU'(s) < 0\} .
\end{gather*} 
Note that one or more of these intervals might be empty or contain just
 a single point.
We now proceed to estimate the lengths of these four intervals.

Consider first the interval $I_\MUU^1$.
If $s \in I_\MUU^1$, then $G(w_\MUU(s))\le G({\MUU} / {2})$ by~\eqref{fNot2} and \eqref{frd1}. Thus 
\begin{equation*}
w_\MUU'(s)^2 \ge \frac{1}{2s^2}(G(\MUU) - G(w_\MUU(s)))
\ge \frac{1}{2b^2}(G(\MUU) - G({\MUU} / {2}))
\ge \frac{\MUU^{\alpha+2}}{8b^2(\alpha + 2)}.
\end{equation*}
(The last inequality requires an explicit calculation using again~\eqref{fNot2} and \eqref{frd1}.)
Consequently, if $s < t$ and $s,t \in I_\MUU^1$, then
\begin{equation*}
|w_\MUU(t) - w_\MUU(s)|
\ge \frac {t-s } {2b\sqrt{2(\alpha + 2)}} \MUU^{1 + \frac {\alpha} {2}}.
\end{equation*}
Since $|w_\MUU(t) - w_\MUU(s)| \le \frac{\MUU}{2}$ if
$s,t \in I_\MUU^1$, this implies
\begin{equation}
\label{frd4}
|I_\MUU^1| \le b\sqrt{2(\alpha + 2)}\MUU^{- \frac{\alpha}{2}}.
\end{equation}
In exactly the same way, we also have
\begin{equation}
\label{frd5}
|I_\MUU^4| \le b\sqrt{2(\alpha + 2)}\MUU^{- \frac{\alpha}{2}}.
\end{equation}

We next turn our attention to $I_\MUU^2$.  To estimate its length,
we partition this interval into further subintervals. 
Differentiating~\eqref{IPE} we get
\begin{equation*} 
4s^2 w _ \MUU  '''-(1-(4\gamma +8)s)w _ \MUU  ''+(4\gamma  -\beta +(\alpha +1)|w _ \MUU  |^\alpha )w _ \MUU  '=0.
\end{equation*} 
Therefore, if $w _ \MUU  ''(s)=0$ for
some $s \in I_\MUU^2$,  we have from~\eqref{frd1} that $w _ \MUU  '''(s)<0$,
and so $w_\MUU '' $ has a local maximum at that point.
It follows that $w _ \MUU''$ can have at most one zero in $I_\MUU^2$.
We may therefore define the intervals
\begin{equation*} 
J_\MUU^1 = \{s \in I_\MUU^2: w_\MUU''(s) \ge 0\}, \quad J_\MUU^2 = \{s \in I_\MUU^2: w_\MUU''(s) < 0\}
\end{equation*} 
so that $I_\MUU^2 = J_\MUU^1 \cup  J_\MUU^2$.

If $s \in J_\MUU^1$, then
\begin{equation*}
w_\MUU'(s) \ge \frac{g(w_\MUU(s))}{1 - 4\gamma s}
\ge \frac{w_\MUU(s)^{\alpha+1}}{2(1 - 4\gamma s)}
\ge \frac{1}{2}w_\MUU(s)^{\alpha+1},
\end{equation*}
where again we use~\eqref{IPE} and~\eqref{frd1}.  In other words,
$(w_\MUU(s)^{-\alpha})' \le -\frac{\alpha}{2}$.
Integrating, we see that if $s < t$ and $s,t \in J_\MUU^1$,
then
\begin{equation*}
t - s \le \frac{2}{\alpha}w_\MUU(s)^{-\alpha}
\le \frac{2^{\alpha + 1}}{\alpha}\MUU^{-\alpha}.
\end{equation*}
This shows that
\begin{equation}
\label{frd6}
|J_\MUU^1| \le \frac{2^{\alpha + 1}}{\alpha}\MUU^{-\alpha}.
\end{equation}

If $s \in J_\MUU^2$, then
\begin{equation*}
(w_\MUU'(s) - w_\MUU(s)^{\frac{\alpha}{2}+1})'
= w_\MUU''(s) - ({\frac{\alpha}{2}+1}) w_\MUU(s)^{\frac{\alpha}{2}} w_\MUU'(s) < 0.
\end{equation*}
Therefore, $J_\MUU^2$ is the union of two intervals,
$J_\MUU^{2+}$ and $J_\MUU^{2-}$, where, respectively,
$w_\MUU'(s) \ge w_\MUU(s)^{\frac{\alpha}{2}+1}$,
and $w_\MUU'(s) < w_\MUU(s)^{\frac{\alpha}{2}+1}$.
On $J_\MUU^{2+}$, we have that
$(w_\MUU(s)^{-\frac{\alpha}{2}})' \le -\frac{\alpha}{2}$,
which integrates to yield, if $s < t$ and $s,t \in J_\MUU^{2+}$,
\begin{equation*}
t - s \le \frac{2}{\alpha}w_\MUU(s)^{-\frac{\alpha}{2}}
\le \frac{2^{1+\frac{\alpha}{2}}}{\alpha}\MUU^{-\frac{\alpha}{2}}.\end{equation*}
This implies that
\begin{equation}
\label{frd7}
|J_\MUU^{2+}| \le  \frac{2^{1+\frac{\alpha}{2}}}{\alpha}\MUU^{-\frac{\alpha}{2}}.
\end{equation}
On the other hand, if $s \in J_\MUU^{2-}$, then
\begin{equation*} 
\begin{split} 
4b^2 w_\MUU '' (s) & \le 4s^2 w_\MUU '' (s) = (1-4\gamma s)w _ \MUU  '(s)-g(w_ \MUU  (s)) \\
& \le w _ \MUU  '(s)-g(w_ \MUU  (s))   \le  w _\MUU (s)^{1+ \frac {\alpha } {2}} - \frac {1} {2} w_ \MUU  (s)^{\alpha +1} \\ 
& \le - \frac {1} {4} w_ \MUU  (s)^{\alpha +1},
\end{split} 
\end{equation*}
where we need to impose an additional largeness condition on $\MUU$, { i.e.}
\begin{equation}
\label{frd8}
\MUU \ge 2^{1 + \frac{4}{\alpha}}.
\end{equation}
Consequently, if $s<t$ and $s, t \in  J_\MUU^{2-}$, then
\begin{equation*} 
  w_\MUU ' (s)  \ge    w_\MUU ' (s) - w_\MUU ' (t) \ge  \frac {1} {16 b^2} \int _s^t w_ \MUU  (\sigma )^{\alpha +1} d\sigma  \ge  \frac {t-s} {16 b^2} 
   w_ \MUU  (s )^{\alpha +1} .
\end{equation*} 
Since $w_\MUU ' (s) \le  w _\MUU (s)^{1+ \frac {\alpha } {2}} $ and $w_\MUU (s) \ge \frac {\MUU} {2}$,  we deduce that 
\begin{equation*} 
1 \ge   \frac {t-s} {16 b^2}  w_ \MUU  (s )^{\frac {\alpha } {2}} \ge   \frac {t-s} {2^{4+ \frac {\alpha } {2}}  b^2}   \MUU ^{\frac {\alpha } {2}} .
\end{equation*} 
Therefore,
\begin{equation}
\label{frd9}
 |J_\MUU^{2-}| \le    2^{ 4 + \frac {\alpha } {2}}b^2 \MUU ^{ - \frac {\alpha } {2} }.
\end{equation}
For future reference, we recall that
\begin{equation}
\label{frd10}
I_\MUU^2 = J_\MUU^1 \cup J_\MUU^{2+}\cup J_\MUU^{2-},
\end{equation}
and so the length of $I_\MUU^2$ is estimated by \eqref{frd6}, 
\eqref{frd7} and \eqref{frd9}.

Finally, we estimate the length of the interval $I_\MUU^3$.
Observe first that since $g(w_\MUU(s)) > 0$ on $I_\MUU^3$
(by \eqref{frd1}) and $w_\MUU'(s) \le 0$ on $I_\MUU^3$,
it follows from \eqref{IPE} that $w_\MUU''(s) < 0$ on $I_\MUU^3$.
Hence, for $s\in I_\MUU^3$, 
\begin{equation*} 
4b ^2 w _\MUU '' (s) +g(w_\MUU (s) )\le 4s^2  w _\MUU '' (s) +g(w_\MUU (s) ) = (1 - 4\gamma s) w_\MUU ' (s)  \le 0
\end{equation*} 
so that
\begin{equation*} 
(2 b ^2 w _\MUU '(s)^2+ G( w _\MUU (s)))'=(4b^2 w _\MUU ''(s)+ g(w _\MUU (s)))w _\MUU '(s) \ge 0.
\end{equation*}
Therefore, if $s < t$ and $s,t \in I_\MUU^3$,
\begin{equation} \label{fTPA08:3} 
2b ^2 w _\MUU '(t)^2 + G(w _\MUU (t)) \ge 2b^2w _\MUU' (s )^2 + 
G ( w _\MUU (s)) \ge G( w _\MUU (s )) .
\end{equation}
Recall that $w_\MUU '' <0$ and $w_\MUU ' \le 0$ on $I_\MUU^3$, so that $w_\MUU (s) > w_\MUU (t)$.
Since $G(w _\MUU (t)) > 0$ by \eqref{frd1}, we deduce that  $G(w _\MUU (s)) > G(w _\MUU (t)) >0$, and it follows
from \eqref{fTPA08:3} that
\begin{equation*} 
 \frac {-\sqrt{2} b w _\MUU '(t)} {(G(w _\MUU (s))-G(w _\MUU (t)))^{\frac {1} {2}}} \ge  1.
\end{equation*}
Integrating, we obtain
\begin{equation} \label{fNP10} 
\begin{split} 
t - s &  \le \int_s^t  \frac {-\sqrt{2} b w _\MUU '(\tau)} {(G(w _\MUU (s))-G(w _\MUU (\tau)))^{\frac {1} {2}}}d\tau\\
& =  \int_{\frac {w _\MUU (t)} {w _\MUU (s)}}^{1} \frac {\sqrt 2 b w _\MUU (s)} {(G(w _\MUU (s))-  G( z w _\MUU (s) ))^{\frac {1} {2}}}\, dz,
\end{split} 
\end{equation}  
where we have set $z = \frac{w _\MUU (\tau)}{w _\MUU (s)}$, for $s \le \tau \le t$.
Furthermore, once again using \eqref{frd1}, 
\begin{equation*} 
\begin{split}
G(w _\MUU (s)) -  G(zw _\MUU ( s) ) & = \frac {w _\MUU (s) ^{\alpha +2} } {\alpha +2} (1-z^{\alpha +2})  - \frac {\beta w _\MUU (s) ^2 } {2} (1-z^2) \\
&\ge  \Bigl( \frac {w _\MUU ( s) ^{\alpha +2}} {\alpha +2}- \frac { |\beta| w _\MUU ( s) ^2} {2}\Bigr) (1-z^{\alpha +2})  \\ & \ge \frac {w _\MUU (s) ^{\alpha +2}} {2(\alpha +2)} (1-z^{\alpha +2}).
\end{split}
\end{equation*} 
Putting this into \eqref{fNP10}, we obtain that
\begin{equation*} 
t - s   \le  \frac {2 b \sqrt{\alpha +2} } { w _\MUU(s) ^{\frac {\alpha } {2}} } \int_0 ^{1} \frac { dz} { (1-z^{\alpha +2}) ^{\frac {1} {2}}}  \le \frac {2^{1+\frac {\alpha } {2}} b \sqrt{\alpha +2} } {  \MUU ^{\frac {\alpha } {2}} } \int_0 ^{1} \frac { dz} { (1-z^{\alpha +2}) ^{\frac {1} {2}}},
\end{equation*} 
so that
\begin{equation}
\label{frd11} 
|I_\MUU^3| \le 2^{1+\frac {\alpha } {2}} b \sqrt{\alpha +2}  \Bigl(   \int_0 ^{1} \frac { dz} { (1-z^{\alpha +2}) ^{\frac {1} {2}}} \Bigr) \, \MUU ^{-\frac {\alpha } {2}}.
\end{equation}

The proof is now complete.  
Indeed, by~\eqref{frd3} and~\eqref{frd10}, the function $\iota(\MUU)$,
for $\MUU$ satisfying \eqref{frd1} and \eqref{frd8}, can
be taken as the sum of the right hand sides of the estimates
\eqref{frd4}, \eqref{frd5}, \eqref{frd6}, \eqref{frd7}, \eqref{frd9} and \eqref{frd11}, where $b$ is given by~\eqref{frd}

\section{Asymptotic behavior of the solutions of~\eqref{IPE}} \label{sec:auxiliary} 

In this section, we give the proof of Proposition~\ref{eAB1}, which concerns the asymptotic behavior as $s\to \infty $ of the solutions $( w_\MUU ) _{ \MUU \ge 1 }$ of~\eqref{IPE} given by Proposition~$\ref{eEX1:b1}$. 
This behavior must be studied separately in each of the four different cases

\begin{enumerate}[{\rm (1)}] 
\item \label{ic1}  $N\ge 3$ and $\beta >0$;

\item \label{ic2}  $N\ge 3$ and $\beta <0$;

\item \label{ic3}  $N\ge 3$ and   $\beta =0$;

\item \label{ic4}  $N=2$.
\end{enumerate} 
Although this section is somewhat technical, the results are not surprising and the methods used are standard techniques. We note that some of the results in this section will also be used in the proof of Proposition~\ref{eEX2}.

Before we begin the detailed analysis, we observe that in all cases, if  $\Ltu( \MUU )$ given by \eqref{eAB4} exists and is finite., then $f(r)=r^{-\frac {2} {\alpha }}w_\MUU(r^{-2})$ solves~\eqref{fpr3}-\eqref{fpr3:1} with $a =\Ltu(\MUU)$. Since $w_\MUU\not \equiv  0$, we deduce that $\Ltu(\MUU)\ne 0$. 

\subsection{Proof of Proposition~$\ref{eAB1}$ in the case $N\ge 3$ and $\beta > 0$} \label{sec:beta_positive}  

\begin{lemma} \label{eEX1:c1}  
Assume $N\ge 3$ and $\beta >0$ and let $\MUU \ge \MUUz$. If $w_\MUU (s) \to 0$ as $s\to \infty $, then the limit~\eqref{eAB4}  exists and is finite.
\end{lemma} 

\begin{proof} 
 It follows from Proposition~\ref{eEX1:b1}~\eqref{eEX1:2:6} that  $w_{\MUU} $ has a constant sign for $s$ large.  The result is now a consequence of~\cite[Proposition~3.5]{SoupletW}. (Note that this last proposition assumes that $w _{\MUU}  $ has a constant sign on $(0,\infty )$, but the argument uses only the fact that $w _{\MUU}  $ has a constant sign for $s$ large.)
\end{proof} 

Proposition~\ref{eAB1}  follows from Proposition~\ref{eEX1:b1}~\eqref{eEX1:3:b1} and Lemma~\ref{eEX1:c1}.

\subsection{Proof of Proposition~$\ref{eAB1}$ in the case $N\ge 3$ and $\beta < 0$} \label{sec:beta_negative} 

Throughout this subsection, we assume that $N\ge 3$ and $\beta <0$, so that $\lambda _1 > \lambda _2 >0$ with the notation~\eqref{fPR1}. 
We define
\begin{equation} \label{fBNb2} 
Lw=4s^2w''+4\gamma sw'-\beta w
\end{equation}  
so that $Lw_\MUU=\Phimu _\MUU$, where 
\begin{equation} \label{fTPA55} 
\Phimu _ \MUU =w_ \MUU '-|w_ \MUU |^\alpha w_ \MUU .
\end{equation} 
Let $\phi _1$ and $\phi _2$ be defined by~\eqref{fPR2} and~\eqref{fPR3}, respectively. 
A simple calculation shows that  $L\phi _1=L\phi _2=0$.  It then straightforward to check that $w_\MUU $ satisfies the following variation of the parameter formula
\begin{equation} \label{fBN3} 
\begin{split} 
w_\MUU = & \left (\frac {\lambda _2 w_\MUU (1)+w_\MUU '(1)} {\lambda _2-\lambda _1}-\frac {1} {2N-4}\int_1^s \tau^{\lambda _1-1} \Phimu _\MUU  \, d\tau \right )\phi _1 \\
& + \left (\frac {\lambda _1 w_\MUU (1)+ w_\MUU '(1)} {\lambda _1-\lambda _2}+\frac {1} {2N-4}\int_1^s \tau^{\lambda _2-1} \Phimu _\MUU \, d\tau \right )\phi _2
\end{split} 
\end{equation}  
for all $s>0$.

\begin{lemma} \label{eBN5} 
Given any $\MUU \ge 1$, the limit~\eqref{eAB2} exists and is finite, the map $s\mapsto s^{\lambda _2-1} \Phimu_ \MUU  (s)$ belongs to $L^1(1,\infty)$, and 
\begin{equation} \label{eBN6}
\Lt( \MUU )=\frac {\lambda _1 w_ \MUU (1)+w_ \MUU '(1)} {\lambda _1-\lambda _2}+\frac {1} {2N-4}\int_1^\infty \tau^{\lambda _2-1}\Phimu_ \MUU (\tau)\, d\tau .
\end{equation} 
\end{lemma}

\begin{proof} 
We first prove that 
\begin{equation} \label{fBN10}
\sup_{s \ge 1} \, [s^{\lambda _2} | w_ \MUU  (s)|+ s^{\lambda _2+1}|w_ \MUU  '(s)| ] ( \log s )^{-1} < \infty .
\end{equation} 
Indeed, $z(s)=s^{\lambda _2 }w_ \MUU   (s)$ satisfies~\eqref{fBN5a} with $\sigma =\lambda _2$, i.e. (using~\eqref{fPR1}) 
\begin{equation*} 
4s^2z''+ 2 N sz'-z' +\frac {\lambda _2 } {s}z+ \frac {1} {s^{ \alpha \lambda _2 }} |z|^\alpha z=0.
\end{equation*}  
Therefore, if
\begin{equation} \label{fBN6} 
 \widetilde{H}  (s)=2s^2|z'|^2+\frac {\lambda _2 } {2s} z^2  +\frac {1} {(\alpha +2) s^{\alpha \lambda _2}}|z|^{\alpha +2}
\end{equation}
then
\begin{equation}  \label{fBN6b} 
 \widetilde{H}  '(s)=-\frac {\alpha \lambda _2 |z|^{\alpha +2}} {(\alpha +2)s^{\alpha \lambda _2 +1}} - \frac {\lambda _2 |z|^2} {2s^2}+(1  - 2 (N - 2 ) s ) | z ' | ^ 2
\end{equation} 
It follows that the right hand side of~\eqref{fBN6b} is negative for $ s>s_0=\frac {1} { 2 (N - 2 ) }$; and so  $  \widetilde{H}  (s) $ is bounded. 
In particular, $ | z ' (s)| \le \frac {C} {s}$, so that $ | z (s)| \le C \log s$. Estimate~\eqref{fBN10} easily follows.

Next, we deduce from~\eqref{fBN10}  and the property $\lambda _2\alpha <1$ that
\begin{equation}  \label{eBN7} 
|\Phimu _\MUU (s)|\le |w_ \MUU '(s)|+|w_ \MUU (s)|^{\alpha +1}\le Cs^{-\lambda_2(\alpha +1)} (\log s)^{\alpha +1} 
\end{equation} 
for $s$ large.
Therefore, $s^{\lambda _2-1} \Phimu_ \MUU  (s)\in L^1(1,\infty)$ and 
\begin{equation}  \label{eBN7u} 
\begin{split} 
\int_1^s \tau^{\lambda _1-1} |\Phimu _ \MUU (\tau)|\, d\tau & \le \int_1^s \tau^{\lambda _1 - \lambda _2 (\alpha +1) -1} (\log \tau )^{\alpha +1}   \, d\tau \\ &\le C\max \{s^{\lambda _1-\lambda_2(\alpha +1)},1 \} (\log s)^{\alpha +1} .
\end{split} 
\end{equation} 
Hence, from~\eqref{fBN3} we see that $\Lt( \MUU )$ is well defined and that~\eqref{eBN6} holds.   
\end{proof} 

\begin{lemma} \label{eBN8} 
If  $\Lt( \MUU )=0 $, then the limit~\eqref{eAB4}  exists and is finite.
\end{lemma}
\begin{proof} 
If $\Lt( \MUU )=0 $ then~\eqref{fBN3} and~\eqref{eBN6} yield
\begin{equation} \label{eBN9} 
\begin{split} 
w_ \MUU  = & \left (\frac {\lambda _2 w_ \MUU  (1)+w_ \MUU  '(1)} {\lambda _2-\lambda _1}-\frac {1} {2N-4}\int_1^s \tau^{\lambda _1-1}\Phimu _ \MUU  (\tau)\, d\tau \right )\phi _1 \\
& - \left (\frac {1} {2N-4}\int_s^\infty \tau^{\lambda _2-1}\Phimu _ \MUU (\tau)\, d\tau \right )\phi _2 
\end{split} 
\end{equation}  
and
\begin{equation} \label{eBN9:b1} 
\begin{split} 
w_ \MUU  '= & \left (\frac {\lambda _2 w_ \MUU  (1)+w_ \MUU  '(1)} {\lambda _2-\lambda _1}-\frac {1} {2N-4}\int_1^s \tau^{\lambda _1-1}\Phimu_ \MUU   (\tau)\, d\tau \right )\phi '_1 \\
& - \left (\frac {1} {2N-4}\int_s^\infty \tau^{\lambda _2-1}\Phimu_ \MUU  (\tau)\, d\tau \right )\phi ' _2 .\end{split} 
\end{equation} 
It follows that 
\begin{equation*} 
\sup_{s >0} \, [ (1+s) ^{\lambda _2} |w_ \MUU  (s)|+  (1+s) ^{\lambda _2 +1}|w_ \MUU  '(s)| ] <\infty .
\end{equation*} 
Indeed, $(1+s) ^{\lambda _2} |w_ \MUU  (s)|$ is bounded by Lemma~\ref{eBN5}. Moreover, since $s^{\lambda _2-1}  |\Phimu_ \MUU  (s)|$ is integrable at infinity, again by Lemma~\ref{eBN5}, we deduce from formula~\eqref{eBN9:b1} that  $(1+s) ^{\lambda _2 +1}|w_ \MUU  '(s)| $ is bounded.
We now set 
\begin{equation} \label{eBN4} 
R= \Bigl\{ \rho >0,\   \sup_{s >0}\ \{ (1+s) ^\rho |w_ \MUU  (s)|+  (1+s) ^{\rho +1}|w_ \MUU  '(s)|\}<\infty \Bigr\}
\end{equation} 
so that $R$ is an interval and $\lambda _2\in R$. We claim that $\lambda _1\in R$.  We prove this by contradiction, so we assume $\lambda _1 \not \in R $. Therefore may choose $\rho >0$ such that
\begin{equation} \label{fTPA77} 
\rho \in R, \quad \rho \ge \lambda _2, \quad \rho < \lambda _1, \quad (\alpha +1) \rho > \sup R, \quad (\alpha +1) \rho \not = \lambda _1.
\end{equation} 
Note that, since $\rho \in R$, we have 
\begin{equation} \label{fTPA99} 
 |\Phimu _\MUU (s) | \le C ( s ^{-\rho -1} + s ^{ -(\alpha +1) \rho } ) 
\end{equation} 
 for $s\ge 1$.  Moreover, $\alpha \rho = \frac {\rho } {\lambda _1}< 1$, so that $  | \Phimu _\MUU (s) | \le C  s ^{ -(\alpha +1) \rho }  $; and so
\begin{equation}\label{eBN10}  
 \int_1^s \tau^{\lambda _1-1}|\Phimu _ \MUU (\tau) | \, d\tau \le C \theta (s)
\end{equation}  
where 
\begin{equation} \label{eBN10b}  
\theta (s)=   
\begin{cases} 
s^{\lambda _1-(\alpha +1)\rho} & (\alpha +1)\rho< \lambda _1 \\
\log (1+s) & (\alpha +1)\rho = \lambda _1 \\
1 &  (\alpha +1)\rho > \lambda _1
\end{cases} 
\end{equation} 
and  
\begin{equation}\label{eBN11} 
 \int_s^\infty \tau^{\lambda _2-1}|\Phimu _ \MUU  (\tau) | \, d\tau \le Cs^{\lambda _2-(\alpha +1)\rho} .
\end{equation}  
Using~\eqref{eBN10}, \eqref{eBN11} and, respectively, \eqref{eBN9} and~\eqref{eBN9:b1} we obtain 
\begin{equation} \label{eBN11b} 
|w_ \MUU  (s)| + s |w_ \MUU  '(s)| \le C \max\{ s^{- \lambda _1} \theta (s), s^{-(\alpha +1) \rho } \} .
\end{equation} 
Since $(\alpha +1) \rho \not \in R$ by~\eqref{fTPA77}, we deduce from~\eqref{eBN11b}  that 
\begin{equation*} 
\max\{ s^{-\lambda _1}  \theta (s) , s^{-(\alpha +1) \rho } \} = s^{-\lambda _1}  \theta (s) 
\end{equation*} 
so that
\begin{equation} \label{fTPA78} 
|w_ \MUU  (s)| + s |w_ \MUU  '(s)| \le  s^{-\lambda _1}  \theta (s) .
\end{equation}  
Since $\lambda _1\not \in R$ by assumption and $(\alpha +1) \rho \not \in R$ (again by~\eqref{fTPA77}), it follows from~\eqref{fTPA78} and~\eqref{eBN10b} that $\theta (s) =\log (1+s)$. This means that $\lambda _1= (\alpha +1) \rho $, which contradicts the last condition in~\eqref{fTPA77}. 
This contradiction establishes that $\lambda _1 \in R$.

To show that  $ \Ltu ( \MUU)$ is well defined, we let $\rho =\lambda _1$ in~\eqref{fTPA99}, so that $|\Phimu _\MUU (s)|\le Cs^{-(\lambda _1+1)}$. Thus, $s^{\lambda _1-1}\Phimu _ \MUU  (s)\in L^1(1,\infty)$. Moreover, 
\begin{equation*} 
\int_s^\infty \tau^{\lambda _2-1}\Phimu _ \MUU (\tau)\, d\tau \le C \int_s^\infty \tau^{- (\lambda _1- \lambda _2)-2}\, d\tau \le C s^{-( \lambda _1- \lambda _2)-1}.
\end{equation*} 
Thus we deduce from~\eqref{eBN9} that  the limit $ \Ltu ( \MUU )$ exists and is finite.
\end{proof} 

Proposition~\ref{eAB1}  follows from  Lemmas~\ref{eBN5} and~\ref{eBN8}.

\subsection{Proof of Proposition~$\ref{eAB1}$ in the case  $N\ge 3$ and $\beta = 0$}
Throughout this subsection, we assume that $N\ge 3$ and $\beta =0$, so that $\lambda _1 =\frac {1} {\alpha } > \lambda _2 =0$ and $\gamma = \frac {1} {\alpha } +1$. Moreover, $\phi _2(s)=(\log s)^{-\frac {1} {\alpha }}$.
We first study the limit~\eqref{eAB2}.

\begin{lemma} \label{eBN2b} 
Given any $\MUU \ge \MUUz$, the limit $\Lt( \MUU )$ given by~\eqref{eAB2} exists and is finite, and 
either $\Lt( \MUU )= 0$ or else $\Lt( \MUU )= \pm  ( \frac {2} {\alpha }  )^{\frac {2} {\alpha }} $.
Moreover,
\begin{equation} \label{fTPA44} 
\frac {s w_\MUU ' (s)} {w_\MUU (s)} \goto  _{ s\to \infty  } 0
\end{equation} 
if and only if $\Lt( \MUU )= \pm  ( \frac {2} {\alpha }  )^{\frac {2} {\alpha }} $.
\end{lemma}  

\begin{proof} 
Since $w_\MUU $ has a finite number of zeros,  we may suppose $w_\MUU (t)>0$ for $t$ large.
Set 
\begin{equation*} 
v(t)=t^{\frac {1} {\alpha }} w _\MUU (s) = \frac { w_\MUU (s)} {\phi _2(s)}
\end{equation*} 
 with $t=\log s$ and $s>1$, so that $v (t) >0$ for $t$ large and $\Lt ( \MUU )= \lim _{ t\to \infty  } v(t)$ if this last limit exists.
Note that
 \begin{equation} \label{fTPA94} 
\frac {sw _\MUU '(s)} {w _\MUU(s)}=\frac {v'(t)} {v(t)}-\frac {1} {\alpha t}. 
\end{equation} 
Since $w_\MUU $ is a solution of~\eqref{IPE}, a  straightforward calculation shows that $v$ satisfies
\begin{equation} \label{fFH18} 
4t v''+\Bigl (\frac {4t} {\alpha }-\frac {8} {\alpha }-te^{-t}\Bigr )v'-\Bigr (\frac {4} {\alpha ^2}-\frac {4(\alpha +1)} {\alpha ^2 t}-\frac {e^{-t}} {\alpha }\Bigl )v+|v|^\alpha v=0 .
\end{equation}  
Set
\begin{equation} \label{fFH19} 
\Psi (v)=\frac {1} {\alpha +2}|v|^{\alpha +2}-\frac {2} {\alpha ^2}|v|^2
\end{equation} 
and
\begin{equation}  \label{fFH19:1} 
V (t) =2t v' (t) ^2+\Bigl (\frac {2(\alpha +1)} { \alpha ^ 2 t}  +\frac {1 } {2\alpha  }e^{-t}\Bigr ) v (t) ^2 + \Psi (v (t) ) 
\end{equation} 
so that in particular $V (t)$ is bounded below. 
It follows using~\eqref{fFH18} that
\begin{equation} \label{fTTPPAA} 
V' (t) =-\Bigl (\frac {4t} {\alpha }-\frac {8} {\alpha }-2-te^{-t}\Bigr ) v' (t) ^2-\Bigl (\frac {2(\alpha +1)} {\alpha ^2 t^2} + \frac {1} {2\alpha }e^{-t}\Bigr )v (t)^2.
\end{equation} 
Therefore, $V' (t) <0$ for $t$ large, and $V$ is decreasing; and so $V(t)$ has a finite limit $V^\infty $ as $t\to \infty $.  
It  follows that $v (t) $  and $ t  |v'(t) |^2 $ are bounded as $t\to \infty$, hence $v'(t) \to 0$ as $t\to \infty $.

 We next show that $v(t) $ has a limit as $t\to \infty $, i.e. that $\Lt (\MUU)$ exists. For this, we set
 \begin{equation} \label{fTTPPAA1} 
 v^+ = \limsup _{ t\to \infty  } v (t), \quad v^- =  \liminf _{ t\to \infty  } v (t) .
 \end{equation} 
 Assuming $v^- < v^+$, it follows that $v$ oscillates asymptotically between $v^-$ and $v^+$, so that there exist $t_n^\pm  \to \infty $ such that $v' (t_n^\pm )   =0$, $v (t_n^\pm ) \to v^\pm$. 
Letting $t= t_n ^\pm$ in~\eqref{fFH19:1}, then $n\to \infty $, we deduce that
\begin{equation}  \label{fTTPPAA2} 
V^\infty = \Psi ( v^- ) = \Psi ( v^+) .
\end{equation} 
Since  $\Psi(v^+)=\Psi(v^-)$, there exists $ \underline{v}\in (v^- ,v^+ )$ such that $\Psi'(\underline{v})=0$. From~\eqref{fFH19} 
we obtain $ \underline{v}=(\frac {2} {\alpha } )^{\frac {2} {\alpha }}$. 
Consider now  $ \underline{v}<a<b<v^+$ and $\tau _n \uparrow \infty$  such that 
$v(\tau_{2n})= a$, $v(\tau_{2n+1})= b$, with $v(t)\in [a,b]$ for $t\in [\tau _{2n}, \tau_{2n+1}]$. Thus $v(t) ^\alpha \ge a^\alpha > \underline{v} ^\alpha =\frac {4} {\alpha ^2}$ for $t\in [\tau _{2n}, \tau_{2n+1}]$. On the other hand, we can write~\eqref{fFH18} in the form
\begin{equation} \label{fFH18:1} 
 \Bigl (v' + \Bigl( \frac {1} {\alpha } - \frac {2} {\alpha t} - \frac {e^{-t}} {4} \Bigr)v \Bigr )'  = \frac {1} {4t} \Bigr ( \underline{v}^\alpha  -  |v| ^\alpha  + \frac {4(\alpha -1)} {\alpha ^2 t }-\frac {1} {\alpha }e^{-t} +  t e^{-t}  \Bigl ) v .
\end{equation}   
We observe that the right-hand side of~\eqref{fFH18:1}  is negative on $ [\tau _{2n}, \tau_{2n+1}] $ for $n$ sufficiently large. Integrating on $ ( \tau _{2n}, \tau_{2n+1} ) $, we obtain
\begin{equation*} 
v' ( \tau_{2n+1} ) + \Bigl( \frac {1} {\alpha } - \frac {2} {\alpha \tau_{2n+1}} - \frac {e^{-\tau_{2n+1}}} {4} \Bigr) b \le  v' ( \tau_{2n } ) + \Bigl( \frac {1} {\alpha } - \frac {2} {\alpha \tau_{2n }} - \frac {e^{-\tau_{2n }}} {4} \Bigr) a.
\end{equation*} 
We now let $n\to \infty $. Since $v'( \tau _n) \to 0$, we obtain $b\le a$, which is absurd.
It follows that there exists $v_{\infty}\ge 0$ such that $v(t)\to v_{\infty}$ as $t\to \infty$.  

We now claim that either $v_\infty =0$ or else $v_\infty =  \underline{v} $. 
Indeed, suppose first $v_\infty \in  ( 0,  \underline{v}  )$, so that $(\underline{v}^\alpha  -  v_\infty  ^\alpha) v_\infty <0$. It follows from~\eqref{fFH18:1} that there exists $\delta >0$ such that
\begin{equation*} 
 \Bigl (v' + \Bigl( \frac {1} {\alpha } - \frac {2} {\alpha t} - \frac {e^{-t}} {4} \Bigr)v \Bigr )'  \le - \frac {\delta  } {t} 
\end{equation*}   
for $t $ large. Integrating the above inequality yields a contradiction with the fact that $v$ and $v'$ are bounded. We obtain in the same way a contradiction if we assume that $v_\infty >  \underline{v} $. 
Thus either $v(t) \to 0$, or else $v(t) \to (\frac {2} {\alpha })^{\frac {2} {\alpha }}$  as $t\to \infty $. This proves the first part of the statement. 

Suppose now $\Lt ( \MUU) = (\frac {2} {\alpha })^{\frac {2} {\alpha }}$, i.e. $v(t) \to (\frac {2} {\alpha })^{\frac {2} {\alpha }}$. Since $v'(t) \to 0$, it follows that 
\begin{equation} \label{fTPA82} 
\frac {v'(t)} {v(t)} \goto _{ t\to \infty  } 0.
\end{equation} 
Applying~\eqref{fTPA94},  we deduce that~\eqref{fTPA44} holds. Conversely, assume~\eqref{fTPA44}, so~\eqref{fTPA82} holds. 
We prove that $v(t) \to (\frac {2} {\alpha })^{\frac {2} {\alpha }}$ by contradiction. Otherwise, $v(t) \to 0$, and it follows from~\eqref{fFH18:1} that $v' + ( \frac {1} {\alpha } - \frac {2} {\alpha t} - \frac {e^{-t}} {4} )v$ is nondecreasing for $t$ large. Since both $v(t)$ and $v'(t)$ converge to $0$ as $t\to \infty $, we deduce that $v' + ( \frac {1} {\alpha } - \frac {2} {\alpha t} - \frac {e^{-t}} {4} )v\le 0$ for $t$ large. It follows that $\frac {v'} {v} \le - \frac {1} {2\alpha }$ for $t$ large, contradicting~\eqref{fTPA82}. This completes the proof. 
\end{proof}

\begin{lemma} \label{eFH1} 
Let $\MUU \ge \MUUz$ and let $s_\MUU >0$ be sufficiently large so that 
\begin{equation*} 
w_\MUU (s) w_\MUU ' (s) < 0, \quad s>s_\MUU.
\end{equation*} 
(See Proposition~$\ref{eEX1:b1}$~\eqref{eEX1:2:7}.) If
\begin{equation} \label{fFH1} 
h_{1, \MUU } (s) = \frac {sw _\MUU ' (s)} {w_\MUU  (s)} + \frac {1} {\alpha }, \quad h_{2, \MUU} (s) = \frac {sw _\MUU ' (s)} {w_\MUU  (s)}
\end{equation} 
for $s>s_\MUU $, then
\begin{align} 
\frac {d} {ds} ( s^{- \frac {N-2} {2} } e^{\frac {1} {4s }} h_{1, \MUU } (s) )  &= - s^{- \frac {N} {2} } e^{\frac {1} {4s }} \Bigl( \frac {1 } {4\alpha s} + \frac { |w_\MUU |^\alpha } {4} + h_{1, \MUU } ^2 (s) \Bigr) \label{fFH2} \\
\frac {d} {ds} ( s^{\frac {N-2} {2}} e^{\frac {1} {4s }} h_{2, \MUU } (s) )  &= - s^{\frac {N-4} {2} } e^{\frac {1} {4s }} \Bigl(  \frac { |w _\MUU |^\alpha } {4} + h_{2, \MUU }^2 (s)  \Bigr) \label{fFH2b}
\end{align} 
for all $s>s_ \MUU$ and $j=1, 2$.
\end{lemma} 

\begin{proof} 
Identity~\eqref{fFH2} is formula~(3.1) in~\cite{SoupletW}, and is straightforward to verify. Formula~\eqref{fFH2b}  can be obtained from~\eqref{fFH2} and the   identity $h_{2, \MUU } (s) = h_{1, \MUU } (s) -\frac {1} {\alpha }=h_{1, \MUU } (s) - \frac {N-2} {2}$.
\end{proof} 

\begin{lemma} \label{eFH2} 
Let $\MUU \ge \MUUz$. With the notation~\eqref{fFH1}, we have 
\begin{equation} \label{fFH3} 
h_{2, \MUU } (s) \goto  _{ s\to \infty  } \ell
\end{equation} 
with either $\ell =0$ or else $\ell = - \frac {N-2} {2}$.
\end{lemma}

\begin{proof} 
We first show that~\eqref{fFH3}  holds for some $\ell \in [-\infty , 0]$.
Note that $h_{ 2, \MUU } (s) <0$ for all $s>s_\MUU $. Suppose there exists $\sigma _\MUU  >s_ \MUU $ such that $h_{ 2, \MUU } ( \sigma _\MUU ) \le  - \frac {N-2} {2}$. This means that $h_{ 1, \MUU } ( \sigma _\MUU ) \le 0$, then by~\eqref{fFH2}, we deduce that $h_{ 1, \MUU } (s) <0$ for all $s> \sigma _\MUU$.   Furthermore, since the right-hand side of~\eqref{fFH2} is negative,
\begin{equation*} 
 s^{- \frac {N-2} {2}} e^{\frac {1} {4s }}    h_{ 1, \MUU } ' (s)   \le - h_{ 1, \MUU } (s) \frac {d} {ds} ( s^{- \frac {N-2} {2}} e^{\frac {1} {4s }}) = h_{ 1, \MUU } (s)  \Bigl( \frac {N-2} {2} + \frac {1} {4s} \Bigr) s^{- \frac {N } {2}} e^{\frac {1} {4s }} <0
\end{equation*} 
for $s> \sigma _\MUU$. Thus $h_{ 1, \MUU }$ is decreasing,  and so is $h_{ 2, \MUU }$. The desired conclusion  then follows.

It remains to consider the case $- \frac {N-2} {2} < h_{ 2, \MUU } (s) < 0$ for $s>s_ \MUU $. Suppose by contradiction that $h_{ 2, \MUU }$ does not have a limit. Let 
\begin{equation*} 
-\frac {N-2} {2}\le  \underline{ \ell } = \liminf  _{ s\to \infty  } h_{ 2, \MUU } (s) <  \overline{ \ell }  = \limsup  _{ s\to \infty  } h_{ 2, \MUU } (s) \le 0
\end{equation*} 
 and let $\ell \in  (  \underline{ \ell} ,  \overline{ \ell} )$.  Consider an increasing sequence $s _n \to \infty $ such that 
\begin{equation} \label{fFH4} 
 h_{ 2, \MUU } (s _n) = \ell \text{ and } h_{ 2, \MUU }( s) \ge \ell \text{ for }s  _{ 2n }\le s\le s  _{ 2n+1 }.
\end{equation} 
Integrating~\eqref{fFH2b} on $(s _{ 2n }, s _{ 2n+1 })$ and applying~\eqref{fFH4} we obtain
\begin{equation*} 
\begin{split} 
( s _{ 2n+1 }^{\frac {N-2} {2}} e^{\frac {1} {4s _{ 2n+1 } }} - s _{ 2n }^{\frac {N-2} {2}} e^{\frac {1} {4s _{ 2n} }} ) \ell &  = - \int  _{ s _{ 2n } }^{s _{ 2n+1 }}  s^{\frac {N-4} {2} } e^{\frac {1} {4s }} \Bigl(  \frac { |w _\MUU |^\alpha } {4} + h_{ 2, \MUU } (s) ^2 \Bigr)\, ds \\ & \ge  - \int  _{ s _{ 2n } }^{s _{ 2n+1 }}  s^{\frac {N-4} {2} } e^{\frac {1} {4s }} \Bigl( \frac { |w _\MUU |^\alpha } {4} + \ell ^2 \Bigr)\, ds .
\end{split} 
\end{equation*} 
Since $w _\MUU (s) \to 0$ as $s \to \infty $, given any $0<\delta  < \ell ^2$, we have $ \frac { |w _\MUU (s)|^\alpha } {4}< \delta   $ for $s$ large, and we deduce that 
\begin{equation*} 
( s _{ 2n+1 }^{\frac {N-2} {2}} e^{\frac {1} {4s _{ 2n+1 } }} - s _{ 2n }^{\frac {N-2} {2}} e^{\frac {1} {4s _{ 2n} }} ) \ell  \ge - \int  _{ s _{ 2n } }^{s _{ 2n+1 }}  s^{\frac {N-4} {2} } e^{\frac {1} {4s }} ( \ell ^2 +\delta  )\, ds .
\end{equation*} 
Since $\ell <0$, we deduce that
\begin{equation*} 
 \frac { s _{ 2n+1 }^{\frac {N-2} {2}} e^{\frac {1} {4s _{ 2n+1 } }} - s _{ 2n }^{\frac {N-2} {2}} e^{\frac {1} {4s _{ 2n} }} } {\int  _{ s _{ 2n } }^{s _{ 2n+1 }}  s^{\frac {N-4} {2} } e^{\frac {1} {4s }} ds } \le \frac { \ell ^2 +\delta  } {- \ell} .
\end{equation*} 
Letting $n\to \infty $ (and since $e^{\frac {1} {4s }} \to 1$ as $s\to \infty $), we conclude that 
\begin{equation*} 
 \frac {N-2} {2} \le  \frac { \ell ^2 +\delta  } {- \ell} .
\end{equation*}
Letting now $\delta  \downarrow 0$, we obtain $\ell \le - \frac {N-2} {2}$, which is absurd since $\ell >  \underline{\ell} \ge - \frac {N-2} {2}$. 

Thus we have shown that~\eqref{fFH3}  holds for some $\ell \in [-\infty , 0]$, and we finally prove that $\ell =0 $ or $\ell = - \frac {N-2} {2}$. 
Integrating~\eqref{fFH2b} on $(1, s)$, we obtain
\begin{equation*} 
 e^{\frac {1} {4s }} h_{ 2, \MUU } (s)   =  s^{- \frac {N-2} {2}} e^{\frac {1} {4 }} h_{ 2, \MUU } (1) -  s^{ - \frac {N-2} {2}} \int _1 ^s \tau ^{\frac {N-4} {2} } e^{\frac {1} {4\tau  }} \Bigl( \frac { |w _\MUU |^\alpha } {4} + h_{ 2, \MUU } (\tau ) ^2 \Bigr) \, d\tau ,
\end{equation*} 
and so, letting $s\to \infty $ and applying~\eqref{fFH3}, 
\begin{equation} \label{fFH5} 
\ell    =  - \lim  _{ s\to \infty  } \frac { \int _1 ^s \tau ^{\frac {N-4} {2} } e^{\frac {1} {4\tau  }} \Bigl( \frac { |w _\MUU |^\alpha } {4} + h_{ 2, \MUU } (\tau ) ^2 \Bigr) \, d\tau } { s^{  \frac {N-2} {2}}} . 
\end{equation} 
  Both the numerator and the denominator in the right-hand side of~\eqref{fFH5} go to $\infty $ with $s$. By l'H\^opital's rule
\begin{equation*} 
\ell    =  - \lim  _{ s\to \infty  } \frac {s ^{\frac {N-4} {2} } e^{\frac {1} {4s  }}  ( \frac { |w _\MUU |^\alpha } {4} + h_{ 2, \MUU } (s ) ^2  ) } {\frac {N-2} {2} s^{  \frac {N-4} {2}}} = - \frac {2} {N-2} \ell ^2.
\end{equation*} 
Therefore, either $\ell = - \frac {N-2} {2}$ or else $\ell=0$.
\end{proof} 

\begin{lemma} \label{eFH3} 
Given $\MUU \ge 1$, if  $\Lt( \MUU )=0 $, then the limit~\eqref{eAB4}  exists and is finite.
\end{lemma}

\begin{proof} 
We follow the argument of~\cite[proof of Proposition~3.5]{SoupletW}. 
Suppose $\Lt( \MUU )=0 $. It follows from Lemma~\ref{eBN2b} that~\eqref{fTPA44} does not hold. Therefore, \eqref{fFH3} holds with $\ell = - \frac {N-2} {2}$, so that $h_{ 1, \MUU }(s) \to 0$ as $s\to \infty $. 
Since the right-hand side of~\eqref{fFH2}  is  negative, we deduce that $h_{ 1, \MUU }(s) >0$ for $s> s_\MUU$.
Moreover, integrating~\eqref{fFH2} on $(s, \infty )$, we obtain
\begin{equation*} 
\begin{split} 
s^{- \frac {N-2} {2} } e^{\frac {1} {4s }} h_{ 1, \MUU } (s) &  = \int  _{ s }^\infty  \tau ^{- \frac {N} {2} } e^{\frac {1} {4\tau  }} \Bigl( \frac { 1 } {4 \alpha \tau } + \frac { |w _\MUU |^\alpha } {4} + h_{ 1, \MUU } (\tau ) ^2 \Bigr) \, d\tau \\
& \le e^{\frac {1} {4s  }}  \Bigl( \frac { 1 } {4 \alpha s } + \frac { |w _\MUU (s) |^\alpha } {4} + \sup  _{ \tau \ge s } h_{ 1, \MUU } (\tau ) ^2 \Bigr) \int _s^\infty \tau ^{- \frac {N} {2} } d\tau  \\
& = \frac {2} {N-2} s^{- \frac {N-2} {2} } e^{\frac {1} {4s  }}  \Bigl( \frac {1 } {4 \alpha s } + \frac { |w _\MUU (s) |^\alpha } {4} + \sup  _{ \tau \ge s } h_{ 1, \MUU } (\tau ) ^2 \Bigr) .
\end{split} 
\end{equation*} 
This implies
\begin{equation*} 
 \sup  _{ \tau \ge s } h_{ 1, \MUU } (\tau ) \le \frac {2} {N-2}  \Bigl( \frac { 1 } {4 \alpha  s } + \frac { |w _\MUU (s) |^\alpha } {4} +  \sup  _{ \tau \ge s } h_{ 1, \MUU } (\tau ) ^2 \Bigr) 
\end{equation*} 
so that (since $h_{ 1, \MUU }(s) \to 0$)
\begin{equation} \label{fFH6} 
 \sup  _{ \tau \ge s } h_{ 1, \MUU } (\tau ) \le  \frac {4} {N-2}    \Bigl( \frac { 1 } {4 \alpha s } + \frac { |w _\MUU (s) |^\alpha } {4}  \Bigr) 
\end{equation} 
for $s$ large. Since $h_{ 1, \MUU }(s) \to 0$, we have $\frac {w_\MUU '} {w _\MUU} + \frac {1} {2\alpha s}\le 0$ for $s$ large, by~\eqref{fFH1}. It follows that $w _\MUU (s) \le C s^{- \frac {1 } {2\alpha } }$, so that by~\eqref{fFH6},  $h_{ 1, \MUU } (s) \le C s^{ -\frac {1} {2} }$. 
Therefore, $\frac {w _\MUU  '(s)} {w _\MUU (s)} + \frac {1} {\alpha s} = O( s^{-\frac {3} {2} })$,
from which the existence of the limit~\eqref{eAB4} easily follows.
\end{proof} 

Proposition~\ref{eAB1}  follows from Lemmas~\ref{eBN2b} and~\ref{eFH3}. 

\subsection{Proof of Proposition~$\ref{eAB1}$ in the case  $N=2$} \label{NeqD} 
We set $\lambda =\frac {1} {\alpha }$, so that $\phi _1(s)=s^{-\lambda }$ and $\phi _2 (s) = s^{-\lambda }\log s$. 
It is straightforward to verify the following variation of the parameter formula
\begin{equation} \label{fNd1} 
w_\MUU (s) =c_{1,\MUU}(s)\phi _1 (s) +c_{2,\MUU}(s)\phi _2 (s), 
\end{equation} 
where 
\begin{equation} \label{fNd1bis}
\begin{split}
&c_{1,\MUU}(s)=w_{\MUU}(1)-\frac {1} {4}\int_1^s \tau^{\lambda -1} (\log \tau ) \Phimu_\MUU(\tau)\, d\tau,\\
&c_{2,\MUU}(s)=\lambda w_{\MUU}(1)+w'_{\MUU}(1) +\frac {1} {4}\int_1^s \tau^{\lambda -1} \Phimu_\MUU(\tau)\, d\tau
\end{split}
\end{equation}  
and $\Phimu _\MUU$ is given by~\eqref{fTPA55}.

\begin{lemma}\label{eNd10} 
Given $ \overline{\MUU} \ge  \MUUz $, there exists $C>0$ such that
\begin{equation} \label{fNd2d} 
\sup_{\MUUz \le \MUU \le  \overline{\MUU} }\, \sup_{s \ge 0} \,  \frac {(1+ s)^{\lambda }} {\log (2+s)}  \, ( (1 + s) | w_\MUU '(s) |+| w_\MUU (s)|)\le C
\end{equation} 
and
\begin{equation} \label{fNd2e} 
\sup_{\MUUz \le \MUU \le  \overline{\MUU} }\, \sup_{s \ge 0} \,  \frac {(1+ s)^{\lambda +1}} { (\log (2+s))^{\alpha +1}}    |\Phimu_\MUU (s)|\le C.
\end{equation} 
\end{lemma} 

\begin{proof} 
$z(s)=s^\lambda  w_\MUU (s)$ is a solution of equation~\eqref{fBN5a1} and it follows from~\eqref{fBN5a3} and~\eqref{fBN5a4} that 
\begin{equation*} 
\sup _{ s\ge 1 }   \Bigl( s  | z' (s)| + \frac {1} {\log (2+s)}  |z (s)|  \Bigr)\le C
\end{equation*} 
where $C$ depends only on $ |z(1)|+  |z'(1)|$. 
This implies that
\begin{equation*} 
 \sup_{s \ge 1} \,  \frac {(1+ s)^{\lambda }} {\log (2+s)}  \, ( (1 + s) | w_\MUU '(s) |+| w_\MUU (s)|)\le C
\end{equation*} 
where $C$ depends only on $ |w_\MUU (1)|+  |w_\MUU '(1)|$. 
Estimate~\eqref{fNd2d}  now follows from the continuous dependence of $w_ \MUU  $ and $w'_ \MUU  $ on $ \MUU  $ given by Proposition~\ref{eEX1:b1}~\eqref{eEX1:6:b1}. 
 Finally, \eqref{fNd2e} is a direct consequence of~\eqref{fNd2d}. 
\end{proof} 

\begin{lemma} \label{eTPA2} 
For all $\MUU \ge \MUUz$, the map $\tau \mapsto \tau^{\lambda -1} (\log \tau ) \Phimu_\MUU(\tau)$ is integrable at $\infty $. Moreover, the functions $c _{ j, \MUU } (s)$ defined by~\eqref{fNd1bis}  have finite limits $c _{ j, \MUU }^\infty $ as $s\to \infty $, given by
\begin{equation} \label{fNd1t}
\begin{split}
&c_{1,\MUU}^\infty =w_{\MUU}(1)-\frac {1} {4}\int_1^\infty  \tau^{\lambda -1} (\log \tau ) \Phimu_\MUU(\tau)\, d\tau,\\
&c_{2,\MUU}^\infty =\lambda w_{\MUU}(1)+w'_{\MUU}(1) +\frac {1} {4}\int_1^\infty  \tau^{\lambda -1} \Phimu_\MUU(\tau)\, d\tau .
\end{split}
\end{equation}  
In addition, the maps $\MUU \mapsto c_{j,\MUU}^\infty $ are continuous $[\MUUz, \infty ) \to \R$ and, 
given $ \overline{\MUU} \ge  \MUUz $, there exists a constant $C$ such that 
\begin{equation} \label{fTPA35} 
 | c _{ j, \MUU } ^\infty  -  c _{ j, \MUU } (s)| \le C s^{- \frac {1} {2}}
\end{equation} 
 for  $j=1,2$, $s\ge 1$ and $\MUUz \le \MUU \le   \overline{\MUU} $.
\end{lemma} 

\begin{proof} 
Fix $ \overline{\MUU} > \MUUz$. 
It follows from~\eqref{fNd2e} that there exists a constant $C$ such that
\begin{equation} \label{fTPA36} 
\tau^{\lambda -1} (\log \tau )  |\Phimu_\MUU(\tau) | \le C \tau ^{-2} ( \log (2 +\tau ))^{\alpha +2}
\end{equation} 
for all $\tau \ge 1$ and all $\MUUz \le \MUU \le   \overline{\MUU} $. This shows the integrability property, the existence of limits $c _{ j, \MUU }^\infty $, formulas~\eqref{fNd1t} and estimates~\eqref{fTPA35}. We now prove the continuity. 
Let $\MUU _1, \MUU _2 \in [ \MUUz ,  \overline{\MUU} )$. Given $s\ge 1$, we deduce from~\eqref{fTPA35} that
\begin{equation*} 
 |c  _{ j, \MUU_1 }^\infty - c  _{ j, \MUU_2 }^\infty | \le   |c  _{ j, \MUU_1 } (s) - c  _{ j, \MUU_2 } (s) | + C s^{- \frac {1} {2}} .
\end{equation*} 
Given $\varepsilon >0$, we first fix $s_0$ sufficiently large so that $C s_0^{- \frac {1} {2}}\le \frac {\varepsilon } {2}$. By continuous dependence (Proposition~\ref{eEX1:b1}~\eqref{eEX1:6:b1}), if $ | \MUU_1 - \MUU_2 |$ is sufficiently small, then $ |c  _{ j, \MUU_1 } (s_0) - c  _{ j, \MUU_2 } (s_0) | \le \frac {\varepsilon } {2}$. Therefore, $|c  _{ j, \MUU_1 }^\infty - c  _{ j, \MUU_2 }^\infty | \le  \varepsilon $, which proves continuity on $[\MUUz ,  \overline{\MUU} )$. Since $ \overline{\MUU} $ is arbitrary, this completes the proof.
\end{proof} 

\begin{lemma}\label{eNd11} 
Let $\MUU \ge \MUUz$, and let $c _{ 1, \MUU }^\infty $ and $c _{ 2, \MUU }^\infty $ be given by~\eqref{fNd1t}.  The limit~\eqref{eAB2} exists and $\Lt ( \MUU) = c _{ 2, \MUU }^\infty $. Moreover, if $\Lt(\MUU)=0$ then the limit~\eqref{eAB4} exists and $\Ltu(\MUU)= c _{ 1, \MUU }^\infty $.
\end{lemma}

\begin{proof} 
By formula~\eqref{fNd1}, 
\begin{equation*} 
\frac {w_\MUU (s)} {\phi _2 (s)} = \frac {c _{ 1, \MUU } (s) } {\log s} + c _{ 2, \MUU  }(s).
\end{equation*}  
It follows from Lemma~\ref{eTPA2} that the limit~\eqref{eAB2}  is well defined and $\Lt ( \MUU ) = c _{ 2, \MUU }^\infty $.
Moreover, it follows from formula~\eqref{fNd1} again that
\begin{equation*} 
\frac {w_\MUU (s)} {\phi _1 (s)} =   c _{ 1, \MUU } (s)   + c _{ 2, \MUU  }(s)   \log s.
\end{equation*} 
Therefore, if $\Lt ( \MUU ) =0$ (hence $c _{ 2, \MUU }^\infty =0$), then by~\eqref{fTPA35}
\begin{equation*} 
 \Bigl|  \frac {w_\MUU (s)} {\phi _1 (s)} -   c _{ 1, \MUU }^\infty   \Bigr| \le  | c _{ 1, \MUU } (s) - c _{ 1, \MUU }^\infty  | + | c _{ 2, \MUU } (s) - c _{ 2, \MUU }^\infty  | \log s \goto  _{ s\to \infty  }0.
\end{equation*}  
Thus we see that  $\Ltu(\MUU)$ is well defined and $ \Ltu(\MUU)=c _{ 1, \MUU }^\infty$. 
\end{proof} 

Proposition~\ref{eAB1}  follows from   Lemma~\ref{eNd11}.

\section{Local behavior of the number of zeros} \label{sFNL} 

This section is devoted to the proof  of Proposition~\ref{eEX2}. We consider separately the cases $N\ge 3$ and $N=2$.

\subsection{The case $N\ge 3$}

The proof is inspired by the proof of~\cite[Lemma~4]{McLeodTW}.
Since  $\lambda _1 > \frac {\lambda _1 + \lambda _2} {2} >  \max \{ 0, \lambda _2  \}$ we may fix
\begin{equation} \label{fNZE1} 
 \quad  \sigma \in (\max \{ 0, \lambda _2 \}, \textstyle { \frac {\lambda _1 + \lambda _2} {2} } ).
\end{equation} 
We set
\begin{equation} \label{fNZE2}  
a_1= 1 + \lambda _1 + \lambda _2 -2 \sigma , \quad a_2 =-2 (\sigma - \lambda _1) (\sigma - \lambda _2)  , \quad a_3 =  4 (\lambda _1 + \lambda _2 -2\sigma )
\end{equation} 
so that by~\eqref{fNZE1} 
\begin{equation}  \label{fNZE4} 
a_1 >1, \quad a_2 >0, \quad a_3 >0.
\end{equation} 
Given $\MUU \ge \MUUz$ we set 
\begin{equation}  \label{fNZE5} 
 z _\MUU (s) = s^\sigma w_\MUU (s)
\end{equation} 
so that (see~\eqref{fBN5a})
\begin{equation} \label{fNZE6} 
4s^2  z _\MUU  ''+ 4 a_1 s  z _\MUU  ' -  z '  - 2 a_2  z _\MUU   + \frac {\sigma } {s}  z _\MUU   + s^{-\alpha \sigma }| z _\MUU  |^\alpha  z _\MUU  =0.
\end{equation}
Moreover, we set
\begin{equation} \label{fNZE7} 
 \widetilde{H}  _\MUU  (s) = 2 s^2  |  z  _\MUU ' |^2 - a_2 z _\MUU  ^2 + \frac {\sigma } {2 s} z _\MUU  ^2 + \frac {s^{-\sigma \alpha }} {\alpha +2}  |  z _\MUU   |^{\alpha +2} .
\end{equation} 
Elementary calculations show that
\begin{equation}  \label{fNZE8} 
 \widetilde{H}  _\MUU '(s)   = -\frac {\alpha \sigma s^{- \alpha \sigma -1} } { \alpha +2}  |  z _\MUU   |^{\alpha +2} -\frac {\sigma } {2s^2}  z _\MUU  ^2 +(1 - a_3 s)| z _\MUU  '|^2 .
\end{equation} 
In particular,
\begin{equation}  \label{fNZE8b1} 
 \widetilde{H}  _\MUU  '(s)  \le 0, \quad s\ge \frac {1} {a_3}.
\end{equation} 
Given $\delta >0$, it follows from~\eqref{fNZE8} that 
\begin{equation*} 
\begin{split} 
s^{ - \delta} \frac {d} {ds} [ s^\delta   \widetilde{H}  _\MUU  (s)    ] & =  \widetilde{H}  _\MUU ' (s)   +\frac {\delta } {s}  \widetilde{H}  _\MUU  (s)   \\
 = - [ ( a_3 - 2\delta ) s -1]  |  z _\MUU  ' |^2 & - \frac { z _\MUU  ^2} { 2 s^2} [ (1- \delta ) \sigma + 2 \delta a_2 s ] - \frac {s^{- \alpha \sigma -1}  |  z _\MUU   |^{\alpha +2}} {\alpha +2} [ \alpha \sigma -\delta ] .
\end{split} 
\end{equation*} 
Therefore, if we fix
\begin{equation} \label{fNZE10} 
0 < \delta < \min  \Bigl\{ \frac {a_3} {4}, 1, \alpha \sigma  \Bigr\}
\end{equation} 
then 
\begin{equation}  \label{fNZE11} 
 \frac {d} {ds} [ s^\delta  \widetilde{H}  _\MUU  (s)   ]   \le - \delta a_2 s^{ \delta -1 } z _\MUU  (s)^2, \quad s\ge \frac {2} {a_3}.
\end{equation} 
Suppose first $\Lt (  \overline{\MUU }) \not = 0 $, for instance $\Lt (  \overline{\MUU })  > 0 $. 
Since $\sigma > \max\{0, \lambda _2 \}$, we see   that $ z  _{  \overline{\MUU}  } (s) \to \infty $, and it follows from~\eqref{fNZE11} that $ \frac {d} {ds} [ s^\delta  \widetilde{H}   _{  \overline{\MUU}  }  (s)   ]  \le - c s^{\delta  -1}$ for $s$ large, with $c>0$. Since $s^{\delta  -1}$ is not integrable at $\infty $, there exists $s_0\ge \frac {2} {a_3}$  such that $\widetilde{H}   _{  \overline{\MUU}  }  (s_0) <0$. By continuous dependence, we have $\widetilde{H}  _\MUU  ( s_0 ) <0 $ for $\MUU $ close to $ \overline{\MUU} $. Therefore, by~\eqref{fNZE11}, $\widetilde{H}  _\MUU  ( s ) <0 $ for $s\ge s_0$ and $\MUU $ close to $ \overline{\MUU} $, and it follows from~\eqref{fNZE7} that $w_\MUU$ does not vanish  for $s\ge s_0$ and $\MUU $ close to $ \overline{\MUU} $. Again by continuous dependence of $w_\MUU$ on $\MUU$ in $C^1([0,  s_0 ]) $, we conclude that $w_\MUU$ has the same number of zeros as $w _{  \overline{\MUU}  }$ provided $\MUU $ is sufficiently close to $ \overline{\MUU} $. This proves Property~\eqref{eEX2:1} of Proposition~\ref{eEX2}.

We now assume $\Lt (  \overline{\MUU} )=0$.
It follows from Proposition~\ref{eAB1}~\eqref{eAB1:2}  that $\Ltu ( \overline{\MUU} )$ exists and is finite,  therefore by~\eqref{fNZE1}  
\begin{equation}  \label{fNZE12} 
z _{  \overline{\MUU }  } (s) \goto _{ s\to \infty  }0. 
\end{equation} 
Moreover, we deduce from equation~\eqref{fNZE6} that for large $s$, if $ z _{  \overline{\MUU }  } '$ vanishes, then $z _{  \overline{\MUU }  } '' $ has the sign of $z _{  \overline{\MUU }  }$. It easily follows that $ z _{  \overline{\MUU }  } '$ has constant sign for $s$ large. Therefore, we may assume without loss of generality that there exists
\begin{equation} \label{fNZE13b2} 
 \overline{s} \ge  \frac {2} {a_3} 
\end{equation} 
such that
\begin{equation}  \label{fNZE13} 
z _{  \overline{\MUU }  } (s) >0  \text{ and }  z _{  \overline{\MUU }  } '(s) <0  \text{ for } s\ge  \overline{s} .     
\end{equation} 
We now let $\varepsilon >0$ to be specified later, and we consider $ \MUU  \ge  \MUUz$ such that
\begin{equation} \label{fNZE13b1} 
 |  \MUU   -  \overline{ \MUU  } | \le \varepsilon .
\end{equation} 
It follows easily from continuous dependence (Proposition~\ref{eEX1:b1}~\eqref{eEX1:6:b1}), and the fact that if $w_\MUU (s)= 0$ then $w_\MUU '(s) \not = 0$, that $w_\MUU $ has  $\wN( \overline{ \MUU } )$ zeros on $(0,  \overline{s} ) $ if $\varepsilon $ is sufficiently small. This means that $z_\MUU $ has $\wN( \overline{ \MUU } ) $ zeros on $(0,  \overline{s}) $. Moreover, by choosing $\varepsilon $ possibly smaller, we have
\begin{equation}  \label{fNZE14} 
z _\MUU  (  \overline{s} ) >0  \text{ and } z _\MUU ' ( \overline{s} ) <0 .
\end{equation} 
If $z _\MUU >0$ on $( \overline{s}, \infty ) $, then $w _\MUU $ has $\wN (  \overline{\MUU } ) $ zeros on $(0,\infty )$. 

Assume now that $z _\MUU $ has a zero on $(  \overline{s} , \infty )$, and let $ s_ \MUU  >  \overline{s} $ be the smallest such zero. By~\eqref{fNZE13} and continuous dependence,  
\begin{equation}  \label{fNZE15} 
 s_ \MUU  \goto  _{ \varepsilon \to 0 } \infty . 
\end{equation} 
In particular, if $\varepsilon $ is sufficiently small, then
\begin{equation} \label{fNZE16} 
 s_ \MUU \ge \frac {\sigma } {a_2}  .
\end{equation} 
Moreover, $ z _\MUU ' (  s_ \MUU  ) <0$; and so we may define $ s_ \MUU  < \tau _\MUU  \le \infty  $ by
\begin{equation} \label{fNZE17} 
\tau _\MUU = \sup \{ s>   s_ \MUU  ; \, z _\MUU ' <0  \text{ on } (  s_ \MUU , s )  \} .
\end{equation} 
If $\tau _\MUU =\infty $, then $z _\MUU $ has $\wN (  \overline{\MUU } )  +1$ zeros on $(0, \infty )$ and $\limsup z _\MUU  (s) <0$ as $s\to \infty $. Therefore, $ \limsup \frac {w _\MUU (s)} {\phi _1(s)}=\infty$, showing that $\Lt( \MUU )\ne 0$.
If $\tau _\MUU  <\infty $, then $z _\MUU  (\tau _\MUU  ) <0$, $z _\MUU ' (\tau _\MUU  )= 0$ and $z  _\MUU '' (\tau ) \ge 0$, so that by~\eqref{fNZE6} 
\begin{equation*} 
 |z _\MUU  (\tau  _\MUU )|^\alpha \ge \tau _\MUU  ^{\alpha \sigma }  \Bigl( 2a_2 - \frac {\sigma } {\tau _\MUU  } \Bigr).
\end{equation*} 
Applying~\eqref{fNZE16}, we deduce that $ |z _\MUU ( \tau _\MUU  )| \ge a_2^{\frac {1} {\alpha }} \tau _\MUU  ^{ \sigma }   $.
Thus we see that there exist $ s _\MUU  < t_1 < t_2\le  \tau  _\MUU  $ such that
\begin{equation} \label{fNZE19} 
 z _\MUU  ( t_1 ) = - \frac {1} {2} a_2^{\frac {1} {\alpha }} \tau _\MUU  ^{ \sigma }  , \quad  z _\MUU  ( t_2 ) = -  a_2^{\frac {1} {\alpha }} \tau _\MUU  ^{ \sigma }  
 \end{equation} 
and
\begin{equation} \label{fNZE20} 
 - \frac {1} {2} a_2^{\frac {1} {\alpha }} \tau _\MUU  ^{ \sigma }\ge z _\MUU (s) \ge -  a_2^{\frac {1} {\alpha }} \tau _\MUU  ^{ \sigma }, \quad   t_1 \le s\le t_2 .
\end{equation} 
It follows from~\eqref{fNZE11}, \eqref{fNZE13b2}  and~\eqref{fNZE20} that
\begin{equation*}  
\begin{split} 
t_2 ^\delta   \widetilde{H}  _\MUU (t_2 ) -  \overline{s} ^\delta  \widetilde{H}  _\MUU (  \overline{s}  ) & \le    - \delta a_2 \int  _{   \overline{s}  }^{t_2}  s^{ \delta -1 } z _\MUU (s)^2   \le   - \delta a_2 \int  _{t_1} ^{t_2}  s^{ \delta -1 } z _\MUU  (s)^2 \\ & \le    - \frac {\delta } {4} a_2^{\frac {\alpha +2} {\alpha }}  \tau _\MUU ^{2\sigma } \int  _{t_1} ^{t_2}  s^{ \delta -1 }    =     - \frac {1 } {4} a_2^{\frac {\alpha +2} {\alpha }}  \tau _\MUU ^{2\sigma } (t_2 ^\delta - t_1^\delta ) .
\end{split} 
\end{equation*} 
Since $s ^\delta  \widetilde{H}  _\MUU (s)$ is nonincreasing on $( \overline{s} , \infty )$ by~\eqref{fNZE11} and~\eqref{fNZE13b2}, and since by continuous dependence, $ \overline{s} ^\delta  \widetilde{H}  _\MUU  (  \overline{s}  )  $ is bounded uniformly on $ \MUU  $ satisfying~\eqref{fNZE13b1}, we have
\begin{equation} \label{fNZE21} 
\begin{split} 
\tau _\MUU  ^\delta   \widetilde{H}  _\MUU (\tau _\MUU ) & \le t_2 ^\delta   \widetilde{H}  _\MUU (t_2 ) \le   \overline{s} ^\delta  \widetilde{H}  _\MUU  (  \overline{s}  )  - \frac {1 } {4} a_2^{\frac {\alpha +2} {\alpha }}  \tau _\MUU ^{2\sigma } (t_2 ^\delta - t_1^\delta )  \\ &  \le  C - \frac {1 } {4} a_2^{\frac {\alpha +2} {\alpha }}  \tau _\MUU ^{2\sigma } (t_2 ^\delta - t_1^\delta ) 
\end{split} 
\end{equation} 
for some constant $C$ independent of $ \MUU  $. 
Moreover, it follows from~\eqref{fNZE19} that
\begin{equation}   \label{fNZE22} 
\frac {1} {2} a_2^{\frac {1} {\alpha }} \tau  _\MUU ^{ \sigma } =  \Bigl| \int  _{ t_1 }^{t_2}  z _\MUU '  \Bigr| .
\end{equation} 
Furthermore,   $ \widetilde{H}  _\MUU $ is nonincreasing on $( \overline{s} , \infty )$  by~\eqref{fNZE8b1} and~\eqref{fNZE13b2}, from which it follows (using~\eqref{fNZE7} 
and~\eqref{fNZE20}) that 
\begin{equation}  \label{fNZE23} 
s^2  |z _\MUU ' (s)|^2  \le \frac {1} {2} \widetilde{H}  _\MUU ( \overline{s} ) + \frac {a_2} {2}  z _\MUU  (s)^2 \le C +  \frac {a_2^{\frac {\alpha +2} {\alpha }}} {8} \tau _\MUU ^{2\sigma } , \quad t_1\le s\le t_1.
\end{equation} 
In particular, 
\begin{equation}  \label{fNZE24} 
\int  _{ t_1 }^{t_2}  |z _\MUU '| \le C \tau  _\MUU ^\sigma \log \frac {t_2 } {t_1}.
\end{equation} 
Formulas~\eqref{fNZE22} and~\eqref{fNZE24} imply that  $\log \frac {t_2 } {t_1}$ is bounded below, so that there exists $\nu >0$ independent of $ \MUU  $ satisfying~\eqref{fNZE13b1} such that $t_2 \ge (1+ \nu ) t_1$. 
Since $t_1 \ge  s_\MUU \to \infty $ as $\varepsilon \to 0$ by~\eqref{fNZE15}, we see that if $\varepsilon $ is sufficiently small, then $t_2 ^\delta - t_1^\delta\ge 1$. Therefore, we deduce from~\eqref{fNZE21} and the fact that $\tau  _\MUU  \ge  s_\MUU \to \infty $ as $\varepsilon \to 0$, that $ \widetilde{H} _\MUU (\tau  _\MUU ) <0$ if $\varepsilon $ is sufficiently small. Applying~\eqref{fNZE8b1} and~\eqref{fNZE13b2}, we conclude that $ \widetilde{H}  _\MUU  (s ) \le  \widetilde{H}  _\MUU  (\tau _\MUU ) <0$  for $s\ge \tau  _\MUU $. It now follows from~\eqref{fNZE7} that there exists $\rho >0$ such that $z  _\MUU (s) \le -\rho $ for $s\ge \tau _\MUU $. Therefore  $z _\MUU $ has $\wN (  \overline{\MUU } )   +1$ zeros on $(0, \infty )$ and $ \frac {w _\MUU (s)} {\phi _1(s)}\to \infty$. By Proposition~\ref{eAB1}~\eqref{eAB1:2} we deduce that $\Lt( \MUU )\ne 0$. This proves Property~\eqref{eEX2:2} of Proposition~\ref{eEX2}, and completes the proof in the case $N\ge 3$.

\subsection{The case $N=2$}

We continue with the notation established in Subsection~\ref{NeqD}.
We first consider $ \overline{\MUU} \ge \MUUz$ such that $\Lt (  \overline{\MUU}  )\not = 0$.
It follows that there exists $ \overline{s} >0$ such that
\begin{equation} \label{fTPAC} 
 | w_{ \overline{\MUU } } (s)|\ge \frac {1} {2}  |\Lt (  \overline{\MUU}  )| \phi _2 (s)
\end{equation} 
for $s\ge  \overline{s} $.
We deduce from formula~\eqref{fNd1} that
\begin{equation*} 
w_\MUU (s) - c _{ 2, \MUU }^\infty \phi _2 (s)  = c_{1,\MUU}(s)\phi _1 (s) + (c_{2,\MUU}(s)- c _{ 2, \MUU }^\infty )\phi _2 (s) 
\end{equation*} 
so that, by continuity of  the map $\MUU \mapsto c_{1,\MUU}^\infty $ (Lemma~\ref{eTPA2}) and estimate~\eqref{fTPA35}
\begin{equation*} 
| w_\MUU (s) - c _{ 2, \MUU }^\infty \phi _2 (s) |  \le  C \phi _1 (s)  + C s^{-\frac {1} {2}}\phi _2 (s) 
\end{equation*} 
for some constant $C$ independent of $s\ge 1$ and $\MUU$ with $1\le \MUU \le  \overline{\MUU} +1 $. 
Thus we see that there exists $  \widetilde{s} \ge 1$ such that
\begin{equation*} 
| w_\MUU (s) - c _{ 2, \MUU }^\infty \phi _2 (s) | \le \frac {1} {10}  |\Lt (  \overline{\MUU}  )| \phi _2 (s)
\end{equation*} 
for $s\ge  \widetilde{s} $ and $\MUUz \le \MUU  \le  \overline{\MUU } +1$. 
Therefore,
\begin{equation*} 
\begin{split} 
 |w_\MUU (s) - w _{ \overline{\MUU } } (s)|& \le | w_\MUU (s) - c _{ 2, \MUU }^\infty \phi _2 (s) |+| c _{ 2, \overline{\MUU } }^\infty- c _{ 2, \MUU }^\infty | \phi _2 (s) +| w_{ \overline{\MUU } } (s) - c _{ 2,  \overline{\MUU }  }^\infty \phi _2 (s) |
 \\& \le \frac {1} {4}  |\Lt (  \overline{\MUU}  )| \phi _2 (s)
\end{split} 
\end{equation*} 
 for $s\ge  \widetilde{s} $ and $\MUU $ close to $ \overline{\MUU} $.
Since
$ |w_\MUU (s)|\ge  |w _{ \overline{\MUU}  } (s)| -  |w_\MUU (s) - w _{ \overline{\MUU } } (s)|$, we deduce by applying~\eqref{fTPAC}  that   if $\MUU $ is close to $ \overline{\MUU} $, then  $w_\MUU$ has no zero on $[  \max \{   \overline{s},  \widetilde{s}  \}, \infty )$. By possibly assuming that $\MUU $ is closer to $ \overline{\MUU} $, it follows from continuous dependence (Proposition~\ref{eEX1:b1}~\eqref{eEX1:6:b1}) that $\wN ( \MUU) = \wN (  \overline{\MUU} )$, which proves Property~\eqref{eEX2:1} of Proposition~\ref{eEX2}.

We next assume $\Lt (  \overline{\MUU}  )=0$. Given $\MUU \ge \MUUz$,  define  $v_\MUU (t)$ for $t >0$ by
\begin{equation} \label{eNZE3:0}
t^{\frac {1} {2}} v_\MUU (t)=s^{\lambda }w_\MUU (s), \quad t= \log s 
\end{equation} 
so that $v_ \MUU$ is a solution of~\eqref{eNZE3:3}.  
We deduce from~\eqref{eNZE3:0} and Proposition~\ref{eAB1} that
\begin{equation} \label{fTPAC1} 
t^{  \frac {1} {2}} v_\MUU (t) \goto  _{ t\to \infty  } \Ltu ( \MUU ) .
\end{equation} 
Moreover, it follows from~\eqref{fNd2d} and~\eqref{eNZE3:0} that
\begin{equation} \label{fTPAC7} 
M: = \sup_{\MUUz \le \MUU \le  \overline{\MUU}+1 }\, \sup_{t \ge 1} \, t^{ - \frac {1} {2}}  |v_\MUU (t)|  <\infty .
\end{equation} 
We deduce from~\eqref{fTPAC1} and~\eqref{eNZE3:3} that  for $t$ large, if $ v _{   \overline{\MUU}   } '(t)=0$ then $ v _{   \overline{\MUU}   } ''(t)$ and $v  _{   \overline{\MUU}   } (t)$ have the same sign. Since $v _{  \overline{\MUU}  } (t) \to 0$ as $t\to \infty $ (by~\eqref{fTPAC1}), it easily follows that  there exists $ \overline{t}_{ \overline{\MUU } } > 1 $  such that
$ v_{  \overline{\MUU} } (t) v_{  \overline{\MUU}  }'(t)<0$ for $ t\ge  \overline{t} _{  \overline{\MUU}  }$.
Thus we may assume without loss of generality that
\begin{equation}  \label{eNZE3:9} 
 v _{  \overline{\MUU }  }( t )>0 \quad  \text{and} \quad   v _{  \overline{\MUU }  } '( t)<0 \quad  \text{for} \quad t\ge  \overline{t} _{  \overline{\MUU}  }.  
\end{equation}  
In particular, $w _{  \overline{\MUU}  }$ has $\wN (  \overline{\MUU}  )$ zeros on $(0, \log  \overline{t}  _{  \overline{\MUU}  }) $, so that by continuous dependence, $w_\MUU$ also has  $\wN (  \overline{\MUU}  )$ zeros on $(0, \log  \overline{t}  _{  \overline{\MUU}  }) $ provided $ |\MUU -  \overline{\MUU} |$ is sufficiently small. 
Therefore, we need only show that, by possibly assuming $ |\MUU -  \overline{\MUU} |$ smaller, either 
\begin{equation} \label{fTPAC2} 
v_\MUU (t) >0, \quad t\ge  \overline{t}  _{  \overline{\MUU}  }
\end{equation} 
or else 
\begin{equation}  \label{fTPAC3} 
 \text{$v_\MUU $ has one zero on $[  \overline{t}_{  \overline{\MUU}  }, \infty ) $ and $\Lt ( \MUU) \not = 0$. }
\end{equation}  
By~\eqref{eNZE3:9}, we have $v_\MUU (  \overline{t} _{  \overline{\MUU}  } )>0 $ if $ | \MUU -  \overline{\MUU} |$ is sufficiently small, and we define $t_\MUU>  \overline{t} _{  \overline{\MUU}  } $ by
\begin{equation*} 
t_\MUU = \sup\{ t>  \overline{t} _{  \overline{\MUU}  } ;\,  v_\MUU > 0 \text{ on }  (  \overline{t} _{  \overline{\MUU}  }, t)\}.
\end{equation*} 
If $t_\MUU =\infty $, then we are in case~\eqref{fTPAC2}. Assume now $t_\MUU <\infty $, so that $v_\MUU (t_\MUU ) =0$ and $v_\MUU ' (t_\MUU ) < 0$, and set
\begin{equation*} 
 \widetilde{t}_\MUU  =\sup\{ t> t_\MUU ;\,  v_\MUU '<0 \text{ on }  (t_\MUU, t)\}.
\end{equation*} 
We claim that $ \widetilde{t} _\MUU =\infty$. Assuming this claim, we see that $v_\MUU$ has one zero on $[  \overline{t}, \infty ) $ and that $\limsup v_\MUU (t) <0$ as $t\to \infty $. Therefore, $t^{\frac {1} {2}}v_\MUU (t) \to -\infty $ as $t\to \infty $. By~\eqref{eNZE3:0}, this means that $s ^\lambda  w_\MUU (s) \to -\infty $ as $s\to \infty $, and Proposition~\ref{eAB1} implies that $\Lt ( \MUU )\not = 0$. Thus we see that we are in case~\eqref{fTPAC3}, and this completes the proof. 
We finally prove that $ \widetilde{t} _\MUU =\infty$. Otherwise, $v_\MUU ' ( \widetilde{t} _\MUU ) =0  \text{ and } v_\MUU '' ( \widetilde{t} _\MUU ) \ge 0$.
Since $v_\MUU  ( \widetilde{t} _\MUU ) <0$, equation~\eqref{eNZE3:3} yields
\begin{equation*} 
\frac { 1 } { \widetilde{t}_\MUU ^2}+e^{- \widetilde{t}_\MUU }\Bigl (\frac {1}  { 2 \widetilde{t}_\MUU }-\frac {1} {\alpha } - \widetilde{t}_\MUU ^{\frac {\alpha } {2} } |v_\MUU |^\alpha  \Bigr) \le 0
\end{equation*} 
so that
\begin{equation} \label{fTPAC4} 
 \widetilde{t}_\MUU ^{\frac {\alpha } {2} } |v_\MUU ( \widetilde{t}_\MUU) |^\alpha  \ge e^{ \widetilde{t}_\MUU }\frac { 1} { \widetilde{t}_\MUU ^2}+ \frac {1}  { 2 \widetilde{t}_\MUU }-\frac {1} {\alpha }  
\end{equation} 
By~\eqref{eNZE3:9} and continuous dependence, $t_\MUU$, hence $ \widetilde{t}_\MUU $ can be made arbitrarily large by assuming $ |\MUU -  \overline{\MUU} |$ sufficiently small. 
In particular, with $M>0$ defined by~\eqref{fTPAC7}, we have
\begin{equation*} 
e^{ \widetilde{t}_\MUU } \frac { 1} { \widetilde{t}_\MUU ^2}+ \frac {1}  { 2 \widetilde{t}_\MUU }-\frac {1} {\alpha } \ge    (2 M   \widetilde{t}_\MUU )^\alpha 
\end{equation*} 
if $ |\MUU -  \overline{\MUU} |$ is sufficiently small. 
Applying~\eqref{fTPAC4}, we obtain  $  \widetilde{t}_\MUU ^{ - \frac {1 } {2} } |v_\MUU ( \widetilde{t}_\MUU) |  \ge 2M$, which 
 contradicts~\eqref{fTPAC7}.  Thus $ \widetilde{t} _\MUU =\infty$, which completes the proof for $N=2$.

\section{Singular profiles} \label{sSinPro} 

In what follows, we prove Theorem~\ref{eMain2}.
We suppose $N\ge 3$ and $\frac {2} {N-2} \le \alpha < \frac {4} {N-2}$, and
we consider the collection $( w_\MUU ) _{ \MUU \ge 1 }$ of solutions of~\eqref{IPE} given by Proposition~$\ref{eEX1:b1}$.
 It follows from Proposition~\ref{eAB1}~\eqref{eAB1:1} and~\eqref{fPR3}  that the corresponding profile 
\begin{equation} \label{fTPA42} 
  \widetilde{f}_\MUU  (r)=r^{-\frac {2} {\alpha }}w_\MUU (r^{-2})
\end{equation} 
behaves like $ \Lt (\MUU) r^{- \frac {2} {\alpha }} \phi _2 ( r ^{-2}) $ as $r\to 0$. 
Thus we see that if $ \Lt (\MUU) \not = 0$ (i.e. if $w_\MUU$ has slow decay), then $  \widetilde{f}_\MUU  (r)$ is singular at $r = 0$. On the other hand, one verifies easily that $  \widetilde{f}_\MUU  \in  L^{\alpha +1}_\Loc (\R^N )$. Furthermore, it follows from~\eqref{fTPA42} that 
\begin{equation}  \label{fTPA42b2} 
r^{1+ \frac {2} {\alpha }}   \widetilde{f}_\MUU ' (r) = - \frac {2} {\alpha } w_\MUU (r^{-2}) - 2 r^{-2}  w_\MUU ' (r^{-2}); 
\end{equation} 
and so
\begin{equation} \label{fTPA42b1} 
r^{ \frac {2} {\alpha }}   (  |\widetilde{f}_\MUU  (r)| + r |\widetilde{f}_\MUU ' (r)| ) \le C \sup  _{ s\ge 0 } \, (  |w_\MUU (s)| + s  |w_\MUU ' (s)| ) \le C
\end{equation}
where we apply  Proposition~\ref{eEX1:b1}~\eqref{eEX1:1:b1} in the last inequality.

We claim that $ \widetilde{f}_\MUU $ is a solution of~\eqref{fpr2} in the sense of distributions, i.e.
\begin{equation} \label{fTPA61} 
\int_{\R^N}  \widetilde{f} _\MUU  \Bigl( \Delta \varphi - \frac {1} {2} \nabla \cdot (x \varphi ) + \frac {1} {\alpha } \varphi  + | \widetilde{f} _\MUU |^\alpha  \varphi \Bigr) \, dx =0 
\end{equation}  
for all  $\varphi \in C^\infty _\Comp (  \R^N)$. 
To see this, we let $\varepsilon >0$. Since $  \widetilde{f}_\MUU  \in  L^{\alpha +1}_\Loc (\R^N )$, we see that
\begin{equation} \label{fTPA62} 
\begin{split} 
 \Bigl| \int_{ \{  |x|<\varepsilon  \}} \widetilde{f} _\MUU    \Bigl( \Delta \varphi - \frac {1} {2} \nabla \cdot (x \varphi ) + \frac {1} {\alpha } \varphi  &+ | \widetilde{f} _\MUU |^\alpha  \varphi \Bigr) \Bigr| \\& \le C   \int_{ \{  |x|<\varepsilon  \}} ( |\widetilde{f}_\MUU  | +  | \widetilde{f} _\MUU |^{ \alpha +1})\goto _{ \varepsilon \downarrow 0 } 0 .
\end{split} 
\end{equation}  
On the other hand, $\widetilde{f}_\MUU $ is a classical solution of~\eqref{fpr2} on $  \R^N  \setminus \{ 0\} $, so that integration by parts yields
\begin{equation} \label{fTPA63} 
\begin{split} 
\int_{ \{  |x| >\varepsilon  \}} \widetilde{f} _\MUU    \Bigl( \Delta \varphi - \frac {1} {2} \nabla \cdot (x \varphi ) + \frac {1} {\alpha } \varphi  &+ | \widetilde{f} _\MUU |^\alpha  \varphi \Bigr)     \\&  = -  \int_{ \{  |x|=\varepsilon  \}}  \Bigl( \widetilde{f} _\MUU  \frac {\partial \varphi } {\partial r}- \varphi \frac {\partial  \widetilde{f} _\MUU } {\partial r} -\frac {\varepsilon } {2}  \widetilde{f}_\MUU \varphi   \Bigr)    .
\end{split} 
\end{equation} 
Therefore, \eqref{fTPA61}  follows from~\eqref{fTPA62} and~\eqref{fTPA63} provided we show that 
\begin{equation*} 
 \int_{ \{  |x|=\varepsilon  \}}  \Bigl(  |\widetilde{f} _\MUU|  + \Bigl|  \frac {\partial  \widetilde{f} _\MUU } {\partial r}  \Bigr| \Bigr)   \goto _{ \varepsilon \downarrow 0 } 0 ,
\end{equation*} 
i.e. $r^{N-1} ( |   \widetilde{f}_\MUU (r) | +  |   \widetilde{f}_\MUU ' (r)|) \to 0$ as $r \downarrow 0$.
If $\alpha > \frac {2} {N-2}$, this is a consequence of~\eqref{fTPA42b1}.  If $\alpha = \frac {2} {N-2}$, this follows from~\eqref{fTPA42b2}, \eqref{fTPA42b1}, and the  fact that  $ |w_\MUU ( s ) | + s |w_\MUU ' ( s ) |  \to 0$ as $s\to \infty $ by Proposition~\ref{eEX1:b1}~\eqref{eEX1:3:b1}. 

We now conclude the proof of Theorem~\ref{eMain2} as follows. By Lemma~\ref{eTPA15}~\eqref{eTPA15:2}, there exist an integer $ \overline{m} $ and an increasing sequence $( \overline{\MUU}_m  ) _{ m\ge  \overline{m}  } \subset [ \MUUz, \infty )$, $\overline{\MUU}  _m \to \infty $ as $m\to \infty $, such that $\wN (  \overline{\MUU} _m ) = m$ and $ \Lt (  \overline{\MUU} _m )\not = 0$. Next, we note that by Proposition~\ref{eEX1:b1}~\eqref{eEX1:3:b1} (for $\frac {2} {N-2} < \alpha < \frac {4} {N-2}$) and Lemma~\ref{eBN2b} (for $\alpha = \frac {2} {N-2}$), $ \Lt (  \overline{\MUU} _m ) = \pm \nu$ where
\begin{equation*} 
\nu =
\begin{cases} 
\beta ^{\frac {1} {\alpha }} & \frac {2} {N-2} < \alpha < \frac {4} {N-2} \\
(\frac {2} {\alpha } )^{\frac {2} {\alpha }} & \alpha = \frac {2} {N-2} .
\end{cases} 
\end{equation*} 
Since  $ w _{  \overline{\MUU}_m  }$ has $m$ zeros, we see that  $ \Lt (  \overline{\MUU} _m ) = (-1)^m \nu $. Consequently,  the profiles $h_m =  \widetilde{f}  _{  \overline{\MUU} _m }$ satisfy the first part of Theorem~\ref{eMain2}. 

Since $ |h_m (r)| + r | h_m '(r) | \le C r^{-\frac {2} {\alpha }}$ by~\eqref{fTPA42b1} and, as noted previously, $ h_m =  \widetilde{f} _{  \overline{\MUU} _m }  \in  L^{\alpha +1}_\Loc (\R^N )$, the second part of Theorem~\ref{eMain2} is an immediate consequence of the following lemma, with $f= h_m$ and $\omega (x) \equiv  \overline{\MUU}_m $.

\begin{lemma} \label{eSP1} 
Let $\alpha \ge \frac {2} {N-2}$, and let $f\in C^2( \R^N \setminus \{0\} ) \cap L^{\alpha +1}_\Loc (\R^N )$ be a solution of~\eqref{fpr2} in the sense of distributions and satisfy
\begin{equation} \label{eSP1:1} 
 |f(x)| +  | x \cdot \nabla f(x)|\le C  |x|^{-\frac {2} {\alpha }} \quad x\not = 0 ,
\end{equation} 
and $ |x|^\frac {2} {\alpha } f(x)- \omega (x) \to 0$ as $ |x| \to \infty $, for some $\omega \in C (\R^N \setminus \{0\})$ homogeneous of degree $0$. 
 If $u \in C^2( (0,\infty )\times (\R^N \setminus \{0\})) $ is defined by~\eqref{fpr1} for $t>0$ and $x\not = 0 $, then  
 \begin{equation}  \label{eSP1:3:1} 
 | u | + t  |u _t | +  |x\cdot \nabla u| \le C  |x|^{- \frac {2} {\alpha }} 
\end{equation}
and
\begin{equation}  \label{eSP1:3:2} 
u, u_t , x\cdot \nabla u\in  C((0,\infty ), L^p (\R^N ) + L^q (\R^N ) ) \text{ for } 1\le p<  \frac {N\alpha } {2}< q .
\end{equation}  
Furthermore,
\begin{equation} \label{eSP1:4} 
 \Delta u,  |u|^\alpha u\in C((0,\infty ), L^p (\R^N ) + L^q (\R^N ) ) 
\end{equation} 
and $u$ is a solution of~\eqref{NLHE}, and also of~\eqref{fNZ2} with $\DI (x)= \omega (x)  |x|^{-\frac {2} {\alpha }}$,  in $C((0,\infty ), L^p (\R^N ) + L^q (\R^N ) ) $ where $p, q$ satisfy~\eqref{eSP1:5}.
Moreover,  $u$ satisfies~\eqref{fEP2} and~\eqref{fEP3}.
 In addition, the map $t \mapsto u (t) - e^{t \Delta } \DI $ is in $ C([0,\infty ), L^r (\R^N ) )$ for all   $\frac {N\alpha } {2(\alpha +1) } < r < \frac {N\alpha } {2}$ if $\alpha > \frac {2} {N-2}$, and  in $ C([0,\infty ), L^1 (\R^N ) + L^r (\R^N ) )$ for all  $r>1$ if $\alpha = \frac {2} {N-2}$.
\end{lemma} 

\begin{proof} 
The proof is similar to that of Lemma~\ref{eQSol}, but some extra care is needed because of the possible singularity of $f$ at $x=0$.
Differentiating~\eqref{fpr1} with respect to $t$, we see that $u$ satisfies equation~\eqref{feQSol4} for $t>0$ and $x\not = 0$. Property~\eqref{eSP1:3:1} follows from \eqref{fpr1} and~\eqref{eSP1:1}, and property~\eqref{eSP1:3:2} is then a consequence of the dominated convergence theorem. 
Next, we deduce from estimate~\eqref{eSP1:1} if $\alpha >\frac {2} {N-2}$, and the assumption $f \in L^{\alpha +1}_\Loc (\R^N )$ if $\alpha =\frac {2} {N-2}$ that $ |u|^\alpha u\in C((0,\infty ), L^p (\R^N ) + L^q (\R^N ) ) $ whenever $p,q$ satisfy~\eqref{eSP1:5}. In addition, it follows from~\eqref{fpr2} that for every $t>0$,
\begin{equation} \label{fEId2} 
\Delta u + \frac {1} {2t} x\cdot \nabla u + \frac {1} {\alpha t} u +  |u|^\alpha u=0
\end{equation} 
in the sense of distributions. In particular, we see that $ \Delta u\in C((0,\infty ), L^p (\R^N ) + L^q (\R^N ) ) $ whenever $p,q$ satisfy~\eqref{eSP1:5}. Moreover, \eqref{fEId2} and~\eqref{feQSol4} imply that  $u$ is a solution of~\eqref{NLHE}. Properties~\eqref{fEP2} and~\eqref{fEP3} follow from the fact that $u (t, x) \to \DI  (x) $ for all $x\not = 0$ by~\eqref{fTR1}, estimate~\eqref{eSP1:3:1}, and dominated convergence.  

We now show that $u$ satisfies~\eqref{feQSol1} for all $0< \varepsilon <t $. To see this, we consider $\varphi \in C^\infty _\Comp (\R^N )$ such that $\varphi (x)=1$ for $ |x| \le 1$. 
Note that $ \Delta (\varphi u) = \varphi \Delta u + 2\nabla u\cdot \nabla \varphi + u\Delta \varphi $, so that
\begin{equation}  \label{fEId3} 
(\varphi u)_t = \Delta (\varphi u) - 2\nabla u\cdot \nabla \varphi - u\Delta \varphi + \varphi  |u|^\alpha u.
\end{equation}  
Therefore, it follows from~\eqref{eSP1:3:2}  that $\varphi u\in C^1 ((0,\infty ), L^p (\R^N ) )$;  and from~\eqref{eSP1:4} and~\eqref{eSP1:3:2}  that $\Delta (\varphi u), \nabla u\cdot \nabla  \varphi , \varphi  |u|^\alpha u \in C ((0,\infty ), L^p (\R^N ) )$, for $1< p< \frac {N\alpha } {2(\alpha +1)}$ if $\alpha >\frac {2} {N-2}$ and $p=1$ if $\alpha =\frac {2} {N-2}$. In the case $\alpha =\frac {2} {N-2}$, note also that $\nabla (\varphi u)\in C((0,\infty ), L^1 (\R^N ) )$ by~\eqref{eSP1:3:1} and dominated convergence. 
Therefore $\varphi u$ has the required regularity so that~\eqref{fEId2} implies
\begin{equation} \label{fEId4} 
\varphi u(t)= e^{ (t- \varepsilon ) \Delta } \varphi u(\varepsilon ) + \int _\varepsilon ^t e^{t-s) \Delta } [ - 2\nabla u\cdot \nabla \varphi - u\Delta \varphi + \varphi  |u|^\alpha u ]
\end{equation} 
in $L^p (\R^N ) $ for all $0<\varepsilon <t$. Next, we have
\begin{equation}  \label{fEId5} 
((1- \varphi) u)_t = \Delta ((1- \varphi ) u) + 2\nabla u\cdot \nabla \varphi + u\Delta \varphi +(1- \varphi ) |u|^\alpha u.
\end{equation}  
Arguing as above, one sees that $(1-\varphi )u$ has the required regularity so that
\begin{equation} \label{fEId6} 
(1-\varphi ) u(t)= e^{(t- \varepsilon ) \Delta } (1-\varphi ) u(\varepsilon ) + \int _\varepsilon ^t e^{t-s) \Delta } [  2\nabla u\cdot \nabla \varphi +u\Delta \varphi + (1- \varphi ) |u|^\alpha u ]
\end{equation} 
in $L^q (\R^N ) $ for all $q >\frac {N \alpha } {2} $ and $0<\varepsilon <t$. Summing up~\eqref{fEId4} and~\eqref{fEId6}, we see that $u$ satisfies~\eqref{feQSol1} in $ L^p (\R^N ) + L^q (\R^N )  $ for all $p, q$ satisfying~\eqref{eSP1:5} and all $0< \varepsilon <t $. 

We now prove the last statement in the lemma, which will imply equation~\eqref{fNZ2}.
For this, we consider separately the cases $\alpha >\frac {2} {N-2}$ and $\alpha =\frac {2} {N-2}$.
Suppose first $\alpha >\frac {2} {N-2}$ and let $\frac {N\alpha } {2(\alpha +1) } < r < \frac {N\alpha } {2}$.
Since $  | \cdot |^{- \frac {2(\alpha +1)} {\alpha }} \in L^1 (\R^N ) +L^r (\R^N ) $, we have $e^{ \Delta }   | \cdot |^{- \frac {2(\alpha +1)} {\alpha }} \in L^r (\R^N ) $, so that by scaling
\begin{equation*} 
 \| e^{(t-s) \Delta }   | \cdot |^{- \frac {2(\alpha +1)} {\alpha }} \| _{ L^r } =   (t-s) ^{- \frac {\alpha +1} {\alpha } + \frac {N} {2r} }   \| e^{ \Delta }   | \cdot |^{- \frac {2(\alpha +1)} {\alpha }} \| _{ L^r }.
\end{equation*} 
Therefore, since $ |f(x)|\le C  |x|^{-\frac {2} {\alpha }}$, we have
\begin{equation*} 
\| e^{(t-s) \Delta }  |u(s)|^\alpha u(s) \| _{ L^r }\le C  \| e^{(t-s) \Delta }   | \cdot |^{- \frac {2(\alpha +1)} {\alpha }} \| _{ L^r } = C  (t-s) ^{- \frac {\alpha +1} {\alpha } + \frac {N} {2r} }.
\end{equation*} 
Since $ - \frac {2(\alpha +1)} {\alpha } >-1$, we see that the integral term in~\eqref{fNZ2} is continuous $[0, \infty ) \to L^r (\R^N ) $, and that we can let $\varepsilon \downarrow 0$ in the  integral term in~\eqref{fNZ2}. Since $e^{(t-\varepsilon ) \Delta } u(\varepsilon ) \to e^{t\Delta } \DI$ as $\varepsilon \downarrow 0$ in $L^q (\R^N )$ for all $q>\frac {N\alpha } {2}$, this completes the proof in the case $\alpha >\frac {2} {N-2}$.

Finally, in the case $\alpha = \frac {2} {N-2}$, recall that $f\in L^{ \alpha +1} _\Loc (\R^N ) $.
We fix $T>0$ and write for $0< t <T$
\begin{equation*} 
 |u|^\alpha u = 1 _{ \{  |x|<\sqrt t \} }  |u|^\alpha u + 1 _{ \{  \sqrt t <  |x| < \sqrt T\} }  |u|^\alpha u  + 1 _{ \{  |x|>\sqrt T \} }  |u|^\alpha u  = : F_1 + F_2 + F_3.
\end{equation*} 
We have by~\eqref{fpr1} 
\begin{equation} 
\begin{split} 
\| e^{(t-s) \Delta }  F_1 (s) \| _{ L^1 } & \le \|  F_1 (s) \| _{ L^1 }= \int  _{  \{  |x|<\sqrt s \} } s^{-\frac {\alpha +1} {\alpha }}   \Bigl| f \Bigl( \frac {x} {\sqrt s} \Bigr) \Bigr|^{\alpha +1} dx \\ & =  \int  _{  \{  |x|< 1 \} }  |f(x)|^{\alpha +1}dx
\end{split} 
\end{equation} 
since $ \frac {\alpha +1} {\alpha } = \frac {N} {2}$. Moreover, since $ |u ( t, x)|^{\alpha +1} \le C  |x|^{-\frac {2(\alpha +1)} {\alpha }} = C  |x|^{-N }$,
\begin{equation} 
\| e^{(t-s) \Delta }  F_2 (s) \| _{ L^1 }  \le  \|  F_2 (s) \| _{ L^1 }   \le   \int  _{  \{ \sqrt s < |x| < \sqrt T \} }  |x| ^{ - N } dx    \le C  \Bigl|  \log \frac {s} {T}   \Bigr| .
\end{equation} 
Furthermore, given $r>1$,
\begin{equation} 
\begin{split} 
\| e^{(t-s) \Delta }  F_3 (s) \| _{ L^r } & \le  \|  F_2 (s) \| _{ L^r }  \\ & \le    \Bigl(  \int  _{  \{  |x| > \sqrt T \} }  |x| ^{ - rN } dx \Bigr)^{\frac {1} {p}} \le C<\infty .
\end{split} 
\end{equation} 
because $r> \frac {N\alpha } {2(\alpha +1)}$. Therefore,
\begin{equation*} 
\| e^{(t-s) \Delta }  |u(s)|^\alpha u(s) \| _{ L^1 + L^r }   \le C (1 +  |\log s| )
\end{equation*} 
and one can easily complete the proof as in the case $\alpha >\frac {2} {N-2}$.
\end{proof} 

\begin{remark} 
In the case $\alpha > \frac {2} {N-2}$, the singular stationary solution $u$ of~\eqref{NLHE}  given by~\eqref{statsingsol} is in particular a self-similar  solution. Its profile (which is $u$ itself) 
satisfies the assumptions of Lemma~\ref{eSP1}, so that  $u$ is a solution of~\eqref{fNZ2}.
Note that $u$ is time-independent, but the other two terms in~\eqref{fNZ2} do depend on time.
\end{remark}

\begin{remark} \label{eRemz} 
If $\alpha  <\frac {2} {N-2}$ and $\Lt ( \MUU)\not = 0$, then the singularity of $ \widetilde{f}_\MUU  (r)$ at $0$ is of order $r^{-(N-2)}$. Therefore $\Delta  \widetilde{f} _\MUU$ has a Dirac mass at the origin, so that $ \widetilde{f} _\MUU$ is not a solution of~\eqref{fpr2}  in the sense of distributions. 
\end{remark}

\section{The case of dimension 1} \label{sDim1} 
The purpose of this section is to prove Theorem~\ref{eMain3}.
Its proof is similar to that of Theorem~\ref{eMain1}, with two major differences. The first one is that the proofs of the results analogous to Propositions~\ref{eAB1} and~\ref{eEX2} are completely elementary. This is due to the fact that equation~\eqref{fpr3} is not singular at $r=0$.  (It seems that there is no simplification in the proof of Proposition~\ref{fFIn1} when $N=1$.) The second major difference is that Theorem~\ref{eMain3} concerns both radially symmetric (i.e. even) profiles, and odd profiles. 

Throughout this section, we suppose $N=1$.
We first observe that if $f\in C^2([0,\infty ))$ is a solution of~\eqref{fpr3} on $(0,\infty )$ and if $f'(0)=0$, then extending $f$  by setting $f(x)= f(-x)$ for $x<0$ yields a solution of~\eqref{fpr2} on $\R$. Similarly, if $f (0)=0$, then extending $f$  by setting $f(x)= -f(-x)$ for $x<0$  also yields a solution of~\eqref{fpr2}. Therefore, we need only construct solutions of~\eqref{fpr3} on $(0,\infty )$ that satisfy either $f '(0)= 0$ (corresponding to an even profile) or $f (0)= 0$ (corresponding to an odd profile).

We consider the collection $( w_\MUU ) _{ \MUU \ge 1 }$ of solutions of~\eqref{IPE} given by Proposition~$\ref{eEX1:b1}$. We recall that, by Proposition~\ref{eEX1:b1}~\eqref{eEX1:6:b1},
\begin{equation} \label{fTPA27b1} 
 \text{the map $\MUU \mapsto  w_\MUU  $ is continuous $[\MUUz, \infty ) \to C^1([0, 1])$.} 
\end{equation} 
Next, we define
\begin{equation} \label{fTPA25} 
 \widetilde{f}_\MUU (r)=   r^{-\frac {2} {\alpha }} w _{ \MUU  }  (r^{-2})
\end{equation} 
for $r>0$. It follows that $ \widetilde{f}_\MUU $ is a solution of the profile equation~\eqref{fpr3} and that
\begin{equation} \label{fTPA26} 
 \widetilde{f}_\MUU (1)=    w _{ \MUU  }  (1), \quad  \widetilde{f}_\MUU ' (1)=  -2  w _{ \MUU  } ' (1) - \frac {2} {\alpha } w _{ \MUU  }  (1) .
\end{equation} 
Moreover, since $N=1$, equation~\eqref{fpr3} is not singular at $r=0$, so that $ \widetilde{f}_\MUU $ can be extended to a solution of~\eqref{fpr3} for all $r\in \R$. (This follows from an obvious energy argument.)
One deduces easily by using~\eqref{fTPA27b1} and~\eqref{fTPA26} that 
\begin{equation} \label{fTPA27} 
 \text{the map $\MUU \mapsto  \widetilde{f}_\MUU  $ is continuous $[\MUUz, \infty ) \to C^1([0, 1])$.} 
\end{equation} 
Thus we see that the study of $w_\MUU $ on $[0,\infty )$ is equivalent to the study of $w_\MUU$ on $[0,1]$ and the study of $ \widetilde{f}_\MUU $ on $[0, 1]$, The problem is therefore reduced to  compact intervals, which considerably simplifies the analysis.

Since $N=1$, we have $\lambda _1=\frac {1} {\alpha }$, $\lambda _2=\frac {1} {\alpha }+\frac {1} {2}$,  $\beta<0$, and $\lambda _2>\lambda _1 >0$. In particular, unlike in the case $N\ge 3$, $s^{-\lambda _1}$ decays more slowly than $s^{-\lambda _2}$, and also represents the generic behavior as $s\to \infty $ of solutions of~\eqref{IPE}. Consequently, in this section we define 
\begin{equation} \label{fL2N1} 
\Lt ( \MUU) = \lim  _{ s\to \infty  } s^{\lambda _1} w_\MUU (s) = \lim  _{ s\to \infty  } s^{\frac {1} {\alpha }} w_\MUU (s) 
\end{equation} 
and 
\begin{equation} \label{fL1N1} 
\Ltu ( \MUU) = \lim  _{ s\to \infty  } s^{\lambda _2} w_\MUU (s) = \lim  _{ s\to \infty  } s^{\frac {1} {\alpha }+ \frac {1} {2}} w_\MUU (s) 
\end{equation} 
whenever these limits exist.

\begin{remark} \label{eTPA3} 
The following properties are simple consequences of the above observation.
\begin{enumerate}[{\rm (i)}] 

\item  \label{eTPA3:1} 
Since both $w_\MUU$ and $ \widetilde{f}_\MUU $ have only a finite number of zeros on the compact interval $[0,1]$, it follows from~\eqref{fTPA25} that $w_\MUU$ has at most a finite number of zeros on $[0,\infty )$. This gives a quick proof of Proposition~\ref{eEX1:b1}~\eqref{eEX1:2:6}  in the case $N=1$. 
Recall that the possible zero of $ \widetilde{f} _\MUU$ at $r=0$ is not counted in $\wN ( \MUU )$, where $\wN ( \MUU )$ is given by~\eqref{fFIn1:1}.

\item  \label{eTPA3:2} 
Formula~\eqref{fTPA25} implies
\begin{equation*} 
 \widetilde{f}_\MUU (0) = \lim  _{ r\downarrow 0 } \widetilde{f}_\MUU (r)  = \lim  _{ s\to \infty  } s^{\frac {1} {\alpha }} w_\MUU (s) .
\end{equation*} 
Therefore, the limit~\eqref{fL2N1} exists and 
\begin{equation}\label{fTPA27:10} 
\Lt ( \MUU) =  \widetilde{f}_\MUU  (0) .
\end{equation} 
Since $ \widetilde{f}_\MUU (0) $ depends continuously on $\MUU$, $\Lt $ is continuous $[\MUUz, \infty ) \to \R$. 
In addition, 
\begin{equation*}  
 \widetilde{f}_\MUU ' (0) = \lim  _{ r\downarrow 0 } \frac {\widetilde{f}_\MUU (r) - \widetilde{f}_\MUU (0)} {r} = \lim  _{ s\to \infty }  s^{\frac {1} {2}} (s^{\frac {1} {\alpha }} w_\MUU (s) -\Lt (\MUU) ).
\end{equation*} 
Thus we see that if $\Lt (\MUU) =0$, then  the limit~\eqref{fL1N1} exists and 
\begin{equation} \label{fTPA27:20} 
\Ltu ( \MUU) = \widetilde{f}_\MUU ' (0) .
\end{equation} 
In particular, $\Ltu ( \MUU) \not = 0$ for otherwise we would have $\widetilde{f}_\MUU  (0) = \widetilde{f}_\MUU ' (0) =0$, hence $ \widetilde{f} _\MUU \equiv 0$. 

\item  \label{eTPA3:3} 
It follows from~\eqref{fTPA25} that
\begin{equation*} 
\lim _{ r\to \infty  } r^{\frac {2} {\alpha }} \widetilde{f}_\MUU (r) = w _\MUU (0) = \MUU \ge \MUUz
\end{equation*} 
so that no zero of $ \widetilde{f}_\MUU $ can appear or disappear at infinity by varying $\MUU$.
Moreover if $r>0$, then $\widetilde{f}_\MUU ' (r)\not = 0$ whenever $\widetilde{f}_\MUU (r) =0$, so that no zero of $ \widetilde{f}_\MUU $ on $( 0, \infty ) $ can appear or disappear  by varying $\MUU$.
Thus we see that $\wN ( \MUU)$ can only change at a $\MUU$ for which $ \widetilde{f}_\MUU  (0)  = 0$, i.e. $\Lt (  \MUU ) = 0$. In particular, if $\MUUz \le \MUU_1 < \MUU_2$, $\wN ( \MUU_1) \not = \wN (\MUU _2)$ and $ \widetilde{f} _{ \MUU_1 } (0)  \widetilde{f} _{ \MUU_2 } (0) \not = 0 $, then there exists $\MUU\in (\MUU_1, \MUU_2)$ such that $ \widetilde{f}_\MUU (0) =0 $.

\item  \label{eTPA3:4} 
It follows from Property~\eqref{eTPA3:3} above that if $\Lt (  \overline{\MUU} ) \not = 0$, then $\wN ( \MUU) = \wN (  \overline{\MUU} )$ for $\MUU $ close to $\ \overline{ \MUU} $. 

\item  \label{eTPA3:5} 
Suppose $\Lt (  \overline{\MUU} ) = 0$ (i.e. $ \widetilde{f} _{  \overline{\MUU}  } (0)= 0$), so that  $  \widetilde{f}_{\overline{\MUU}} ' (0) \not =  0$. Without loss of generality we may suppose $  \widetilde{f}_{\overline{\MUU}} ' (0) >  0$. It follows from~\eqref{fTPA27} that there exist $\delta , \rho >0$ such that if $ |\MUU -  \overline{\MUU} |\le \delta $, then
\begin{gather} 
\frac {1} {2}   \widetilde{f}_{ \overline{\MUU } } ' (0) \le   \widetilde{f} _\MUU '(r) \le 2   \widetilde{f}_{ \overline{\MUU } } ' (0) , \quad 0\le r\le \rho , \label{eTPA3:5:1}  \\
|   \widetilde{f} _\MUU (0) |\le \frac {\rho } {4}  \widetilde{f}_{ \overline{\MUU } } ' (0)  . \label{eTPA3:5:2}
\end{gather} 
Therefore, $   \widetilde{f}_ \MUU$ is  increasing in $[0,\rho ]$ and
\begin{equation} \label{eTPA3:5:3}
    \widetilde{f} _\MUU (r) \ge     \widetilde{f} _\MUU (0) + \frac {r} {2}   \widetilde{f}_{ \overline{\MUU } } ' (0) \quad 0\le r\le \rho .
\end{equation} 
It follows in particular from~\eqref{eTPA3:5:2} and~\eqref{eTPA3:5:3} that
\begin{equation}  \label{eTPA3:5:4}
  \widetilde{f} _\MUU (\rho ) \ge   \frac {\rho } {4}   \widetilde{f}_{ \overline{\MUU } } ' (0) >0.
\end{equation} 
If $  \widetilde{f} _\MUU (0) \ge 0$, then $  \widetilde{f} _\MUU (r) > 0$ on $(0, \rho ]$
by~\eqref{eTPA3:5:3}, so that $\wN ( \MUU) = \wN (  \overline{\MUU} ) $. If $  \widetilde{f} _\MUU (0) < 0$, then $  \widetilde{f} _\MUU $ has exactly one zero on $(0, \rho ]$ by~\eqref{eTPA3:5:4} and the fact that $ \widetilde{f}_\MUU $ is increasing, so that $\wN ( \MUU) = \wN (  \overline{\MUU} ) +1$. 

\item  \label{eTPA3:6} 
It follows from Properties~\eqref{eTPA3:4} and~\eqref{eTPA3:5} above that Proposition~\ref{eEX2} holds as well in the case $N=1$, where $\Lt ( \MUU)$ and $\Ltu ( \MUU)$  are now given by~\eqref{fL2N1} and~\eqref{fL1N1}. 
\end{enumerate} 
\end{remark} 

In order to prepare the proof of Theorem~\ref{eMain3}, we make the following observation.

\begin{lemma} \label{eLFN1} 
Let $\wN ( \MUU) $ be as defined by~\eqref{fFIn1:1}. Set $ \overline{m} =  \wN  (2) +1$,
\begin{equation} \label{eNZn1:9}
E_m= \{ \MUU \ge 2 ;\,  \wN ( \MUU )\ge  m+1  \}
\end{equation} 
and
\begin{equation} \label{eNZn1:10}
\MUU _m = \inf  E _m \in [2, \infty ) 
\end{equation} 
for $m\ge  \overline{m} $.
It follows that $E_m \not = \emptyset$ and $2\le \MUU _m <\infty $ for all $m\ge  \overline{m} $.
Moreover,  the sequence $(\MUU _m) _{ m\ge  \overline{m}  }$ is increasing and satisfies
\begin{equation} \label{fTPA90} 
\MUU  _m \goto  _{ m\to \infty  } \infty, \quad \wN ( \MUU _m) =m, \quad \Lt( \MUU  _m)=0 .  
\end{equation} 
\end{lemma} 

\begin{proof} 
 It follows from Proposition~\ref{fFIn1} that $E_m \not = \emptyset$, so that $2\le \MUU _m <\infty $. 
 The remaining properties follow from the argument used  in the proof of Lemma~\ref{eTPA15} (using Remark~\ref{eTPA3}~\eqref{eTPA3:6}  where the proof of Lemma~\ref{eTPA15} uses Proposition~\ref{eEX2}).
\end{proof} 

\begin{proof} [Proof of Theorem~$\ref{eMain3}$~\eqref{eMain3:1}]
We consider the sequence $(\MUU _m) _{ m\ge  \overline{m}  }$ given by Lemma~\ref{eLFN1}.
We  construct an increasing sequence $ (\widetilde{ \MUU}_m)  _{ m\ge  \overline{m} +1 } \subset ( \MUUz , \infty ) $ such that, with the notation~\eqref{fTPA25}
\begin{equation} \label{fTPA90:1} 
 \widetilde{\MUU}   _m \goto  _{ m\to \infty  } \infty,  \quad 
  \wN (  \widetilde{\MUU } _m) =m  ,   \quad 
  \widetilde{f}  _{  \widetilde{\MUU}  _m } ( 0) \not = 0  ,  \quad 
  \widetilde{f}  _{  \widetilde{\MUU}  _m } ' ( 0)  = 0.
\end{equation} 
Since  $w _{ \MUU _m } (0) =   \MUU_m >0$ and $ \wN  (\MUU _m) = m$, we deduce that $(-1)^m w _{ \MUU _m }(s) >0$ for $s$ large. Since $\Ltu ( \MUU _m )\not = 0$ is well defined (because $\Lt ( \MUU _m )=0$), we conclude that $(-1)^m \Ltu ( \MUU _m ) > 0$. 
Therefore,  it follows from~\eqref{fTPA27:10} and~\eqref{fTPA27:20}  that
\begin{gather} 
 \widetilde{f}  _{ \MUU _m } ( 0) =0 , \label{fTPA17} \\
(-1)^m \widetilde{f}  _{ \MUU _m } ' ( 0)> 0 . \label{fTPA27:2} 
\end{gather} 
We now consider $m\ge  \overline{m} +1$. Note that for $1\le \MUU < \MUU_m$, we have $\wN (\MUU)\le m$ by~\eqref{eNZn1:9}-\eqref{eNZn1:10}. 
Since $\widetilde{f}  _{ \MUU _{m-1} } ' ( 0)  \widetilde{f}  _{ \MUU _m } ' ( 0)<0$ by~\eqref{fTPA27:2},  the map $\MUU \mapsto  \widetilde{f}_\MUU '(0) $ has at least one zero on $(\MUU _{ m-1 }, \MUU_m)$. We denote by $ \widetilde{\MUU} _m$ the largest such zero, and it follows that
\begin{gather} 
\wN ( \MUU) \le m,\quad \widetilde{\MUU}_m \le  \MUU \le \MUU_m \label{ffC1}  \\
(-1)^m \widetilde{f}  _{ \MUU  } ' ( 0)> 0, \quad  \widetilde{\MUU}_m < \MUU \le \MUU_m \label{ffC2}  \\
\widetilde{f}  _{  \widetilde{\MUU}_m   } ' ( 0)= 0 .  \label{ffC3} 
\end{gather} 
We claim that
\begin{equation} \label{ffC4} 
\wN ( \widetilde{\MUU}_m   ) = m .
\end{equation} 
Assuming~\eqref{ffC4}, the conclusion easily follows. Indeed, $\MUU _{ m-1 }<  \widetilde{\MUU}_m < \MUU_m $, so that $ (\widetilde{ \MUU}_m)  _{ m\ge  \overline{m} +1 }$ is increasing and $ \widetilde{\MUU} _m \to \infty $. Moreover, since $ \widetilde{f}  _{  \widetilde{\MUU}_m   } \not \equiv 0$,  \eqref{ffC3} implies that $ \widetilde{f}  _{  \widetilde{\MUU}_m   }  ( 0)\not = 0$. Together with~\eqref{ffC4}, this shows the the sequence $ (\widetilde{ \MUU}_m)  _{ m\ge  \overline{m} +1 } $ has the desired properties. 

We now prove~\eqref{ffC4}. 
It follows from Remark~\ref{eTPA3}~\eqref{eTPA3:5} that there exists $\varepsilon >0$ such that for every $\MUU_m- \varepsilon \le \MUU \le \MUU_m$, we have either $(-1)^m  \widetilde{f}_\MUU (0) < 0 $ and $\wN (\MUU) =m+1$ or $(-1)^m  \widetilde{f}_\MUU (0) \ge  0 $ and $\wN (\MUU) =m$. Applying~\eqref{ffC1}, we deduce that $(-1)^m  \widetilde{f}_\MUU (0) \ge  0 $ and $\wN (\MUU) =m$ for $\MUU_m- \varepsilon \le \MUU \le \MUU_m$. We now decrease $\MUU$, and we note that, as long as $(-1)^m  \widetilde{f}_\MUU ' (0) >  0 $, we may keep applying Remark~\ref{eTPA3}~\eqref{eTPA3:5} (if $  \widetilde{f}_\MUU  (0) =  0 $) or Remark~\ref{eTPA3}~\eqref{eTPA3:4} (if $(-1)^m  \widetilde{f}_\MUU  (0) >  0 $), so that $(-1)^m  \widetilde{f}_\MUU (0) \ge  0 $ and $\wN (\MUU)= m$. 
It follows that $(-1)^m  \widetilde{f}_\MUU ' (0) >  0 $,  $(-1)^m  \widetilde{f}_\MUU (0) \ge  0 $ and $\wN (\MUU)= m$ for all $ \widetilde{\MUU}_m < \MUU \le \MUU_m $. 
Property~\eqref{ffC4} now follows from Remark~\ref{eTPA3}~\eqref{eTPA3:4}.

We next show that there exists a map $m : (0,\infty ) \to \N$  such that, given any $\MUU >0$, there exist  four sequences $(a _{ \MUU , m }^\pm ) _{ m \ge  m _ \MUU  }\subset (0,\infty )$ and $(b_{ \MUU , m }^\pm ) _{ m \ge  m_\MUU   }\subset (-\infty , 0)$ satisfying Properties~\eqref{eMT2:1}--\eqref{eMT2:4} of Theorem~\ref{eMT2}. 
This is done exactly as in the proof of Theorem~\ref{eMT2}, starting with the sequence $(  \widetilde{\MUU}_m)  _{ m\ge  \overline{m} +1 } $ defined above instead of the sequence $( \MUU _m) _{ m\ge  \overline{m}  }$ of Lemma~\ref{eTPA15}. 

Finally, we may assume $\MUU >0$ without loss of generality, and we see that the profiles $f= f_{a^\pm} $ with $a^\pm = a^\pm  _{ \MUU, \frac {m} {2} }$ if $m\ge m_\MUU $ is even and $a^\pm =  b^\pm  _{ \MUU, \frac {m-1} {2} }$ if $m\ge m_\MUU$ is odd, are two different radially symmetric solutions of~\eqref{fpr2} with $m$ zeros on $[ 0, \infty )$. Moreover, $r^{\frac {2} {\alpha }} f_{a^\pm} (r) \to \MUU$ as $r\to \infty $, and it follows from~\cite[Proposition~3.1]{HarauxW} that $ |f_{a^\pm}  (r)| + r  |f'_{a^\pm}  (r)|\le C( 1+ r^2)^{-\frac {1} {\alpha }}$. Theorem~\ref{eMain3}~\eqref{eMain3:1} is now an immediate consequence of Lemma~\ref{eQSol}, where $f (x) =f_{a^\pm} ( |x|)$ and $\omega (x) \equiv \MUU$. 
\end{proof} 

Part~\eqref{eMain3:2} of Theorem~\ref{eMain3} concerns odd solutions of~\eqref{IPE}. Therefore, given any $b\in \R$, we consider the solution $g_b$ of 
\begin{equation} \label{fTPA60} 
\begin{cases} 
\displaystyle g_b '' + \frac {r} {2} g_b ' + \frac {1} {\alpha } g_b +  |g_b |^\alpha g_b =0 \quad r\ge 0\\
g_b  (0)= 0, \quad g'_b  (0)= b .
\end{cases} 
\end{equation} 
The solutions $g_b$ of~\eqref{fTPA60} have properties similar to the solutions $f_a$ of~\eqref{fpr3}--\eqref{fpr3:1} (with $N=1$). We summarize some of these properties in the following proposition.

\begin{proposition} \label{eTPA4} 
Problem~\eqref{fTPA60} is globally well posed,
\begin{equation}  \label{fTPA1b61}
\sup  _{ r\ge 0 } \, (1+ r^2) ^{\frac {1} {\alpha }} (  | g_b (r) | + r | g_b '  (r) | ) < \infty ,
\end{equation} 
 and the limit
\begin{equation} \label{fTPA161}
L _1 (b) = \lim_{r\to\infty}r^{\frac{2}{\alpha}}g_b(r) \in \R
\end{equation}
exists, and is a continuous function of $b \in \R$.  Moreover, if $b \neq 0$, then $g_b$ has at most finitely many zeros on $(0,\infty )$, and we set
\begin{equation} \label{fTPA162} 
N_1(b )  = \text{the number of zeros of the function $g_b$ on $(0,\infty )$} .
\end{equation}
If $L_1 (b) \not = 0$ for some $b\not = 0$, then $N_1 $ is constant in some neighborhood of $b$.
\end{proposition} 

\begin{proof} 
Let  $b\in \R$ and let $g_b$ be the solution of~\eqref{fTPA60}. It follows from standard energy arguments that $g_b$ exists globally. We set
\begin{equation} \label{fPLE1} 
z_b (s) = s^{-\frac{1} {\alpha}} g_b \Bigl(  \frac{1} {\sqrt{s}}  \Bigr),\quad 0<s< \infty 
\end{equation}
so that $z_b$ is a solution of~\eqref{IPE} on $(0,\infty )$. It follows from~\cite[Proposition~2.4]{SoupletW} that $z_b \in C^1([0,\infty ))$; and so $z_b (0) = \lim  _{ s\to 0 } z_b (s)$ exists and is finite. Applying~\eqref{fPLE1}, we obtain the existence and finiteness of the limit~\eqref{fTPA161} with $L _1 (b) = z_b (0)$. 
Moreover, $ |z_b (s) | + s |z_b '  (s) |$ is bounded on $[0,1]$, so that by~\eqref{fPLE1},  $(1+ r^2) ^{\frac {1} {\alpha }} (  | g_b (r) | + r | g_b '  (r) | )$ is bounded on $[1,\infty )$. Since $(1+ r^2) ^{\frac {1} {\alpha }} (  | g_b (r) | + r | g_b '  (r) | )$ is clearly bounded on $[0, 1]$, we obtain~\eqref{fTPA1b61}. 

The continuous dependence of $L _1 (b) = z_b (0)$ on $b$ follows from arguments in~\cite{SoupletW}. More precisely, let
\begin{equation} \label{fPLE0} 
\sigma = \frac {1} {4(\gamma -1)},\quad  \tau = \frac {1} {\sqrt \sigma }= 2 \sqrt{\gamma -1}
\end{equation} 
and note that, given any $\tau <T<\infty $, the map $ b\mapsto g_b $ is continuous $ \R\to C^1([\tau , T])$. Applying~\eqref{fPLE1}, this implies that, given any $0<\varepsilon <\sigma $,
\begin{equation} \label{fPLE1b} 
 \text{the map }b\mapsto z_b  \text{ is continuous } \R\to C^1( [ \varepsilon , \sigma  ]).
\end{equation} 
Moreover, the map $s\mapsto 2s^2 z_b '(s)^2 + G (z_b (s))$ is nondecreasing on $(0,\sigma )$ by~\eqref{fENE1:b2} and~\eqref{fENE3:b1}, so that $z_b$ is bounded on $(0,\sigma )$ in terms of $z_b(\sigma )$ and $z_b ' (\sigma )$. Therefore, by~\eqref{fPLE1}, 
\begin{equation}  \label{fPLE2} 
\sup _{ 0<s<\sigma  }  |z_b (s)| \le A ( g_b( \tau ), g_b ' (\tau ))
\end{equation} 
where $A$ is a continuous function of its arguments. 
Let $b\in \R$. It follows from~\cite[formula~(2.3)]{SoupletW} that  for all $0\le  s \le  \sigma $
\begin{equation}  \label{fPLE2b} 
\begin{split} 
z_b (s) = & z_b (\sigma ) -  \sigma ^\gamma e^{\frac {1} {4\sigma }}  \Bigl( \int _s^\sigma  t ^{-\gamma } e^{-\frac {1} {4t }} \, dt  \Bigr) z_b'(\sigma )  \\ &  - \frac {1} {4} \int _s^\sigma  t ^{-\gamma } e^{-\frac {1} {4t }} \int _t ^\sigma r^{\gamma -2} e^{\frac {1} {4r}}  g (z _{ b } (r) ) \, dr dt
\end{split} 
\end{equation} 
Given $b_0, b \in \R$ and applying~\eqref{fPLE2b}, we have for all $0\le  s \le  \sigma $
\begin{equation*} 
\begin{split} 
 |z _{ b_0 } (s)- z _{ b } (s) |  \le & |z _{ b_0 } (\sigma )- z _{ b } (\sigma ) | \\ &  + \sigma ^\gamma e^{\frac {1} {4\sigma }}  \Bigl( \int _0^\sigma  t ^{-\gamma } e^{-\frac {1} {4t }} \, dt  \Bigr)  |z _{ b_0 } ' (\sigma ) - z _{ b } ' (\sigma ) | \\ &  + \frac {1} {4} \int _0^\sigma  t ^{-\gamma } e^{-\frac {1} {4t }} \int _t ^\sigma r^{\gamma -2} e^{\frac {1} {4r}}  | g (z _{ b_0 } (r) ) - g (z _{ b } (r) ) | \, dr dt
\end{split} 
\end{equation*} 
Setting
\begin{equation*} 
C_1 = \sigma ^\gamma e^{\frac {1} {4\sigma }}   \int _0^\sigma  t ^{-\gamma } e^{-\frac {1} {4t }} \, dt <\infty ,
\end{equation*} 
we deduce that, given any $0< \nu  <\sigma $
\begin{equation*} 
\begin{split} 
\sup _{ 0\le s\le \sigma  } & |z _{ b_0 } (s)- z _{ b } (s) |  \le  |z _{ b_0 } (\sigma )- z _{ b } (\sigma ) |   + C_1 |z _{ b_0 } ' (\sigma ) - z _{ b } ' (\sigma ) | \\ &  + \Bigl( \frac {1} {4} \int _0^\nu   t ^{-\gamma } e^{-\frac   {1} {4t }} \int _t ^\sigma r^{\gamma -2} e^{\frac {1} {4r}} \, dr dt  \Bigr) \sup _{ 0\le r\le \nu  } | g (z _{ b_0 } (r) ) - g (z _{ b } (r) ) | \\ &  +  \Bigl( \frac {1} {4} \int _\nu ^\sigma  t ^{-\gamma } e^{-\frac {1} {4t }} \int _t ^\sigma r^{\gamma -2} e^{\frac {1} {4r}} \, dr dt  \Bigr)  \sup _{ \nu \le r\le \sigma  }  | g (z _{ b_0 } (r) ) - g (z _{ b } (r) ) | .
\end{split} 
\end{equation*} 
We observe that (see~\cite[Lemma~2.1]{SoupletW})
\begin{equation*} 
C_2= \sup  _{ 0<s<\sigma  } \frac {1} {s}\int _0^s t ^{-\gamma } e^{-\frac {1} {4t }} \int _t ^\sigma r^{\gamma -2} e^{\frac {1} {4r}}   \, dr dt <\infty 
\end{equation*} 
and
\begin{equation*} 
C_3 = \int _0 ^\sigma  t ^{-\gamma } e^{-\frac {1} {4t }} \int _t ^\sigma r^{\gamma -2} e^{\frac {1} {4r}} \, dr dt  <\infty 
\end{equation*} 
so that
\begin{equation} \label{fPLE3} 
\begin{split} 
\sup _{ 0\le s\le \sigma  } & |z _{ b_0 } (s)- z _{ b } (s) |   \le  |z _{ b_0 } (\sigma )- z _{ b } (\sigma ) | + C_1  |z _{ b_0 } ' (\sigma ) - z _{ b } ' (\sigma ) | \\ &  + C_2 \nu  \sup _{ 0< r < \nu   }  | g (z _{ b_0 } (r) ) - g (z _{ b } (r) ) |  + C_3 \sup _{ \nu < r < \sigma    }  | g (z _{ b_0 } (r) ) - g (z _{ b } (r) ) |  .
\end{split} 
\end{equation} 
We now fix $b_0\in \R$ and $\delta >0$. 
We deduce from~\eqref{fPLE1b} that if $\eta >0$ is sufficiently small, then
\begin{equation} \label{fPLE6} 
\sup  _{  |b- b_0| \le \eta }  \{   |z _{ b_0 } (\sigma )- z _{ b } (\sigma ) | + C_1  |z _{ b_0 } ' (\sigma ) - z _{ b } ' (\sigma ) |\} \le \frac {\delta } {3} .
\end{equation} 
Moreover, it follows from~\eqref{fPLE2} that there exists $0 < C_4< \infty $ such that
\begin{equation} \label{fPLE4} 
\sup  _{  |b- b_0| \le \eta } \sup _{ 0< r < \sigma  }  | g (z _{ b_0 } (r) ) - g (z _{ b } (r) ) |  \le C_4,
\end{equation} 
and we fix $0< \nu  < \sigma $ sufficiently small so that 
\begin{equation} \label{fPLE5} 
C_2 C_4 \nu  \le \frac {\delta } {3}. 
\end{equation} 
For this fixed value of $\nu$, it follows from~\eqref{fPLE1b} that, assuming $\eta >0$ possibly smaller
\begin{equation} \label{fPLE7} 
 C_3 \sup _{ \nu < r < \sigma    }  | g (z _{ b_0 } (r) ) - g (z _{ b } (r) ) |   \le \frac {\delta } {3}. 
\end{equation} 
Estimates~\eqref{fPLE3}, \eqref{fPLE6}, \eqref{fPLE4}, \eqref{fPLE5} and~\eqref{fPLE7} yield
\begin{equation*} 
\sup  _{  |b- b_0| \le \eta }   \sup _{ 0\le s\le \sigma  }  |z _{ b_0 } (s)- z _{ b } (s) | \le \delta .
\end{equation*} 
Therefore,
\begin{equation} \label{fPLE1b2} 
 \text{the map }b\mapsto z_b  \text{ is continuous } \R\to C( [ 0 ,\sigma  ]).
\end{equation} 
In particular, $L_1 (b)= z_b (0)$ depends continuously on $b$. 

We next show that if $b\not = 0$, then $g_b$ has finitely many zeroes on $(0,\infty )$. It is clear that $g_b$ has finitely many zeros on $[0,\tau ]$, where $\tau $ is defined by~\eqref{fPLE0}. Applying~\eqref{fPLE1}, it remains to show that $z_b$  has finitely many zeros on $[0,\sigma ]$. This is clear if $z_b (0) \not = 0$. Thus we now assume $z_b (0)=0$, and it follows from~\cite[Proposition~2.7~(i)]{SoupletW} that there exists $0<\varepsilon  \le \sigma $ such that $z_b '$ has at most one zero on $(0, \varepsilon )$, so that $z_b$ has at most two zeros on $[0,\varepsilon ]$. Since equation~\eqref{IPE} is not singular on $[\varepsilon ,\sigma ]$, $z_b$ has finitely many zeros on $[\varepsilon ,\sigma ]$, which proves the desired property.

Finally, if $L (  \overline{b} ) \not = 0$, i.e. $z_ {\overline{b} } (0) \not = 0$, then if follows from~\eqref{fPLE1b2} that if $ |b - b_0|$ is sufficiently small, then there exists $\varepsilon >0$ such that $ | z_b (s) |\ge \varepsilon $ for all $0\le s\le \sigma $. By~\eqref{fPLE1}, this means that $  |g_b (r) |\ge \varepsilon r^{-\frac {2} {\alpha }} $ for $r\ge \tau $, with $\tau $ defined by~\eqref{fPLE0}.
Since $g_b$ has the same number of zeros as $g _{ b_0 }$ on $[0, \tau )$ for $ |b-b_0|$ sufficiently small, we deduce that $N_1 (b)= N_1 (b_0)$ provided $ |b-b_0|$ sufficiently small.
\end{proof} 

\begin{proof} [Proof of Theorem~$\ref{eMain3}$~\eqref{eMain3:2}]
We first claim that the sequence $(\MUU _m) _{ m\ge  \overline{m}  }$ of Lemma~\ref{eLFN1} gives rise to a sequence $(\beta  _m) _{ m\ge  \overline{m}  } \subset (0,\infty ) $ satisfying (with the notation~\eqref{fTPA161}-\eqref{fTPA162}) 
\begin{gather} 
\beta _m \goto _{ m\to \infty  }\infty   \label{eNZ3n1:1} \\
 (-1)^m  L_1 (\beta _m ) \goto _{ m\to \infty  }\infty   \label{eNZ3n1:1b1} \\
N_1 (\beta _m ) = m . \label{eNZ3n1:3} 
\end{gather} 
Indeed, we have $w _{ \MUU _m } (0) =   \MUU_m >0$. 
Moreover, $\Lt ( \MUU _m )=0$, so that by Remark~\ref{eTPA3}~\eqref{eTPA3:2}  the limit~\eqref{fL1N1} exists and $\Ltu ( \MUU _m )\not = 0$. On the other hand, $ \wN  (\MUU _m) =m$, so that $(-1)^m w _{ \MUU _m } (s) >0$ for $s$ large.  Thus we see that  $ (-1)^m \Ltu ( \MUU _m ) > 0$. 
We set 
\begin{equation*} 
\beta   _m=  (-1)^m  \Ltu ( \MUU _{ m } ), 
\end{equation*} 
so that 
\begin{equation*} 
g _{ \beta _m }(r)=  (-1)^m  r^{-\frac {2} {\alpha }} w _{ \MUU _m } ( r^{-2} )
\end{equation*} 
for all $r>0$.  Moreover $N_1 (\beta _m) =  \wN ( \MUU  _m ) =m $, $L_1 (\beta _m) = (-1)^m w _{ \MUU _m  } (0)= (-1)^m  \MUU_m $.
In particular, $(-1)^m  L _1 (\beta _m) \to \infty $ as $m\to \infty $, and since $L_1$ is continuous $[0, \infty ) \to \R$, we see that $\beta _m \to \infty $ as $m\to \infty $. This establishes the claim.

Let $\MUU >0$. By~\eqref{eNZ3n1:1b1}, we may choose an integer $m_ \MUU  $ sufficiently large so that
\begin{equation} \label{fMT2n1:1} 
L_1 (\beta _{2m}) >  \MUU  \text{ and }  L_1 (\beta _{2m-1 }) < -  \MUU    \text{ for all } m\ge m_ \MUU  .
\end{equation} 
It follows from Proposition~\ref{eTPA4} (continuity of $L_1$) that there exist 
\begin{equation}  \label{fMT2n1:2} 
0 < \beta _{ 2m-1 }< c_m^- < \beta _{2m} < c_m ^+ < \beta  _{ 2m+1 }
\end{equation} 
such that 
\begin{equation}  \label{fMT2n1:3} 
 L_1 ( c_m^- )= L_1 (c_m ^+ ) =0  \text{ and } L_1( c ) >0  \text{ for all }   c_m^- < c < c_m^+ .
\end{equation} 
Next, since  $N_1 (\beta _{2m}) =2m  $ (by~\eqref{eNZ3n1:3}), we deduce from~\eqref{fMT2n1:2}, \eqref{fMT2n1:3}   and Proposition~\ref{eTPA4} that
\begin{equation}  \label{fMT2n1:4} 
N_1 ( c) = 2m  \text{ for all } c_m^- < c < c_m^+ .
\end{equation} 
From~\eqref{fMT2n1:1} and~\eqref{fMT2n1:3}, it follows that there exist 
\begin{equation} \label{fMT2n1:6} 
c_m^- < c _{  \MUU  , m }^- < \beta _{2m} < c _{  \MUU  ,m }^+ < c_m^+
\end{equation} 
such that $L_1 (c _{  \MUU  , m }^\pm ) =  \MUU  $.
We deduce from~\eqref{fMT2n1:4} and~\eqref{fMT2n1:6} that $N_1 ( c _{  \MUU  ,m }^\pm ) = 2 m$, and from~\eqref{fMT2n1:2}, \eqref{eNZ3n1:1}  and~\eqref{fMT2n1:6} that $c _{  \MUU  ,m }^\pm \to \infty $ as $m\to \infty $.  Thus we see that the sequences $(c _{  \MUU  ,m }^\pm ) _{ m \ge m_ \MUU   }$ satisfy 
\begin{enumerate}[$\bullet$] 

\item 
$c_{ \MUU , m }^- <c _{ \MUU , m }^+$ and $ c_{ \MUU , m }^\pm \to \infty $  as $m\to \infty $; 

\item 
$L_1 (c_{ \MUU , m }^\pm )  =\MUU $ for all $m\ge  m _ \MUU  $;

\item 
 $ N_1 (c_{ \MUU , m }^\pm ) =2m$ for all $m\ge  m _ \MUU  $. 
\end{enumerate} 

By considering $\beta  _{ 2m+1 } \to \infty $ (instead of  $\beta _{2m}$), one constructs as above a sequence $( d _{  \MUU  , m }^\pm ) _{ m \ge  m_ \MUU   } $ such that 
\begin{enumerate}[$\bullet$] 

\item 
$d _{ \MUU , m }^+ < d_{ \MUU , m }^-$ and $ d_{ \MUU , m }^\pm \to -\infty $  as $m\to \infty $;

\item  
$L_1 (d_{ \MUU , m }^\pm )  =\MUU $ for all $m\ge  m _\MUU  $;

\item 
$ N_1 (d_{ \MUU , m }^\pm ) =2m +1$ for all $m\ge  m _\MUU  $. 
\end{enumerate} 

Finally, we may assume $\MUU >0$ without loss of generality, and we define $ \widetilde{g} _{b ^\pm }$ by
\begin{equation*} 
 \widetilde{g}_{b ^\pm } (x)= 
 \begin{cases} 
 g_{b ^\pm } (x) & x\ge 0 \\ - g_{b ^\pm } (-x) & x <0
 \end{cases}  
\end{equation*} 
with ${b ^\pm }= c^\pm  _{ \MUU, \frac {m} {2} }$ if $m\ge m_\MUU $ is even and ${b ^\pm }=  d^\pm  _{ \MUU, \frac {m-1} {2} }$ if $m\ge m_\MUU$ is odd. We observe that $ \widetilde{g}_{b ^\pm } $ are two different solutions of the profile equation~\eqref{fpr2}, which are odd and have $m$ zeros on $(0,\infty )$, hence $2m+1$ zeros on $\R$. Moreover, $ |x|^{\frac {2} {\alpha }}  \widetilde{g}_{b ^\pm } (x)  \to  \MUU$ as $x\to \infty $ and $ |x|^{\frac {2} {\alpha }}  \widetilde{g}_{b ^\pm } (x)  \to  -\MUU$ as $x\to -\infty $. Since $(1 +  |x|^2)^{\frac {1} {\alpha }} (  |\widetilde{g}_{b ^\pm } (x)| + | x \widetilde{g}_{b ^\pm } '(x)|)$ is bounded on $\R$ by~\eqref{fTPA1b61}, Theorem~\ref{eMain3}~\eqref{eMain3:2} now follows from Lemma~\ref{eQSol}, where $f=  \widetilde{ g } _{b^\pm }$ and $\omega (x)= 1$ for $x>0$ and $\omega (x)=-1$ for $x<0$.
\end{proof} 

\appendix

\section{Nonexistence of local, nonnegative solutions} \label{sNEX} 

As explained in the introduction, the main achievement of this paper is the construction of solutions of equation~\eqref{NLHE} with positive initial values, i.e. $\MUU  |x|^{-\frac {2} {\alpha }}$, for which no local in time nonnegative solution exists.  This last assertion is a consequence of~\cite[Theorem~1]{Weissler4}. However, the result in~\cite{Weissler4} only concerns nonnegative solutions of the integral equation~\eqref{fNZ2}. On the other hand, if $0< \alpha \le \frac {2} {N}$, the solutions constructed in Theorems~\ref{eMain1} and~\ref{eMain3} do not  satisfy the integral equation. For completeness, we state and prove below a proposition which establishes nonexistence of nonnegative solutions for both classical solutions of~\eqref{NLHE}, and solutions of the integral equation~\eqref{fNZ2}. 

\begin{proposition} \label{eNEX1} 
Let $\MUU >0 $, $\alpha >0$ and let $\DI (x)= \MUU  |x|^{-\frac {2} {\alpha }}$. 
\begin{enumerate}[{\rm (i)}] 
\item \label{eNEX1:1} 
Suppose either $0<\alpha \le \frac {2} {N}$, or else $\alpha >\frac {2} {N}$ and $\MUU >\MUU_0$ where
\begin{equation} \label{fNEX1} 
\frac {1} {\MUU _0}=  \alpha ^{\frac {1} {\alpha  } } 2^{- \frac {2} {\alpha  } } \pi ^{ - \frac {N} {2}} \int  _{ \R^N  } e^{- |y|^2}  | y |^{-\frac {2} {\alpha }} dy.
\end{equation} 
It follows that for all $T>0$, there is no solution $u \in L^\infty _\Loc ((0,T), L^\infty  (\R^N ) )$ of~\eqref{NLHE}, $u\ge 0$, which is a classical solution on $(0,T)$, and satisfies the initial condition in the sense that $u(t) \to \DI$ in $L^1_\Loc  (\R^N \setminus \{0\}) $.

\item \label{eNEX1:2} 
If $\alpha > \frac {2} {N}$ and $\MUU >\MUU_0$ with $\MUU_0$ defined by~\eqref{fNEX1}, then for all $T>0$, there is no measurable, almost everywhere finite, function $u:(0,T)\times \R^N \to \R$, $u\ge 0$, which satisfies the integral equation~\eqref{fNZ2}. (Note that all terms in~\eqref{fNZ2} are integrals of nonnegative, measurable functions, which are well defined, possibly infinite.)
\end{enumerate} 
\end{proposition} 

\begin{proof} 
Suppose $\alpha >\frac {2} {N}$, so that $\DI\in  L^1 (\R^N ) +L^\infty  (\R^N ) $. 
It follows, since $\DI$ is radially symmetric, radially decreasing, and homogeneous, that
\begin{equation} \label{fNEX2} 
 (\alpha t)^{\frac {1} {\alpha  } }  \| e^{t \Delta } \DI \| _{ L^\infty  } =  (\alpha t)^{\frac {1} {\alpha  } }  e^{t \Delta } \DI (0) =  \alpha ^{\frac {1} {\alpha  } }  e^{ \Delta } \DI (0)  = \frac {\MUU } {\MUU_0}
\end{equation} 
where $\MUU_0 $ is defined by~\eqref{fNEX1}. 
Therefore, if $\MUU >\MUU_0$, then $  (\alpha t)^{\frac {1} {\alpha  } }  \| e^{t \Delta } \DI \| _{ L^\infty  } >1 $ for all $t>0$, and Property~\eqref{eNEX1:2} follows from~\cite[Theorem~1]{Weissler4}.

We next prove Property~\eqref{eNEX1:1}. Suppose there exist $T>0$ and a solution $u$ on $(0,T)$ in the sense of~\eqref{eNEX1:1}. 
Since $u$ is a classical solution on $(0 ,T)$, it is also a solution of the integral equation starting at $u(\tau )$ for $0<\tau <T$, i.e.
\begin{equation}  \label{feFBW2:3} 
u(t)= e^{(t- \tau ) \Delta  } u( \tau ) + \int _\tau ^t e^{(t- s ) \Delta  }  |u(s)|^\alpha u(s) \,ds 
\end{equation} 
for all $0< \tau <t<T$.
Suppose first $\alpha >\frac {2} {N}$ and $\MUU >\MUU_0$.
We deduce from~\cite[Theorem~1]{Weissler4} that for all $0<\tau <T$,
\begin{equation*} 
\sup _{ 0<t<T-\tau    } ( \alpha t )^{\frac {1} {\alpha }} \| e^{t\Delta  } u(\tau ) \| _{ L^\infty  } \le 1.
\end{equation*} 
Given a compact subset $K$ of $\R^N  \setminus \{0\}$, it follows that 
\begin{equation*} 
\sup _{ 0<t<T-\tau } ( \alpha t )^{\frac {1} {\alpha }} \| e^{t\Delta  } {\mathbbm 1}_K u(\tau ) \| _{ L^\infty  } \le 1.
\end{equation*} 
Since ${\mathbbm 1}_K u(\tau )\to {\mathbbm 1}_K \DI$ as $\tau \downarrow 0$ in $L^1 (\R^N ) $, we obtain
\begin{equation*} 
\sup _{ 0<t<T   } ( \alpha t )^{\frac {1} {\alpha }} \| e^{t\Delta  } {\mathbbm 1}_K \DI \| _{ L^\infty  } \le 1.
\end{equation*} 
Since $K\subset \subset \R^N  \setminus \{0\}$ is arbitrary, we conclude that
\begin{equation*} 
\sup _{ 0<t<T   } ( \alpha t )^{\frac {1} {\alpha }} \| e^{t\Delta  }  \DI \| _{ L^\infty  } \le 1
\end{equation*} 
which is in contradiction with~\eqref{fNEX2}. Suppose now $0<\alpha \le \frac {2} {N}$.  It follows from~\eqref{feFBW2:3} that 
\begin{equation} \label{fSL1} 
u(t)\ge  e^{(t- \tau ) \Delta  } u( \tau )
\end{equation} 
for all $0< \tau <t<T$.
We fix $0<t<T$ and $M>0$. Since $\DI$ is not integrable at $x=0$, there exists $0<\varepsilon <1$ such that if $K= \{ \varepsilon  < |x| < 1\}$, then  $e^{t \Delta  }({\mathbbm 1}_K  \DI ) (0) \ge M$. Since $u(\tau )\ge {\mathbbm 1}_K u(\tau )$, we have $u(t, 0)\ge  e^{(t- \tau ) \Delta  } ({\mathbbm 1}_K u( \tau )) (0)$ by~\eqref{fSL1}, and we deduce by letting $\tau \to 0$ that $u(t, 0) \ge M$. Since  $M>0$ is arbitrary, we deduce that $u(t, 0) =\infty $, which is absurd. 
\end{proof} 

\section{Number of oscillations of the profile} \label{sZEROS} 

In this appendix, we comment on  the smallest value of $m_0$ possible in Theorems~\ref{eMain1} and~\ref{eMain3}, which we still call $m_0$.
Recall that by Proposition~\ref{fFIn1}, $m_0 \to \infty $ as $ | \MUU | \to \infty $. 
In addition, we can make the following observations about $m_0$.

If $\alpha >\frac {2} {N}$ and $(N-2) \alpha <4$, then there exist positive self-similar solutions of~\eqref{NLHE}.  
We conjecture that for every $\MUU \not = 0$ with $ | \MUU|$ sufficiently small, we may take $m_0=0$ in Theorem~\ref{eMain1} and in Theorem~\ref{eMain3}~\eqref{eMain3:1} and~\eqref{eMain3:2}. 

If $\alpha \le \frac {2} {N}$, then~\eqref{NLHE} does not have any global, positive solution. In particular, the profile of any nontrivial self-similar solution must have at least one zero. It follows that for every $\MUU \not =  0$, $m_0 \ge 1$ in Theorem~\ref{eMain1} and in Theorem~\ref{eMain3}~\eqref{eMain3:1}. 
More can be said, and we consider separately the cases $N=1$ and $N\ge 2$.

First consider equation~\eqref{NLHE} in one dimension.
Suppose $\alpha \le 2$ and let the integer $k\ge 0$ be defined by $\frac {2} {k+2}< \alpha \le  \frac {2} {k+1}$.
It follows from~\cite{MizoguchiY} that every nontrivial solution of~\eqref{NLHE} must have at least $k+1$ zeros on $\R$ for every $t>0$. In particular, the profile of every nontrivial self-similar solution must have at least $k+1$ zeros on $\R$. 
Thus we see that $m_0$ in Theorem~\ref{eMain3}~\eqref{eMain3:1} and~\eqref{eMain3:2} satisfies $m_0 \ge \frac {k+1} {2}$ and $m_0 \ge \frac {k-1} {2}$, respectively. Note that these lower estimates are independent of $\MUU $ and go to $\infty $ as $\alpha \to 0$. 

We now consider the case $N\ge 2$ and we claim that for all $0< \alpha \le \frac {2} {N}$ there exists an integer $m(\alpha )\ge 1$ with $m(\alpha ) \to \infty $ as $\alpha \to 0$, such that $m_0$ in Theorem~\ref{eMain1} must satisfy $m_0\ge m(\alpha )$ for all $\MUU \not = 0$. 
Indeed, the self-similar solutions in Theorem~\ref{eMain1} all have profiles which are solutions of~\eqref{fpr3}-\eqref{fpr3:1} for some $a\not = 0$. Therefore, if we set $m(\alpha ) = \inf \{ N(a);\, a\not = 0 \}$, where $N(a) $ is the number of zeros of the solution of~\eqref{fpr3}-\eqref{fpr3:1} (see~\eqref{nz}), then $m_0$ in Theorem~\ref{eMain1} satisfies $m_0\ge m(\alpha )$ for all $\MUU \not = 0$. 
We claim that 
\begin{equation} \label{fZ2} 
m(\alpha ) \goto _{ \alpha \downarrow 0} \infty 
\end{equation} 
To prove~\eqref{fZ2}, we note that by symmetry $m(\alpha ) = \inf \{ N(a);\, a> 0 \}$. Furthermore, we
may assume $(N-2)\alpha <2$, and  it follows from~\cite[Proposition~1]{Yanagida} that $N (a)$ is a nondecreasing function of $a$. 
Thus we see that $m (\alpha )=\inf\{N (a);\, a \in (0,1)\}$. 
Given $a\in (0,1)$ and $f_a$ the solution of~\eqref{fpr3}-\eqref{fpr3:1},  we let $z(s)=\frac {1} {a}f _a (r)$, where $ s=\alpha ^{-\frac {1} {2}}r$, so that $z(0)=1$, $z'(0)=0$ and 
\begin{equation*}
z''(s) + \frac{N-1}{s}z'(s) + z(s) +\alpha  \frac{s}{2}z'(s) +\alpha a^{\alpha }|z(s)|^{\alpha}z(s) = 0.
\end{equation*}
An obvious energy argument shows that $z$ and $z'$ are bounded, uniformly in $s\ge 0$, $\alpha >0$ and in $a\in (0,1)$. 
In particular, $\alpha  \frac{s}{2}z'(s) +\alpha a^{\alpha }|z(s)|^{\alpha}z(s) \to 0$ as $\alpha \to 0$, uniformly for  $s$ bounded  and  $a\in (0,1)$. Since the solution $Z$ of the limiting problem as $\alpha \to 0$, i.e. $Z''+\frac {N-1} {r}Z'+ Z=0 $, $Z(0)=1$ and $Z'(0) =0$,  is well known to have infinitely many zeros on $(0,\infty )$ (see e.g.~\cite{Barrett}), the claim~\eqref{fZ2} follows from a standard perturbation argument.


\begin{thebibliography}{99}

\bibitem{Baras}{P. Baras,} {Non unicit\'e des solutions d'une \'equation d'\'evolution
non lin\'eaire, Ann Fac. Sci. Toulouse Math.  {\bf 5} (1983), no. 3-4, 287--302.}
\MScN{MR0747196} \LINK{http://www.numdam.org/item?id=AFST_1983_5_5_3-4_287_0}

\bibitem{Barrett}{J.H. Barrett} {Behavior of solutions of second order self-adjoint differential equations.  Proc. Amer. Math. Soc.  {\bf 6}  (1955). 247--251.}
\MScN{MR0070804} \DOI{10.1090/S0002-9939-1955-0070804-9}

\bibitem{BrezisC}{H. Brezis and T. Cazenave,} {A nonlinear heat equation with
singular initial data, J. Anal. Math. {\bf 68} (1996), 277-304.}
\MScN{MR1403259} \DOI{10.1007/BF02790212}

\bibitem{CDNW2}{T. Cazenave T., F. Dickstein, I. Naumkin and F. B. Weissler,} {Perturbations of self-similar solutions. Dyn. Partial Differ. Equ.  {\bf 16} (2019), no. 2, 151--183.}
 \MScN{MR3928036} \DOI{10.4310/DPDE.2019.v16.n2.a3}

\bibitem{CW}{T. Cazenave and F. B. Weissler,} {Asymptotically self-similar
global solutions of the nonlinear Schr\"o\-din\-ger and heat  equations,
Math. Z.  {\bf 228} (1998), no. 1, 83--120.}
\MScN{MR1617975} \DOI{10.1007/PL00004606}

\bibitem{DohmenH}{C. Dohmen and M. Hirose,} {Structure of positive radial solutions to the Haraux-Weissler equation, Nonlinear Anal.  {\bf 33}, n.1 (1998), 51--69.}
\MScN{MR1623045} \DOI{10.1016/S0362-546X(97)00542-7}

\bibitem{FujishimaI} {Fujishima Y. and Ioku N.} {Existence and nonexistence of solutions for the heat equation with a superlinear source term. J. Math. Pures Appl. (9)  {\bf 118} (2018), 128-158.}
\MScN{MR3852471} \DOI{10.1016/j.matpur.2018.08.001}

\bibitem{HarauxW}{A. Haraux and F. B. Weissler,} {Non uniqueness for a
semilinear initial value problem, Indiana Univ. Math. J. {\bf 31} (1982), no. 2,
167--189.}
\MScN{MR0648169} \DOI{10.1512/iumj.1982.31.31016}

\bibitem{LaisterRSVL}{R. Laister, J.C. Robinson, M. Sier\. z\c ega and A. Vidal-L\'opez,} {A complete characterisation of local existence for semilinear heat equations in Lebesgue spaces. Ann. Inst. H. Poincar\'e Anal. Non Lin\'eaire  {\bf 33}  (2016), no. 6, 1519--1538.}
\MScN{MR3569241} \DOI{10.1016/j.anihpc.2015.06.005}

\bibitem{McLeodTW}{K. McLeod, W.C. Troy  and F.B. Weissler,} {Radial solutions of $\Delta 
u+f(u)=0$ with prescribed number of zeros, J. Differential Equations  {\bf 83}  (1990), no. 2, 368--378.}
\MScN{MR1033193} \DOI{10.1016/0022-0396(90)90063-U}

\bibitem{MizoguchiY}{N. Mizoguchi and E. Yanagida,} {Critical exponents for the blowup of
solutions with sign changes in a semilinear parabolic equation. II, J. Differential Equations  {\bf
145}  (1998),  no. 2, 295--331.}
\MScN{MR1621030} \DOI{10.1006/jdeq.1997.3387}

\bibitem{PTW}{L.A. Peletier, D. Terman and F.B. Weissler,} {On the equation $\Delta  u+ \frac{1}{ 2} x\cdot\nabla u-u+f(u)=0$.  Arch. Rational Mech. Anal.  {\bf 94}  (1986), no. 1, 83--99.}
\MScN{MR0831771} \DOI{10.1007/BF00278244}

\bibitem{QuittnerS}{P. Quittner and P. Souplet,} {{\it Superlinear parabolic problems. 
Blow-up, global existence and steady states.} Second edition. Birkh\"auser Advanced Texts: Basler Lehrb\"ucher.  Birkh\"auser/Springer, Cham, 2019.}
\MScN{MR3967048} \DOI{10.1007/978-3-030-18222-9}

\bibitem{SoupletW}{P. Souplet and F.B. Weissler,} {Regular self-similar solutions to
the nonlinear heat equation with initial data above the singular steady state,
Ann. Inst. H.~Poin\-ca\-r\'e Anal. Non Li\-n\'e\-ai\-re  {\bf 20} (2003),
213--235.}
\MScN{MR1961515} \DOI{10.1016/S0294-1449(02)00003-3}

\bibitem{Weissler1}{F.B. Weissler,} {Semilinear evolution equations in Banach
spaces, J. Funct. Anal. {\bf 32} (1979), no. 3, 277--296.}
\MScN{MR0538855} \DOI{10.1016/0022-1236(79)90040-5}

\bibitem{Weissler2}{F.B. Weissler,} {Local existence and nonexistence for semilinear
parabolic equations in $L^{p}$, Indiana Univ. Math. J. {\bf 29} (1980), no.~1,
79--102.}
\MScN{MR0554819} \DOI{10.1512/iumj.1980.29.29007}

\bibitem{Weissler4}{F.B. Weissler,} {$L^{p}$ energy and blow-up for a semilinear heat
equation, Proc. Symp. Pure Math. {\bf 45}, part 2 (1986), 545--552.}
\MScN{MR0843641} \DOI{10.1090/pspum/045.2/843641}

\bibitem{Weissler6}{F.B. Weissler,} {Asymptotic analysis of an ordinary differential equation and non-uniqueness for a
semilinear PDE, Arch. Ration. Mech. Anal.  {\bf 91} (1985), no. 3, 231--245.}
\MScN{MR0806003} \DOI{10.1007/BF00250743}

\bibitem{Yanagida}{E. Yanagida,} {Uniqueness of rapidly decaying solutions to the 
Haraux-Weissler equation, J. Differential Equations  {\bf 127} (1996), no. 2, 561--570.}
\MScN{MR1389410} \DOI{10.1006/jdeq.1996.0083}

\end{thebibliography}
\end{document}